\documentclass{amsart}
\usepackage{amsthm}
\usepackage{graphicx}
\usepackage[colorlinks]{hyperref}
\usepackage{subcaption}

\title{The Cahn-Hilliard equation on an evolving surface}
\author{D. O'Connor}
\email{d.g.m.o-connor@warwick.ac.uk}
\author{B. Stinner}
\email{bjorn.stinner@warwick.ac.uk}
\address{
Mathematics Institute and Centre for Scientific Computing\\
University of Warwick\\
UNITED KINGDOM
}

\def\mdt{\partial^{\bullet}_t}
\def\ndt{\partial^{\circ}_t}
\def\sgrad{\nabla_{\Gamma}}
\def\slap{\Delta_{\Gamma}}
\def\sgradt{\nabla_{\Gamma(t)}}
\def\slapt{\Delta_{\Gamma(t)}}
\def\eps{\varepsilon}
\def\bbb{\boldsymbol}
\def\mmm{\mathcal}
\def\mbR{\mathbb{R}}

\begin{document}

\begin{abstract}
We study the asymptotic limit of the Cahn-Hilliard equation on an evolving surface with prescribed velocity. The method of formally matched asymptotic expansions is extended to account for the movement of the domain. We consider various forms for the mobility and potential functions, in particular, with regards to the scaling of the mobility with the interface thickness parameter. Mullins-Sekerka, but also surface diffusion type problems, can be derived featuring additional terms which are due to the domain evolution. The asymptotic behaviour is supported and further explored with some numerical simulations. 
\end{abstract}

\maketitle


\section{Introduction}

The subject of this paper is the Evolving Surface Cahn-Hilliard (ESCH) equation
\begin{align}
& \mdt u + u\sgradt \cdot v= - \sgradt \cdot \bbb{j}, \label{ESCH1}\\
& \bbb{j} = -M(u) \sgradt w, \qquad w = -\eps \Delta_{\Gamma(t)} u + \frac{1}{\eps} f(u). \label{ESCH2} 
\end{align}
Here, $\{ \Gamma(t) \}_t \subset \mbR^n$ is a smoothly evolving surface, $\mdt = \partial_t + v \cdot \nabla$ is the material time derivative associated with the velocity $v(t) : \Gamma(t) \to \mbR^n$ of material points of $\Gamma(t)$, $\sgradt$ and $\Delta_{\Gamma(t)}$ are the surface gradient and Laplace-Beltrami operator, the function $f(u) = F'(u)$ is the derivative of a double-well potential, and $M(u)$ is a mobility function. Note that $w$ is the variation of the Ginzburg-Landau energy 
\begin{equation} \label{GLEF} 
\mmm{E}_{\eps}(u(t)) = \int_{\Gamma(t)} \frac{\eps}{2}|\sgradt u(t)|^2 + \frac{1}{\eps} F(u(t)) dx
\end{equation}
and in the case $v=0$ the system \eqref{ESCH1}, \eqref{ESCH2} is the $M(u)$-weighted $H^{-1}$ gradient flow of \eqref{GLEF}.

With regards to $F(u)$ we consider potentials with two stable non-degenerate minima denoted by $u_a<u_b$ which are twice continuously differentiable on an interval $(\alpha, \beta)$ containing $[u_a,u_b]$. We consider second order phase transitions such that $F(u_a)= F(u_b)$. Specifically, we have a quartic potential and a logarithmic potential in mind defined by 
\begin{align}
F_{log}(u) &= \tfrac{\theta}{2k_1} \left(  \left(\beta - u\right)\log\left(\beta - u\right)  + \left(u - \alpha\right)\log\left(u - \alpha\right) \right) & & \nonumber \\
& \quad - \tfrac{\theta_c}{2k_2} (\beta - u)(u - \alpha) & \quad & \mbox{(logarithmic)}, \label{eq:pot_log} \\
F_{q}(u)   &= \frac{1}{4}\left(u_b - u\right)^2\left(u - u_a\right)^2 & \quad & \mbox{(quartic)}, \label{eq:pot_q}
\end{align}
where $\theta, \theta_c,k_1,k_2 > 0$ are parameters, but we stress that the results are not restricted to these two cases. 

We assume that the mobility $M(u)$ is Lipschitz on $[\alpha,\beta]$ and positive and continuously differentiable on $(\alpha, \beta)$ (the latter for simplicity, a slightly smaller open interval comprising $[u_a,u_b]$ would be sufficient). We have in mind the two specific mobility functions:
\begin{align}
M_{deg}(u) &= | \bar{M} (u-\alpha)(\beta-u) | & \quad & \mbox{(degenerate)}, \nonumber \\
M_{c}(u) &= \bar{M} & \quad & \mbox{(constant)}, \label{eq:mob_const}
\end{align}
where $\bar{M} > 0$ is a constant. Let us introduce the pairings $(F_q,M_c)$ and $(F_{log},M_{deg})$, the former we refer to as the \textit{constant mobility ESCH equation} and the later as the \textit{degenerate ESCH equation}.

On a stationary, flat domain the Cahn-Hillard equation has been introduced to model phase separation in binary alloy systems \cite{CH58,Cahn61}. As a prototype model for segregation of two components in a mixture it has been applied in many areas beyond materials science. We refer to \cite{anc} for a review of the equation and highlight \cite{HKNZ15} as one example of a recent application of relevance. The field $u$ usually stands for the (mass or volume) concentration of one of the components, sometimes also their difference. Cahn and Hilliard motivated the logarithmic double-well potential \eqref{eq:pot_log} in their original works \cite{CH58,Cahn61} by theories of mixing. The parameter $\theta>0$ is the (constant) temperature of the system and $\theta_c>0$ is a critical temperature dependent on the material which determines the onset of phase separation. In the shallow quench limit ($\theta \nearrow \theta_c$), the logarithmic potential can be well approximated by the quartic potentials of the form \eqref{eq:pot_q}. Non-constant mobilities were motivated by Cahn and Hilliard in the original derivation, see also \cite{G96}. But also the case of a constant mobility \eqref{eq:mob_const} has been of interest \cite{EZ86, B93, anc85}. 

Partial differential equations describing phase separation on evolving surfaces or domains occur, for example, in de-alloying of binary alloys \cite{EE08}, in two-phase flow \cite{HH77} (potentially with soluble surfactants \cite{GLS13}), in pattern formation on growing organisms \cite{LB03}, or in phase separation on biomembranes \cite{ES10,ES,ES13}. In contrast to the usual notion of $u$ as a concentration we here take an abstract point of view choosing not to physically interpret the phase field variable.  We only assume that $u$ is a conserved quantity for which \eqref{ESCH1} is a mass balance on the moving surface $\Gamma(t)$. In particular, $\int_{\Gamma(t)}u(t) = \int_{\Gamma(0)} u(0)$ at all times $t$. The essential difference to the standard Cahn-Hilliard equation is the $u \sgrad \cdot v$ term which accounts for local stretching if $\sgrad \cdot v > 0$ (or compressing in case of the opposite sign). 

After the initial stage of separation, solutions $u$ to the Cahn-Hilliard equation exhibit large domains (or \emph{phases}) in which $u$ is almost constant and close to one of the minima $u_a$, $u_b$ of $F$. These phases are separated by moving layers with a thickness that scales with $\eps$. This behaviour of solutions is a general feature of phase field models. We refer to \cite{RS92} for an analysis of a phase field model with regards to the different time scales at which the phase separation and the movement of the interfacial layers take place. In the latter solution regime, by using formally matched asymptotics expansions, limiting free boundary problems (or \emph{sharp interface models}) as $\eps \to 0$ have been derived. For the Cahn-Hilliard equation in a stationary, flat domain, the pairing $(F_q,M_c)$ has been considered by \cite{Pego} whilst \cite{CENC} have studied the pairing $(F_{log},M_{deg})$ including the deep quench limit $\theta \searrow 0$. The method has also been applied to elliptic equations on fixed hypersurfaces in \cite{ES} where also the underlying surface depends on the solution and, thus, on $\eps$. In some cases such expansions have been rigorously shown to converge, for instance, see \cite{AHM08, ABC94}. In \cite{ABC94} it is required that the resultant free boundary problem admits a smooth solution, thus imposing regularity assumptions on the initial condition. In \cite{Stoth96} these regularity assumptions are relaxed but with the restriction to radially symmetric solutions. Regarding other approaches to assess the sharp interface limit, the $H^{-1}$-gradient flow (of the Ginzburg-Landau energy \eqref{GLEF}) structure has been used in the context of $\Gamma$-convergence to show asymptotic convergence to the Mullins-Sekerka problem in \cite{N08} for the pairing $(F_q,M_c)$. However, when working with a deformable domain, without some relation coupling the surface velocity to the solution, the system not necessarily has a gradient flow structure. 

The general aim of this paper is to investigate the impact of the motion of the underlying domain $\Gamma(t)$. Via a formal asymptotic analysis (for instance, see \cite{FP95a}) we investigate the effects of the surface motion on the limiting problem that is obtained as $\eps \to 0$. The methodology has been applied to surface phase field models in the stationary case where also the surface depends on $\eps$ in \cite{ES}. We have further extended the technique so that we can deal with moving surfaces and can apply it to the ESCH equation. As usual, a coordinate change using the signed distance function to the limiting moving phase interface is performed in the narrow interfacial region which blows up its thickness to unit length. But since the underlying space, $\Gamma(t)$, is time dependent, the scaled distance function must take account of transport due to the surface velocity. Technically, the challenge is to expand the material time derivative $\mdt$ in the new coordinates. The analysis is carried out for the case of hypersurfaces in the three-dimensional space ($n=3$) but the ideas should carry through to the case $n>3$. The only difficulty should consist in dealing with several tangential coordinates along the limiting phase interface rather than one.

The scaling of $M(u)$ (or rather $\bar{M}$) with respect to $\eps$ turns out to be crucial when attempting to derive limiting free boundary problems. In the case of a stationary, flat domain ($v=0$) it is equivalent to study the Cahn-Hilliard equation at different time scales as in \cite{Pego}. Specific scalings have been considered in \cite{anc} where $\bar{M} \sim \eps^1$ and \cite{ER13} where $\bar{M} \sim \eps^0$. The former appears as a model for early time phase separation when the interfaces form and the latter as a long time model for interface evolution. The scaling $\bar{M} \sim \eps^{-1}$ appears in \cite{CENC} for the degenerate mobility in the regime of the deep quench limit, $\theta \searrow 0$, of the logarithmic potential. Each of these time scales has been considered in \cite{DD14}. Different scales have also been discussed in \cite{C89} in the context of a more general phase field model. We here consider a fixed time scale given by the evolution of the surface, namely one given by (a) a typical velocity at which the domains evolves and (b) a length scale given by the typical size of the surface. Different scalings of $\bar{M}$ in $\eps$ then relate to the speed at which diffusion effects are taking place in comparison with transport effects. 

Not for all scalings were we able to identify sensible limiting free boundary problems. If the mobility is too small, i.e., $\bar{M}$ is of a high order in $\eps$, then the limiting problems do not see the long time behaviour resulting from the evolution of the phase field variable, whence the dynamics is purely governed by the transport with the given velocity field $v$. If the mobility is too high so that $\bar{M}$ is of a low (negative) order in $\eps$ then the asymptotic limits are forced towards equilibrium states with respect to the phases which are barely affected by the transport. 

In the interesting intermediate case in which $\bar{M} \sim \eps^0$ we obtain the following general limiting free boundary problem:
\begin{align}
\left. \begin{array}{rcl}
u &=& u_i \\
\sgradt \cdot \left( M(u) \sgradt w(t) \right) &=& u \sgradt \cdot v(t) 
\end{array} \right\} & \text{ in } \Gamma^{i}(t), \, i=a,b, \label{eq:FBP_bulk} \\
\left. \begin{array}{rcl}
\left[ w(t) \right]^+_- &=& 0 \\
w(t) &=& S \kappa_\Lambda(t) \\
\frac{1}{u_b-u_a} [M(u) \sgradt w(t)]^+_- \cdot \mu_\Lambda(t) &=& \big{(} v(t) - v_\Lambda(t) \big{)} \cdot \mu_\Lambda(t) 
\end{array} \right\} & \text{ on } \Lambda(t). \label{eq:FBP_int}
\end{align}
Here, $\Lambda(t)$ is the moving boundary separating the bulk phases $\Gamma^{b}(t) = \Gamma^+(t)$ and $\Gamma^{a}(t) = \Gamma^-(t)$, $[ \cdot ]^+_- = (\cdot)^+ - (\cdot)^-$ stand for the jump across $\Lambda(t)$, $S>0$ is a calibration constant depending on the double-well potential $F$, $\kappa_\Lambda(t)$ is the geodesic curvature of $\Lambda(t)$ with respect to $\Gamma(t)$, $v_\Lambda(t)$ is the velocity of $\Lambda(t)$, and $\mu_\Lambda(t)$ is the co-normal of $\Lambda(t)$ with respect to $\Gamma(t)$ which points into $\Gamma^+(t)$. 

Observe that, in general, $u \sgrad \cdot v$ is a non-trivial right hand side in the elliptic equation for the chemical potential, $w$, in \eqref{eq:FBP_bulk}. This fact causes problems when attempting to pass to the deep quench limit $\theta \searrow 0$ for the degenerate ESCH equation. In that limit the degenerate mobility switches off the elliptic equation in the bulk. On a stationary domain a purely geometric equation is obtained in the sharp interface limit, namely surface diffusion \cite{CENC}. But in the present case a non-trivial term persists in the bulk if $\sgrad \cdot v \neq 0$, and there is no mechanism to account for the mass density changes due to this local stretching or compressing. 

We present our findings as follows. In the next section we introduce the notation that we use throughout our analysis, discussing the notions of surface calculus as well as how we handle curves on evolving surfaces. In Section 3 we present our assumptions for performing an asymptotic analysis on an evolving surface. In particular we discuss the necessary expansion of the material time derivative, $\mdt$, in the inner co-ordinate system. In Section 4 we perform the asymptotic analysis of the ESCH equation for the slow mobility, when $\bar{M} \sim \eps^0$, and interpret the results for specific mobility and potential functions. In Section 5 we turn our attention to the fast mobility, when $\bar{M} \sim \eps^{-1}$, and compare the result with the slow mobility. Finally, in Section 6 we present supporting numerical experiments for the theoretical findings and display some interesting behaviour due to a non-trivial velocity.

\section{Notation and some facts on the Evolving Surface Cahn-Hilliard Equation}

From now on we consider the spatial dimension $n=3$ (unless stated otherwise). If we still use $n$ then this indicates that the concepts and facts which are presented hold true in an arbitrary dimension.  

\subsection{Calculus on evolving surfaces}

Regarding calculus and transport identities on moving hypersurfaces we will now collect some essential concepts and basic facts where we refer the reader to \cite{SFEM} and \cite{ESFEM} for more detail. In this work we will focus on the case of a two-dimensional hypersurface evolving in a three dimensional ambient space but we mention that the facts presented in this section can be generalised in a relatively straightforward way to an arbitrary dimension.

The system \eqref{ESCH1}, \eqref{ESCH2} is stated for a smooth, closed, and connected evolving $(n-1)$-dimensional submanifold $\{\Gamma(t)\}_{t \in [0,T]}$ where $\Gamma(t)$ is embedded in $\mathbb{R}^n$ for $t \in [0,T]$. We assume that it is orientable and denote by $\nu(\cdot, t): \Gamma(t) \to \mbR^n$, $t \in [0,T]$, a spatial unit normal vector field. By $\Gamma_0 := \Gamma(0)$ we denote the initial hypersurface. The space-time manifold for the moving surface is denoted by 
\begin{equation}\label{spacetimegraph}
\mmm{G}_T := \bigcup_{t \in[0,T]} \Gamma(t) \times {t}.
\end{equation}

We have time-dependent \emph{material surfaces} in mind, i.e., a material particle $p$ located at $x_p(t) \in \Gamma(t)$ at a time $t \in [0,T)$ has a velocity $\dot{x}_p(t)$ with normal contributions (which determine the evolution of the shape) and tangential contributions (which are related to transport of material along the surface). We assume that there is a smooth velocity field $v(\cdot, t) : \Gamma(t) \to \mbR^n$, $t \in [0,T)$, such that $\dot{x}_p(t) = v(x_p(t),t)$. 

For a function $f:\mmm{G}_T \to \mbR$ the \emph{material time derivative} in a point $(x,t)$ with $x \in \Gamma(t)$, $t \in [0,T)$, is defined by
\begin{equation*}
 \mdt f (x,t) := \frac{d}{dt} f(x_p(t),t) = \frac{\partial \tilde{f}}{\partial t}(x,t) + v(x,t) \cdot \nabla \tilde{f}(x,t)
\end{equation*}
where $x = x_p(t)$ for a material particle $p$ located at $x_p(t)$ at time $t$. Note that for the expressions on the right hand side to be well-defined a smooth extension $\tilde{f}$ of $f$ to a neighbourhood of $\mmm{G}_T$ is required. The tangential or surface gradient is defined as the projection of the standard derivative onto the tangent plane of the surface so that
\[
\sgrad f := \nabla \tilde{f} - (\nabla \tilde{f} \cdot \nu)\nu.
\]
Denoting by $D_i$ the $i$'th component of the surface gradient we can define the Laplace-Beltrami operator as
\[
\slap f := \sgrad \cdot \sgrad f = \sum_{i=1}^{n} D_i D_i f.
\]
For an arbitrary subdomain $V(t) \subset \Gamma(t)$ with a smooth boundary and a function $f \in C^1(\overline{V(t)})$ \emph{integration by parts} reads \cite{SFEM}
\begin{equation}\label{IBPT}
\int_{V} \sgrad f = -\int_{V} f H \nu + \int_{\partial V} f \mu_{ext}.
\end{equation}
Here, $\mu_{ext}$ is the exterior co-normal on the boundary $\partial V(t)$ that is tangent to $\Gamma(t)$, pointing away from $V(t)$ and orthogonal to $\nu$, and $H := -\sgrad \cdot \nu$ is the mean curvature. 

\emph{Reynold's Transport Formula}, also referred to as the \emph{Leibniz Formula}, enables us to compute the time derivative of a time dependent surface integral. For a material test volume $V(t) \subset \Gamma(t)$ it reads \cite{ESFEM}
\begin{equation} \label{Leibniz}
\frac{d}{dt} \int_{V(t)} f = \int_{V(t)} \mdt f + f \sgrad \cdot v.
\end{equation}
If $V(t) \subset \Gamma(t)$ is not a material volume, i.e., the boundary moves with a speed $v_{\partial V}$ which is different from $v$, then \cite{B86}
\begin{equation} \label{LeibnizB}
\frac{d}{dt} \int_{V(t)} f = \int_{V(t)} \big{(} \ndt f - f H v \cdot \nu \big{)} + \int_{\partial V(t)} f v_{\partial V} \cdot \mu_{ext}
\end{equation}
where $\ndt f = \partial_t \tilde{f} (x,t) + v(x,t) \cdot \nu(x,t) \, \nu(x,t) \cdot \nabla \tilde{f}(x,t)$ is the normal time derivative.

\subsection{Curves on evolving surfaces}
\label{sec:curve}

Let $\{ \Lambda(t) \}_{t \in [0,T]}$ denote a smooth, closed, and connected evolving curve on $\{ \Gamma(t) \}_{t \in [0,T]}$. For all $t \in [0,T]$ it splits the surface $\Gamma(t)$ into two domains which we denote $\Gamma^+(t)$ and $\Gamma^-(t)$. Using the notion of the intrinsic distance between points $x,y \in \Gamma(t)$,
\[
 d_I(x,y,t) := \inf \left\{\int_0^1 \|g'\| \, \bigg |  \, g \in C^1([0,1],\Gamma(t)), \, g(0)=x, \, g(1)=y \right\},
\]
we can define the distance to the curve $\Lambda(t)$ for a point $x \in \Gamma(t)$ as 
\begin{equation} \label{eq:def_distLambda}
 d_{\Lambda(t)}(x) := \inf_{y \in \Lambda(t)} d_I (x,y,t)
\end{equation}
and then the \emph{signed distance function} by 
\[
 d(x,t) := 
 \begin{cases}
   d_{\Lambda(t)}(x) & \quad \mbox{if } x \in \Gamma^+(t), \\
  -d_{\Lambda(t)}(x) & \quad \mbox{if } x \in \Gamma^-(t).
 \end{cases}
\]
This also explains the notation of the domains: We have $d(\cdot,t)>0$ in $\Gamma^+(t)$ and $d(\cdot,t)<0$ in $\Gamma^-(t)$. 

By the smoothness assumptions on $\Lambda(t)$ and $\Gamma(t)$ there is a narrow tubular region of thickness $\bar{\eps} > 0$ independent of $t$ such that for all points in this region there is a unique (modulo reparametrisation) geodesic which realises the distance. In the following, the expression \emph{close to $\Lambda(t)$} refers to this tubular region. Define now the unit tangent vector along the geodesic by
\[
 \mu(x,t) := \sgradt d(x,t), \quad x \in \Gamma(t)
\]
which is a smooth function close to $\Lambda(t)$. Its derivative along the geodesic, $(\sgradt \mu) \mu = \sgradt(\sgradt d) \mu$, then is normal to $\Gamma(t)$. We now choose the unique $\tau(x,t)$ such that $(\tau,\mu,\nu)$ is a positively oriented orthonormal basis of $\mbR^3$ on $\Gamma(t)$ close to $\Lambda(t)$. Then 
\begin{equation} \label{eq:propmu}
 \mu \cdot (\sgradt \mu) \mu = 0, \quad \tau \cdot (\sgradt \mu) \mu = 0.
\end{equation}
For the restrictions of $\tau$ and $\mu$ to $\Lambda$ we will write
\[
 \tau_\Lambda(x,t) := \tau(x,t), \quad \mu_\Lambda(x,t) := \mu(x,t), \qquad x \in \Lambda(t).
\]

Let now $\lambda(s,t)$, $s \in R_\Lambda(t) S^1$, for all $t \in [0,T]$ denote a smooth parametrisation of $\Lambda(t)$ by arc-length. Here, $R_\Lambda(t) S^1$ is the circle around the origin of radius $R_\Lambda(t)$ which is such that $2 \pi R_\Lambda(t)$ is the length of $\Lambda(t)$. Assume that the orientation of the parametrisation is such that $\lambda_s(s,t) = \tau(\lambda(s,t),t)$. Let us introduce
\[
 \tau_\lambda(s,t) := \tau_\Lambda(\lambda(s,t),t), \quad \mu_\lambda(s,t) := \mu_\Lambda(\lambda(s,t),t).
\]
The curvature vector of $\Lambda(t)$ is given by $h(s,t) := \partial_s \tau_\lambda(s,t)$ and, as $\tau_\lambda \cdot \partial_s \tau_\lambda = \frac{1}{2} \partial_s | \tau_\lambda |^2 = 0$, can be split up into a portion $h(s,t) \cdot \nu(\lambda(s,t),t)$ normal to $\Gamma(t)$ and a tangential portion
\begin{equation} \label{eq:kappageo}
 \kappa_\Lambda(\lambda(s,t),t) := h(s,t) \cdot \mu_\lambda (s,t) = \partial_s \tau_\lambda(s,t) \cdot \mu_\lambda(s,t) = - \tau_\lambda(s,t) \cdot \partial_s \mu_\lambda(s,t)
\end{equation}
which is known as the \emph{geodesic curvature} of $\Lambda(t)$ with respect to $\Gamma(t)$. One can show that $\kappa_\Lambda(t) : \Lambda(t) \to \mbR$ is independent of the parametrisation. 

We may parametrise $\Gamma(t)$ close to $\Lambda(t)$ as $x_{\Gamma(t)}(s,r,t)$ by extending the parametrisation $\lambda(s,t)$ where $x_{\Gamma(t)}(s,r,t)$ is the solution of
\[
 \tilde{x}(s,0,t) = \lambda(s,t), \qquad \tilde{x}_r(s,r,t) = \mu(\tilde{x}(s,r,t),t), \quad r \in [-\bar{\eps},\bar{\eps}]. 
\]
For fixed $s$ and $t$ the curve $r \mapsto x_{\Gamma(t)}(s,r,t)$ then is a geodesic and
\begin{equation} \label{eq:identify_r_d}
d(x_{\Gamma(t)}(s,r,t)) = r.
\end{equation}

With $v_\Lambda(t): \Lambda(t) \to \mbR^3$ we denote the (intrinsic) normal velocity of $\Lambda(t)$, i.e., it can have a portion in direction $\nu(t)$ and in direction $\mu_\Lambda(t)$ but $v_\Lambda(x,t) \cdot \tau_\Lambda(x,t) = 0$ for all $x \in \Lambda(t)$, $t \in [0,T]$. Note that as $\Lambda(t) \subset \Gamma(t)$ for all $t \in [0,T]$ the velocity of $\Lambda(t)$ in the direction $\nu(t)$ normal to the surface coincides with the one of the surface,
\[
v_{\Lambda}(x,t) \cdot \nu(x,t) = v(x,t) \cdot \nu(x,t) \quad \forall x \in \Lambda(t).
\]
However, the portion of $v_\Lambda(t)$ which is tangential to $\Gamma(t)$ may be different from the tangential portion of $v(t)$. Observe that 
\begin{equation} \label{eq:velLam}
 \lambda_t(s,t) \cdot \mu_\lambda(s,t) = v_\Lambda(\lambda(s,t),t) \cdot \mu_\Lambda((\lambda(s,t),t).
\end{equation}

\subsection{Motivation of and remarks on the ESCH equation}

For completeness and for the convenience of the reader we follow the lines of \cite{ER13} and briefly derive the Cahn-Hilliard equation in the form \eqref{ESCH1}, \eqref{ESCH2}. Let $u(\cdot,t) : \Gamma(t) \to \mbR$, $t \in [0,T]$, be the density of some scalar conserved quantity which means that for any test volume $V(t) \subset \Gamma(t)$ with external co-normal $\mu_{ext}$:
\begin{equation*}
\frac{d}{dt}\int_{V(t)} u = - \int_{\partial V(t)} j \cdot \mu_{ext}
\end{equation*}
with a (spatially) tangential flux $j(\cdot,t) : \Gamma(t) \to \mbR^n$, $t \in [0,T]$. Using \eqref{IBPT} and the transport formula \eqref{Leibniz} yields
\[
\int_{V(t)} \mdt u + u \sgrad \cdot v + \sgrad \cdot j = 0.
\]
As this must hold for any choice of $V(t)$ we obtain \eqref{ESCH1}. One may now postulate that the flux is driven by the gradient of the chemical potential $w$ given as the first variation of the Ginzburg-Landau energy functional \eqref{GLEF} so that
\[
j = -M(u) \sgrad w.
\]

Many results in the literature on the Cahn-Hilliard equation are obtained for a \emph{dimensionless version} where the minima of the double well potential are located at $\pm 1$. Our system can be transformed to such a setting as follows. Setting 
\[
\tilde{u} = \frac{u-u_b}{u_b-u_a}+\frac{u-u_a}{u_b-u_a} \quad \Leftrightarrow \quad u = \tfrac{1}{2} \big{(} (1+\tilde{u}) u_b + (1-\tilde{u}) u_a \big{)}
\]
we define $\tilde{F}(\tilde{u}) := F(u)$ and $\tilde{M}(\tilde{u}) := M(u)$. Then $f(u) = F'(u) = \frac{2}{u_b - u_a} \tilde{F}'(\tilde{u}) = \frac{2}{u_b - u_a} \tilde{f}(\tilde{u})$, 
and a short calculation shows that \eqref{ESCH1}, \eqref{ESCH2} takes the form
\begin{align}
\mdt \tilde{u} + \tilde{u} \sgrad \cdot v + c_1 \sgrad\cdot v =& \, \sgrad \cdot \left(  \tilde{M}(\tilde{u}) \sgrad \frac{\tilde{w}}{c_2^2} \right) \label{ESCH_DL} \\
\tilde{w} &= \, -\eps c_2 \slap \tilde{u} +\frac{\tilde{f}(\tilde{u})}{c_2\eps}
\end{align}
where $c_1 = \frac{u_b+u_a}{2}$, $c_2 = \frac{u_b-u_a}{2}$ and $\tilde{w}$ is the chemical potential corresponding to the first variation of the energy $\tilde{\mmm{E}}_{\eps}(\tilde{u}) := \mmm{E}_{\eps}(u)$. 

We remark that in \cite{ER13} the case $c_1 = 0$, $c_2 = 1$ is considered. If $c_1 \neq 0$ then the essential difference is the source term $c_1 \nabla \cdot v$ in \eqref{ESCH_DL}.

\section{Assumptions for the asymptotic analysis}

The goal is now to identify the sharp interface problem from the diffuse interface problem by matching appropriate asymptotic $\eps$-expansions. The technique has been carefully detailed in \cite{FP95a}. We can also make use of an extension to elliptic problems on stationary surfaces \cite{ES}. A novel extension to the technique concerns the parabolic case and, in particular, consists of accounting for the time dependence of the domain and re-writing the material time derivative $\mdt$ in inner coordinates close to the phase interface. 

\subsection{Solution regime}
\label{sec:solregime}

We consider solution regimes to \eqref{ESCH1}. \eqref{ESCH2} where phases have formed, in each of which $u$ is close to one of the two minima of $F$ and which are separated by layers with a thickness that scales with $\eps$. Let $(u_\eps,w_\eps)_{\eps > 0}$ denote a family of such solutions and assume that it converges to some pairing $(u_0,w_0)$ such that, at each time $t$, the spatial domain $\Gamma(t)$ is split up into domains $\Gamma^a(t) = \{ u_0(t) = u_a \}$ and $\Gamma^b(t) = \{ u_0(t) = u_b \}$ which is separated by a smooth, closed, and connected evolving curve $\Lambda(t)$ to which the level set $\{ u_\eps(t) = (u_b + u_a)/2 \}$ converges. We note that the asymptotic analysis below, in principle, also works for several smoothly evolving curves as long as there is no topological change. 

The aim is now to identify the equations that govern the evolution of $\Lambda(t)$, $u_0(t)$, and $w_0(t)$. In order to use the notation and setting that has been introduced in Section \ref{sec:curve} we identify $\Gamma^b(t)$ with $\Gamma^+(t)$ and $\Gamma^a(t)$ with $\Gamma^-(t)$. The upper index notation of $^+$ and $^-$ then will denote whether $\Lambda(t)$ is approached from $\Gamma^+(t)$ or $\Gamma^-(t)$, and by $[\cdot]_-^+ = (\cdot)^+-(\cdot)^-$ we denote the jump across $\Lambda(t)$.

\subsection{Outer expansions}

We assume that away from the interfacial layer around the curve $\Lambda(t)$ we can expand the phase field variable and the chemical potential in the form
\begin{equation} \label{oe} 
 u(x,t) = \sum_{i} u_i (x,t) \eps^i, \quad w(x,t) = \sum_{i} w_i (x,t) \eps^i
\end{equation}
in each domain $\Gamma^{a,b}(t)$.

\subsection{Inner coordinates}
\label{sec:innercoords}

As the thickness of the interfacial layer scales with $\eps$ it makes sense to blow it up to unit length in order to be able to study the limit of fields and functions as $\eps \to 0$ in a meaningful way. We therefore introduce the scaled (geodesic) distance on $\Gamma(t)$ to the interface $\Lambda(t)$ by 
\begin{equation*} 
z := \frac{r}{\eps}.  
\end{equation*}
In the interfacial layer we work with the new coordinates $(s,z,t)$. But before we state the (inner) expansions of the fields in these coordinates and state the matching conditions with the outer expansions in the adjacent domains we need to discuss how the differential operators transform by the change of coordinates. 

With regards to the spatial differential operators we can proceed as in \cite{ES}. For fixed $t$ consider the inversion of the map $R_\Lambda(t) S^1 \times [-\bar{\eps}, \bar{\eps}] \ni (s,r) \to x_{\Gamma(t)}(s,r,t) \in \Gamma(t)$. This gives rise to writing $s(x,t) \in R_\Lambda(t) S^1$ and $r(x,t) \in [-\bar{\eps}, \bar{\eps}]$ for a point $x \in \Gamma(t)$ with a distance of order $\eps$ to $\Lambda(t)$. We denote the closest point to $x$ on $\Lambda(t)$ with the distance measured along $\Gamma(t)$ by
\begin{equation} \label{def:xLam}
 x_\Lambda(x,t) = \lambda(s(x,t),t).
\end{equation}
The identity \eqref{eq:identify_r_d} implies that $\eps \sgradt z(x,t) = \sgradt r(x,t) = \mu(x,t)$. Taylor expanding in $x_\Lambda$ then yields
\[
 \sgradt z(x,t) = \frac{1}{\eps} \mu_\Lambda(x_\Lambda,t) + \sgradt \mu(x_\Lambda,t) \mu_\Lambda(x_\Lambda,t) z(x,t) + \mmm{O}(\eps). 
\]
Similarly, we can see that 
\[
 \sgradt s(x,t) = \tau_\Lambda(x_\Lambda,t) + \mmm{O}(\eps). 
\]

For a scalar field $\phi: \Gamma(t) \to \mbR$ and a vector field $b: \Gamma(t) \to \mbR^3$ define $\phi(x,t) = \Phi(s,z,t)$ and $b(x,t) = B(s,z,t)$ close to $\Lambda(t)$. Then we obtain for the surface gradient and the surface divergence in the new coordinates
\begin{align}
\sgradt \phi(x,t) =& \, \Phi_s(s,z,t) \sgradt s + \Phi_z (s,z,t) \sgradt z \nonumber \\
=& \, \tfrac{1}{\eps} \Phi_z(s,z,t) \mu_\Lambda(x_\Lambda,t) + \Phi_z(s,z,t) \sgradt \mu(x_\Lambda,t) \, \mu_\Lambda(x_\Lambda,t) z  \nonumber \\
 &\, + \Phi_s(s,z,t) \tau_\Lambda(x_\Lambda,t) + \mmm{O}(\eps), \label{cov1} \displaybreak[0] \\
\sgradt \cdot b(x,t) =& \, B_s(s,z,t) \cdot \sgradt s + B_z(s,z,t) \cdot \sgradt z \nonumber \\
=& \, \tfrac{1}{\eps} B_z(s,z,t) \cdot \mu_\Lambda(x_\Lambda,t) + B_z(s,z,t) \cdot \sgradt \mu(x_\Lambda,t) \mu_\Lambda(x_\Lambda,t) z  \nonumber \\
& \, + B_s(s,z,t) \cdot \tau_\Lambda(x_\Lambda,t) + \mmm{O}(\eps). \label{cov2}
\end{align} 
Using these identities, \eqref{eq:propmu}, and \eqref{eq:kappageo}, a short calculation shows that we can write for the Laplace-Beltrami operator
\begin{align}
 \slapt \phi(x,t) 
 &= \, \big{(} \sgradt z \cdot \partial_z + \sgradt s \cdot \partial_s \big{)} \big{(} \Phi_z \sgradt z + \Phi_s \sgradt s \big{)} \nonumber \\ 
 &= \, \frac{1}{\eps^2} \Phi_{zz}(s,z,t) + \frac{1}{\eps} \tau_\lambda(s,t) \cdot \partial_s \mu_\lambda(s,t) \Phi_z(s,z,t) + \mmm{O}(\eps^0) \nonumber \\
 &= \, \frac{1}{\eps^2} \Phi_{zz}(s,z,t) - \frac{1}{\eps} \kappa_\Lambda(x_\Lambda,t) \Phi_z(s,z,t) + \mmm{O}(\eps^0) \label{laplacecov}.
\end{align}

With regards to the operator $\mdt$ it will turn out that knowledge of the term to lowest order in $\eps$ is sufficient for the asymptotic analysis. As 
\[
 \mdt \phi(x,t) = \Phi_s(s,z,t) \mdt s(x,t) + \Phi_z(s,z,t) \mdt z(x,t)
\]
and $\mdt z = \frac{1}{\eps} \mdt r$ we need to focus on computing the leading order term of $\mdt r$. 

Without loss of generality, let us consider the case $z(x,t)>0$ and consider again a point $x \in \Gamma(t)$ with a distance of order $\eps$ to $\Lambda(t)$. Let $\tilde{t} \mapsto x_p(\tilde{t})$ be the path of a material particle such that $x = x_p(t)$. For all $\tilde{t}$ in a small open interval containing $t$ denote by $g^m(\cdot,\tilde{t})$ a geodesic which realises the distance $d_{\Lambda(\tilde{t})}(x_p(\tilde{t})$ defined in \eqref{eq:def_distLambda} for the point $x_p(\tilde{t})$. After a suitable reparametrisation of the geodesic and accounting for the closeness of $x$ to $\Lambda(t)$ we may write 
\[ 
\eps z(x_p(\tilde{t}),\tilde{t}) = r(x_p(\tilde{t}),\tilde{t}) = \int_{0}^{\eps} \| g^m_\rho(\rho,\tilde{t}) \|_2 \, d\rho
\] 
where the integrand is $\mmm{O}(\eps^0)$. Then 
\begin{multline} \label{eq:md_r1}
\mdt r(x,t) = \frac{d}{d\tilde{t}} r(x_p(\tilde{t}),\tilde{t}) \big{|}_{\tilde{t} = t} = \int_0^{\eps} \frac{g^m_\rho(\rho,\tilde{t})}{\|g^m_\rho(\rho,\tilde{t})\|_2} \cdot g^m_{\tilde{t}\rho}(\rho,\tilde{t}) \, d\rho \big{|}_{\tilde{t} = t} \\
= - \underbrace{\int_0^{\eps} \left(\frac{g^m_\rho}{\|g^m_\rho\|_2}\right)_\rho(\rho,t) \cdot g^m_t(\rho,t) \, d\rho}_{I_1} 
+ \underbrace{\frac{g^m_\rho(\eps,t)}{\|g^m_\rho(\eps,t)\|_2} \cdot g^m_t(\varepsilon,t)}_{I_{2}} 
- \underbrace{\frac{g^m_\rho(0,t)}{\|g^m_\rho(0,t)\|_2} \cdot g^m_t(0,t)}_{I_3}.
\end{multline}

As the integrand of $I_1$ is $\mmm{O}(\eps^0)$ we see that $I_1$ is $\mmm{O}(\eps)$. 

Denote by $\sigma(\tilde{t}) \in R_\Lambda(\tilde{t}) S^1$ the tangential co-ordinate of the point in $\Lambda(\tilde{t})$ where $g^m$ starts so that $g^m(0,\tilde{t}) = \lambda(\sigma(\tilde{t}),\tilde{t})$ and, in particular, $g^m(0,t) = \lambda(\sigma(t),t) = x_\Lambda(x,t)$. Recalling that $\rho \to g^m(\rho,\tilde{t})$ is a geodesic and that $g^m_\rho(\rho,\tilde{t})$ thus points into the direction of $\mu(g^m(\rho,\tilde{t}),\tilde{t})$ we obtain that 
\begin{multline*}
I_3 = - \mu_\Lambda(\lambda(\sigma(t),t),t) \cdot \big{(} \lambda_s(\sigma(t),t) \dot{\sigma}(t) + \lambda_t(\sigma(t),t) \big{)} \\
= - \mu_\lambda(\sigma(t),t) \cdot \lambda_t(\sigma(t),t) = - \mu_\Lambda(x_\Lambda,t) \cdot v_\Lambda (x_\Lambda,t).
\end{multline*}
by orthogonality of $\lambda_s = \tau_\Lambda$ and $\mu_\Lambda$ and using \eqref{eq:velLam}. 

With regards to $I_2$ we use again that the distance of $x$ to $\Lambda(t)$ is $\mmm{O}(\eps)$ which yields that  
\[
\frac{g_\rho^m (\eps,t)}{\|g^m_\rho (\eps,t)\|_2} = \mu(g^m(\eps,t),t) = \mu(x,t) = \mu_\Lambda(x_\Lambda,t) + \mmm{O}(\eps)
\]
and that 
\[
g^m_t(\eps,t) = \dot{x}_p(t) = v(x,t) = v(x_\Lambda,t) + \mmm{O}(\eps). 
\]

Altogether, we obtain from \eqref{eq:md_r1} that
\[
 \mdt r(x,t) = \big{(} v(x_\Lambda,t) - v_\Lambda(x_\Lambda,t) \big{)} \cdot \mu_\Lambda(x_\Lambda,t) + \mmm{O}(\eps)
\]
so that 
\begin{equation} \label{tcov}
\mdt \phi(x,t) = \frac{1}{\eps} \Phi_z(s,z,t) \big{(} v(x_\Lambda,t) - v_\Lambda(x_\Lambda,t) \big{)} \cdot \mu_\Lambda(x_\Lambda,t) + \mmm{O}(\eps^0).
\end{equation}

\subsection{Inner expansions}

In conjunction with the outer region we will employ two $\eps$-expansions in the inner region. However, in contrast with the outer region, we will use the inner variables discussed in the previous section so that the expansions take the form
\begin{equation} \label{ie} 
 u(x,t) = \sum_{i=0}^{\infty} U_i(s,z,t) \eps^i, \quad w(x,t) = \sum_{i=0}^{\infty} W_i(s,z,t) \eps^i. 
\end{equation}

\subsection{Matching conditions}

The above two expansions valid in the inner and outer regions should match in some intermediary region. Given an arbitrary outer field, $\phi$, with expansion functions $\varphi_i$ and $\Phi_i$ there are a set of matching conditions that these functions should satisfy. These conditions are related to the spatial coordinates only and, thus, are independent of the movement of the domain. Therefore, and because a full derivation can be found in the literature (for instance, see \cite{GS06}) we only state them here: In the limit as $z \to \pm\infty$ 
\begin{subequations}
\begin{align}
\label{matching1A} \Phi_0(s,z,t) \sim& \, \varphi_0^{\pm}(x_\Lambda,t), \\
\label{matching1B} \partial_z \Phi_0(s,z,t) \sim& \, 0, \\
\label{matching2A} \Phi_1(s,z,t) \sim& \, \varphi_1^{\pm}(x_\Lambda,t) \pm \sgradt \varphi_0^{\pm}(x_\Lambda,t) \cdot \mu_\Lambda(x_\Lambda,t) z, \\
\label{matching2B} \partial_z \Phi_1(s,z,t) \sim& \, \pm \sgradt \varphi_0^{\pm}(x_\Lambda,t) \cdot \mu_\Lambda(x_\Lambda,t), \\
\label{matching3B} \partial_z \Phi_2(s,z,t) \sim& \, \pm \sgradt \varphi_1^{\pm}(x_\Lambda,t) \cdot \mu_\Lambda(x_\Lambda,t) + \big{(} \mu_\Lambda(x_\Lambda,t) \cdot \sgradt \big{)}^2 \varphi_0^{\pm}(x_\Lambda,t) z
\end{align}
\end{subequations}
where we recall that the superscripts $\pm$ indicate the limit of the field when approaching $\Lambda(t)$ from $\Gamma^+ = \Gamma^b$ or $\Gamma^- = \Gamma^a$, respectively.

\section{Slow Mobility} \label{SlowMobility}

We begin identifying free boundary problems with the case $\bar{M} \sim \eps^0$. As we will briefly discuss below this is the highest scaling of the mobility in $\eps$ (or the slowest mobility) for which a sensible free boundary problem occurs. 

\subsection{Outer solutions}

Inserting the expansions \eqref{oe} into \eqref{ESCH1} and \eqref{ESCH2}, we match orders of $\eps$. To order $\eps^{-1}$ \eqref{ESCH2} yields
\begin{equation} \label{eq:bulk_u0}
f(u_0) = 0,
\end{equation}
which has $u_0 = u_a$ and $u_0 = u_b$ as stable stationary solutions. Motivated by the assumptions on the setting at the beginning of Section \ref{sec:solregime} we can conclude that $u_0 = u_a$ in $\Gamma^a = \Gamma^-$ and $u_0 = u_b$ in $\Gamma^b = \Gamma^+$ which is the first equation of \eqref{eq:FBP_bulk}. To order $\eps^0$ combining \eqref{ESCH1} with the flux term in \eqref{ESCH2} we obtain a bulk problem for the leading order term of the chemical potential:
\begin{equation}\label{standardbulkproblem}
u_0 \sgrad \cdot v = \sgrad \cdot \left( M(u_0) \sgrad w_0  \right).
\end{equation}
This is the PDE in \eqref{eq:FBP_bulk}. It remains to derive the interface conditions \eqref{eq:FBP_int}. For being able to apply the matching conditions we need to ensure that $u_1$ can have trace values in $\Lambda$. So we briefly look at the equation to next order of \eqref{ESCH2} which reads 
\[
w_0 = f'(u_0) u_1,
\]
whence $w_0$ and $u_1$ coincide up to a constant thanks to the assumption that $F$ has non-degenerate minima. It thus remains to investigate the trace values of $w_0$ in $\Lambda$. 

\subsection{Inner solutions}

We now insert the expansions \eqref{ie} into \eqref{ESCH1} and \eqref{ESCH2} and employ the change of variables formula \eqref{cov1}, \eqref{cov2}. To the lowest order, $\eps^{-2}$, \eqref{ESCH1} yields
\begin{equation} \label{eq:inner_w0}
 0 = \partial_z \left( M(U_0) \partial_z W_0 \right).
\end{equation}
Thus there exists a function $\delta(s,t)$ such that $M(U_0)\partial_zW_0 = \delta(s,t)$. Using the matching condition \eqref{matching1B} and that $M > 0$ on $(u_a,u_b)$ we see that $\delta = 0$ and, thus, $\partial_z W_0 = 0$. This implies that $w_0$ is continuous across the interface $\Lambda(t)$ in the limiting problem which is the first condition of \eqref{eq:FBP_int}. 

To the order $\eps^{-1}$ \eqref{ESCH2} yields
\begin{equation} \label{phaseprofile}
0 = -\partial_{zz} U_0 + f(U_0).
\end{equation}
The matching condition \eqref{matching1A} implies that $U_0 \to u_a$ as $z \to -\infty$ and $U_0 \to u_b$ as $z \to \infty$. The solution is the phase field profile. Well-posedness of the boundary value problem is discussed in \cite{F79} and its references. 

At the same order $\eps^{-1}$ \eqref{ESCH1} gives thanks to the new expansion \eqref{tcov}
\begin{equation} \label{eq:forinterfaceevolution}
\partial_z U_0 (v-v_\Lambda) \cdot \mu_\Lambda + U_0 \partial_z v \cdot \mu_\Lambda = \partial_z \left( M(U_0)\partial_z W_1 \right).
\end{equation}
Here, $v$, $v_\Lambda$, $\mu_\Lambda$, and $\partial_z v$ are evaluated at $(x_\Lambda,t)$ with $x_\Lambda$ defined in \eqref{def:xLam}. Using that $\mu_\Lambda$ and $v_\Lambda$ are independent of $z$ the left hand side reads $\partial_z (U_0 (v - v_\Lambda) \cdot \mu_\Lambda)$. We may integrate with respect to $z$ over the interfacial region, i.e., from $-\infty$ to $+\infty$, to obtain the last condition of \eqref{eq:FBP_int},
\begin{equation} \label{standardinterfaceevolution}
(u_b - u_a) \big{(} v - v_\Lambda \big{)} \cdot \mu_\Lambda = [M(u_0)\sgrad w_0]_-^+.
\end{equation}
Note that we have applied the matching conditions \eqref{matching1A} and \eqref{matching2B} to $U_0$ and $\partial_z W_1$, respectively. 

To the order $\eps^0$ \eqref{ESCH2} gives thanks to \eqref{laplacecov}
\begin{equation}\label{cpic}
W_0 =-\partial_{zz} U_1 + \partial_z U_0 \kappa_\Lambda + f'(U_0) U_1
\end{equation}
where $\kappa_\Lambda$ is evaluated at $(x_\Lambda,t)$. We multiply by $\partial_z U_0$ and integrate over the interfacial region. By differentiating \eqref{phaseprofile} with respect to $z$ we see that $\partial_z U_0$ lies in the kernel of the operator $\partial_{zz} - f'(U_0)$. Using this after an integration by parts we obtain the following solvability condition for \eqref{cpic}:
\begin{equation} \label{standardsolvcond}
w_0 = S(U_0) \kappa_\Lambda
\end{equation}
where 
$$S(U_0) = \Big{(} \int_\mbR (\partial_z U_0)^2 \Big{)}/(u_b-u_a)$$ 
is a constant depending on the phase profile of $U_0$ and, thus, on the double-well potential $F$. This is the last condition of \eqref{eq:FBP_int} so that we have derived the complete free boundary problem \eqref{eq:FBP_bulk}, \eqref{eq:FBP_int}.

\subsection{Discussion}
\label{sec:interpret_slow}

Let us discuss the limiting problem \eqref{eq:FBP_bulk}, \eqref{eq:FBP_int} for some specific choices of mobilities and potentials and compare with previous results for a stationary, flat domain in the literature. Also the case of an even slower mobility scaling with $\eps^1$ is briefly discussed.

\begin{itemize}
 \item {\bf Mass conservation:} In the limiting problem \eqref{eq:FBP_bulk}, \eqref{eq:FBP_int} the total mass $\mmm{M}(t) = \int_{\Gamma(t)} u(t)$ is preserved (as it is in the ESCH equation):
 \begin{align}
  & \frac{d}{dt} \Big{(} \int_{\Gamma^+} u_b + \int_{\Gamma^-} u_a \Big{)} \nonumber \\
  & \quad \overset{\eqref{LeibnizB}}{=} -\int_{\Gamma^+} u_b H v \cdot \nu + \int_\Lambda u_b v_\Lambda \cdot (-\mu_\Lambda)- \int_{\Gamma^-} u_a H v \cdot \nu + \int_\Lambda u_a v_\Lambda \cdot \mu_\Lambda \nonumber \\
  & \quad \overset{\eqref{IBPT}}{=} \int_{\Gamma^+} u_b \sgrad \cdot v + \int_\Lambda u_b (v_\Lambda - v) \cdot (-\mu_\Lambda) + \int_{\Gamma^-} u_a \sgrad \cdot v + \int_\Lambda u_a (v_\Lambda - v) \cdot \mu_\Lambda \nonumber \\
  & \quad \overset{\eqref{eq:FBP_bulk}}{=} \int_{\Gamma^+} \sgrad \cdot (M(u_b) \sgrad w) + \int_{\Gamma^-} \sgrad \cdot (M(u_a) \sgrad w) + \int_\Lambda (u_a - u_b) (v_\Lambda - v) \cdot \mu_\Lambda \nonumber \\
  & \quad \overset{\eqref{IBPT},\eqref{eq:FBP_int}}{=} \int_\Lambda [M(u) \sgrad w]_-^+ \cdot (-\mu_\Lambda) + \int_\Lambda [M(u) \sgrad w]_-^+ \cdot \mu_\Lambda \quad = 0. \label{eq:FBP_masscons} 
 \end{align}
 Thus, if $u_b>u_a>0$ there is a bound on the maximal and minimal surface area where the bounds depend on the initial mass. This implies a restriction on the surface velocity $v$. \\
 Observe that such a restriction also applies to the phase field model if the logarithmic potential \eqref{eq:pot_log} is used as then the value of $u$ is bounded from above by $\beta$ and from below by $\alpha$ so that the total mass has to remain between $\int_{\Gamma} \alpha$ and $\int_{\Gamma} \beta$. However, there is no such restriction in the case of a smooth, globally defined potential such as \eqref{eq:pot_q}. \\
 In turn, there is no restriction at all if $u_{a} <0 <u_b$. 

 \item {\bf Constant mobility:} For the case of a constant mobility and a smooth double-well potential such as $F_q$, \cite{Pego} has shown that the sharp interface limit of the Cahn-Hilliard equation is the Mullins-Sekerka problem \cite{MS63}. It corresponds to \eqref{eq:FBP_bulk}, \eqref{eq:FBP_int} with a flat and stationary surface. One difference is that the curvature, $\kappa_\Lambda$, now is the \emph{geodesic curvature} of the interface. Another difference is the addition of the transport term $v \cdot \mu_\Lambda$ in the evolution law for the interface given in \eqref{standardinterfaceevolution}. 
 The most important difference to the Mullins-Sekerka problem is the surface divergence of the surface velocity in \eqref{standardbulkproblem}. In general, the chemical potential is no longer harmonic, and changes over time can occur due to the time dependence of the surface velocity. 

 \item {\bf Non-constant mobility:} With a non-constant but positive (on $(\alpha,\beta)$) mobility we obtain a limiting Mullins-Sekerka type problem where the diffusivities of the chemical potential in the bulk can differ (see \eqref{standardbulkproblem}) which also impacts on the jump term in \eqref{standardinterfaceevolution}. This result is independent of the choice of the double-well potential as long as the smoothness assumptions on $(\alpha, \beta)$ are met and the minima are located at $u_a$ and $u_b$. However, the choice of $F$ influences the leading order profile (solution to \eqref{phaseprofile}) and, thus, the values of $S(U_0)$ in \eqref{standardsolvcond}. But by appropriate choice of coefficients such as $k_1$ and $k_2$ in $F_{log}$ (or a suitable prefactor for $F_q$) one can ensure that $S(U_0) = 1$.

 \item {\bf Slower mobility:} Let us briefly consider the case of an even slower mobility $\bar{M} \sim \eps^1$. Equation \eqref{eq:bulk_u0} still holds true while \eqref{ESCH1} yields to leading order that $u_0 \sgrad \cdot v = 0$. Within the solution regime defined in Section \ref{sec:solregime}, which implies that $u_0$ is constant in the bulk, we thus obtain the solvability condition $\sgrad \cdot v = 0$. This is a strong restriction on the motion of the surface as it corresponds to local incompressibility. In \eqref{eq:forinterfaceevolution} then $W_0$ features instead of $W_1$. With the matching condition \eqref{matching1B} we then see that $v_\Lambda \cdot \mu_\Lambda = v \cdot \mu_\Lambda$. So the interface is simply transported with the surface velocity and any subtle behaviour due to the Cahn-Hilliard dynamics is lost. We remark that this is no contradiction to the results in \cite{Pego} where, for the slow mobility, a Stefan type problem is shown to emerge because that limit is established at the next higher order in $\eps$. 
 
\end{itemize}

\section{Fast mobility}

A fast mobility scaling $\bar{M} \sim \eps^{-1}$ has been used in \cite{CENC} to derive surface diffusion in the deep quench limit $\theta \searrow 0$ of the Cahn-Hilliard equation with $(F_{log},M_{deg})$ on a flat and stationary domain. We will discuss this problem below but first consider the general, non-degenerate case $\theta > 0$ of $(F_q,M_c)$. 

\subsection{Asymptotic analysis}

As previously, we insert the expansions \eqref{oe} and \eqref{ie} into \eqref{ESCH1} and \eqref{ESCH2} and match orders of $\eps$. 

From the outer expansion of \eqref{ESCH2} to order $\eps^{-1}$ we obtain again that $u_0 = u_b$ or $u_0 = u_a$ in $\Gamma^b = \Gamma^+$ and $\Gamma^a = \Gamma^-$, respectively. Combining \eqref{ESCH1} with the flux term in \eqref{ESCH2} we obtain to order $\eps^{-1}$
\begin{equation}\label{ospm1}
0 = \sgrad \cdot \left( M(u_0) \sgrad w_0 \right).
\end{equation}
Multiplying by $w_0$ and integrating over $\Gamma^+(t) \cup \Gamma^-(t)$ we obtain using \eqref{IBPT}
\begin{align}
 0 &= \int_{\Gamma^+(t)} w_0 \sgrad \cdot \big{(} M(u_0) \sgrad w_0 \big{)} + \int_{\Gamma^-(t)} w_0 \sgrad \cdot \big{(} M(u_0) \sgrad w_0 \big{)} \label{eq:energy_w0} \\
 &= -\int_{\Gamma^+(t)} M(u_0) | \sgrad w_0 |^2 - \int_{\Gamma^-(t)} M(u_0) | \sgrad w_0 |^2 - \int_{\Lambda(t)} \big{[} w_0 M(u_0) \sgrad w_0 \big{]}_-^+ \cdot \mu_\Lambda. \nonumber 
\end{align}
To get an idea of the jump term we match with the inner solutions. 

The inner expansion of equation \eqref{ESCH2} yields the equation \eqref{phaseprofile} to order $\eps^{-1}$ that $U_0$ is the phase transition profile again. From \eqref{ESCH1} we obtain to order $\eps^{-3}$ the equation \eqref{eq:inner_w0} for $W_0$ again, and as before using the matching conditions \eqref{matching1B} and \eqref{matching1A} we can conclude that
\begin{equation} \label{eq:res_w0A}
 \partial_z W_0 = 0 \quad \mbox{and} \quad [ w_0 ]_-^+ = 0.
\end{equation}
Using this and the orthogonality of $\mu_\Lambda$ and $\tau_\Lambda$, to order $\eps^{-2}$ the same equation yields
\[
 0 = \partial_z \big{(} M(U_0) \partial_z W_1 \big{)}.
\]

Similarly, we can conclude that $\partial_z W_1 = 0$ and, using the matching conditions \eqref{matching2B} and \eqref{matching2A}, 
\begin{equation} \label{pmbulkproblemBC1}
 0 = [ M(u_0) \sgrad w_0 ]_-^+ \cdot \mu_\Lambda \quad \mbox{and} \quad [ w_1 ]_-^+ = 0.
\end{equation}
Together with \eqref{eq:res_w0A} we see that the last term of \eqref{eq:energy_w0} vanishes, and we can conclude that $\sgrad w_0 = 0$ in $\Gamma^+(t)$ and $\Gamma^-(t)$ so that
\begin{equation} \label{eq:res_w0B}
 w_0(t) \mbox{ is constant on } \Gamma^+(t) \cup \Gamma^-(t).
\end{equation}
We have explicitly noted the time dependence to clarify that $w_0$ can and, in general, will change over time (see below).

Continuing with the outer expansions, \eqref{ESCH2} to order $\eps^0$ yields $w_0 = f'(u_0) u_1$ so that also $u_1$ is constant where we recall that $F''(u_0) = f'(u_0) \neq 0$ for $u_0 \in \{ u_a,u_b \}$ thanks to the assumption that $F$ has non-degenerate minima. Using that $\sgrad w_0 = 0$, equation \eqref{ESCH1} to order $\eps^0$ yields the following elliptic bulk problem for $w_1$:
\begin{equation}\label{pmbulkproblemPDE}
u_0 \sgrad \cdot v = \sgrad \cdot \left( M(u_0)\sgrad w_1  \right).
\end{equation}

One boundary condition is given by \eqref{pmbulkproblemBC1}. In order to determine a second one, consider the inner expansion of \eqref{ESCH1} to order $\eps^{-1}$. Using \eqref{tcov} and that $\partial_z W_0 = 0$, $\partial_s W_0 = 0$ (thanks to \eqref{eq:res_w0B}), and $\partial_z W_1 = 0$ as well as the orthogonality of $\mu_\Lambda$ and $\tau_\Lambda$, a short calculation shows that it greatly simplifies to
\begin{equation} \label{forpmbulkproblemBC2}
\partial_z U_0 (v-v_\Lambda) \cdot \mu_\Lambda + U_0 \partial_z v \cdot \mu_\Lambda = \partial_z \left( M(U_0) \partial_z W_2 \right).
\end{equation}
It reads as \eqref{standardinterfaceevolution} except that $W_1$ is replaced by $W_2$. Integrating with respect to $z$ over $\mbR$, treating the left hand side in the same manner as done for \eqref{standardinterfaceevolution}, and applying \eqref{matching3B} to the right hand side where we use that $\sgrad w_0 = 0$ we arrive at 
\begin{equation} \label{pmbulkproblemBC2}
(u_b - u_a) \big{(} v - v_\Lambda \big{)} \cdot \mu_\Lambda = [M(u_0)\sgrad w_1]_-^+.
\end{equation}
Returning to the higher order inner expansions, from equation \eqref{ESCH2} to order $\eps^{0}$, we obtain \eqref{cpic} again, and conclude as before that \eqref{standardsolvcond} holds true. With \eqref{eq:res_w0B} we obtain that also
\begin{equation} \label{pmsolvcond}
 \kappa_\Lambda(t) = \frac{1}{S(U_0)} w_0(t) \mbox{ is constant along } \Lambda(t) \mbox{ at all times } t. 
\end{equation}
Since $W_0 = S(U_0) \kappa_\Lambda$, writing $U_1 = \tilde{u} \kappa_\Lambda$ and substituting into \eqref{cpic}, then $\tilde{u}$ can be determined as the unique function solving
\begin{equation} \label{eq:tu_U0}
-\partial_{zz}\tilde{u} + f'(U_0) \tilde{u} = S(U_0) + \partial_zU_0
\end{equation}
subject to the boundary condition $\lim_{z \to \pm\infty} \partial_z\tilde{u}=0$ from \eqref{matching1B}.

Finally at order $\eps$ we obtain
\begin{equation}
W_1 = -\partial_{zz} U_2 + \partial_z U_1 \kappa_\Lambda + f'(U_0)U_2 +f''(U_0) \frac{U_1^2}{2}
\end{equation}
This gives us a method to determine the interface condition for the first order term of the chemical potential. Multiplying by $\partial_z U_0$ and integrating as before we can determine $w_1$ to be:
\begin{equation*}
w_1 = \frac{\kappa^2}{u_b-u_a} \int_{-\infty}^{\infty}\partial_z\tilde{u} \partial_zU_0 - \frac{\tilde{u}^2}{2}\partial_zf'(U_0).
\end{equation*}
We may express this in a short from as
\[
w_1 = T(U_0) \kappa_\Lambda^2,
\]
where $T(U_0)$ is a constant depending on the leading order phase profile in the inner region. We have suppressed the dependence on $\tilde{u}$ by noting the dependence of $\tilde{u}$ on the phase profile $U_0$ (see \eqref{eq:tu_U0}).

\subsection{Discussion}

To summarise the findings of the preceding section: The phase interface is in spatial equilibrium in the sense that the geodesic curvature is constant, see \eqref{pmsolvcond}. In the thus split domain we have the set of equations: 
\begin{align}
\left. \begin{array}{rcl}
u &=& u_i \\
\sgrad \cdot \left( M(u) \sgrad \tilde{w}(t) \right) &=& u \nabla_{\Gamma} \cdot v(t) 
\end{array} \right\} & \text{ in } \Gamma^{i}(t), i=a,b, \label{eq:HOP_bulk}\\
\left. \begin{array}{rcl}
\left[ \tilde{w}(t) \right]^+_-= 0 \\
\tilde{w}(t) = T\kappa_\Lambda^2(t)\\
\frac{1}{u_b-u_a} \left[ M(u) \sgrad \tilde{w}(t) \right]^+_- \cdot \mu_\Lambda(t) &=& \big{(} v(t) - v_\Lambda(t) \big{)} \cdot \mu_\Lambda(t)
\end{array}\right\}&  \text{ on } \Lambda(t). \label{eq:HOP_int}
\end{align}

\begin{itemize}
 \item {\bf Mass conservation:} First, observe that the total mass still is preserved in the sharp interface limit potentially implying a restrictions on the velocity $v$. In the identity \eqref{eq:FBP_masscons} $w$ has to be replaced by $\tilde{w}$ for this purpose. 
 
 \item {\bf Interface evolution:} The solvability condition \eqref{pmsolvcond} is an equilibrium condition with respect to the phase separation. This restriction seems reasonable since the fast scaling of the mobility acts to blow up the effects of the Cahn-Hilliard dynamics. But the equilibrium condition \eqref{pmsolvcond} alone doesn't tell us much about the evolution of $\Lambda(t)$. In fact, at a given time $t$ there may be several possible curves $\Lambda(t)$ of constant geodesic curvature such that the mass side condition is satisfied. For instance, if $\Gamma(t)$ is a sphere one will find an infinite number. By the assumptions in Section \ref{sec:solregime} the interface is approximated by level sets of the phase field solutions. Thus, one may expect it to evolve smoothly, and one will also expect that a specific curve is picked in the sharp interface limit. We leave this question open for future studies. 
 
 \item {\bf Deep quench limit of the degenerate equation:} The deep quench limit of \eqref{ESCH1} and \eqref{ESCH2} for the degenerate ESCH equation corresponds to the limit as $\theta \searrow 0$. Then $u_a \to \alpha$ and $u_b \to \beta$ so that the degenerate mobility $M_{deg}(u)$ is switched off in the bulk. In the case of a stationary, flat domain the limiting problem is surface diffusion and has been derived in \cite{CENC}. There, the flux $\bbb{j}$ is expanded in addition to the fields and some matching conditions are replaced by assumptions on the limits of the fluxes when approaching the boundaries of the interfacial layer. This is due to a lack of equations for the bulk fields. \\
 Indeed, also in our case, \eqref{ospm1} does not exist so that we have no equation for $w_0$ in the bulk. In particular, we cannot conclude any more that $\sgrad w_0 = 0$. Similarly, there is no bulk equation for $w_1$: Equation \eqref{pmbulkproblemPDE} reduces to $u_0 \sgrad \cdot v = 0$. Within the solution regime defined in Section \ref{sec:solregime} this means necessarily that $$\sgradt \cdot v(t) = 0$$ in the bulk phases, the implication of which has been discussed in the context of a very slow mobility already (see Section \ref{sec:interpret_slow}). As we also cannot conclude any more that $\partial_s W_0 = 0$ another term of the form $M(u_0) \partial_{ss} W_0$ appears on the right hand side of \eqref{forpmbulkproblemBC2}. Integrating and using suitable assumptions for the flux $M(u_0) \partial_z W_2$ as in \cite{CENC} we obtain 
 \begin{equation} \label{eq:surfdiff}
  (u_b - u_a) \big{(} v(t) - v_\Lambda(t) \big{)} \cdot \mu_\Lambda(t) =\tilde{S}(U_0) \Delta_{\Lambda(t)} \kappa_\Lambda(t)
 \end{equation}
 instead of \eqref{pmbulkproblemBC2}. Here, $\Delta_{\Lambda(t)}$ corresponds to $\partial_{ss}$ after parametrisation and stands for the Laplace-Beltrami operator on the curve $\Lambda(t)$, and $\tilde{S}(U_0) = S(U_0)\int_\mbR M(U_0)$. Equation \eqref{eq:surfdiff} is \emph{surface diffusion} for a curve on a moving surface where the velocity $v$ of the underlying surface manifests by an additional transport term. 

\end{itemize}

\section{Numerical Experiments}\label{sec:Numerics}

Using numerical simulations, the aims of this section are: (1) to support the theoretical findings on the convergence as $\eps \to 0$ stated in the previous sections, and (2) to illustrate and display some of the possible effects due to the motion of the surface. The computational method is based on the evolving surface finite element method \cite{ESFEM} which has been applied to the ESCH equation in \cite{ER13}. It has been implemented in MATLAB for $1$D simulations and in DUNE \cite{dunegridpaperI:08, dunegridpaperII:08, dunefempaper:10} for $2$D simulations where the dimension refers to the manifold. In addition, for the experiment in Section \ref{sec:MovingSphere} we have implemented the Arbitrary Lagrangian Eulerian evolving surface finite element method via the finite element library AMDiS \cite{VV07}. The ALE-ESFEM method was first proposed in \cite{ES-ALE} and analysed further in \cite{EV-ALE}.

With regards to the 1D simulations, we have produced results on a bounded interval with Neumann type boundary conditions, that is $u_x = w_x = 0$, which contradicts the setting of the analysis where we assumed a closed surface (see around equation \eqref{spacetimegraph}). However, we can double the (time dependent) interval and reflect the solution to make it symmetric with respect to the centre. The thus obtained setting can be further extended periodically to the whole real line so that we may think of a solution on an object which topologically is a circle. 

We only carried out computations with the quartic potential \eqref{eq:pot_q} and the constant mobility \eqref{eq:mob_const}. For the Cahn-Hilliard equation on the real line there exists an equilibrium profile given by
\begin{equation} \label{eq:num_eq_profile}
 \frac{u_b + u_a}{2} + \frac{u_b - u_a}{2} \tanh \Big{(} \frac{u_b-u_a}{2\sqrt{2}} \frac{y}{\eps} \Big{)}, \qquad y \in \mathbb{R}.
\end{equation}
We use this profile to specify initial conditions $u_{IC}(x) = u(x,0)$, $x \in \Gamma(0)$, unless stated otherwise. 

\subsection{Stretching and Compression}
\label{sec:deform_int}

\begin{table}
 \tabcolsep=3pt
 \begin{tabular}{llll}
  \hline
  Parameter & Data for Figure \ref{stretchpm1} & Data for Figure \ref{shrinkpm1} & Data for Figure \ref{stretchpos} \\
  \hline
  $u_a$, $u_b$; $\bar{M}; T$ & -1,  1;  1;  10 & -1,  1;  1;  10 & 0.2,  0.8;  1;  2 \\
  $\Gamma(t)$ & $\left \{ \begin{array}{ll} (0,1+t) & t \leq 2 \\ (0,3) & t > 2 \end{array} \right.$ & $\left \{ \begin{array}{ll} (0,3-t) & t \leq 2 \\ (0,1) & t > 2 \end{array} \right.$ & $(0,1+t)$ \\
  $v(x,t)$, $x \in \Gamma(t)$ & $\left \{ \begin{array}{ll} \frac{x}{(t+1)} & t \leq 2 \\ 0 & t > 2 \end{array} \right.$ & $\left \{ \begin{array}{ll} -\frac{x}{(t+1)} & t \leq 2 \\ 0 & t > 2 \end{array} \right.$ & $\frac{x}{(t+1)}$ \\
  $u_{IC}(x)$ & $0.9 \tanh(10x-5)$ & $0.9 \tanh(10x-15)$ & $0.3 \tanh(\frac{x-0.5}{\eps})+0.5$ \\
  \hline
 \end{tabular}
 \caption{Simulation data for Section \ref{sec:deform_int}.} \label{tbl:deform_int} 
\end{table}

\begin{figure}
\begin{subfigure}{0.24\textwidth}
\includegraphics[clip, trim=0.5cm 6cm 0.5cm 6cm,width = \textwidth]{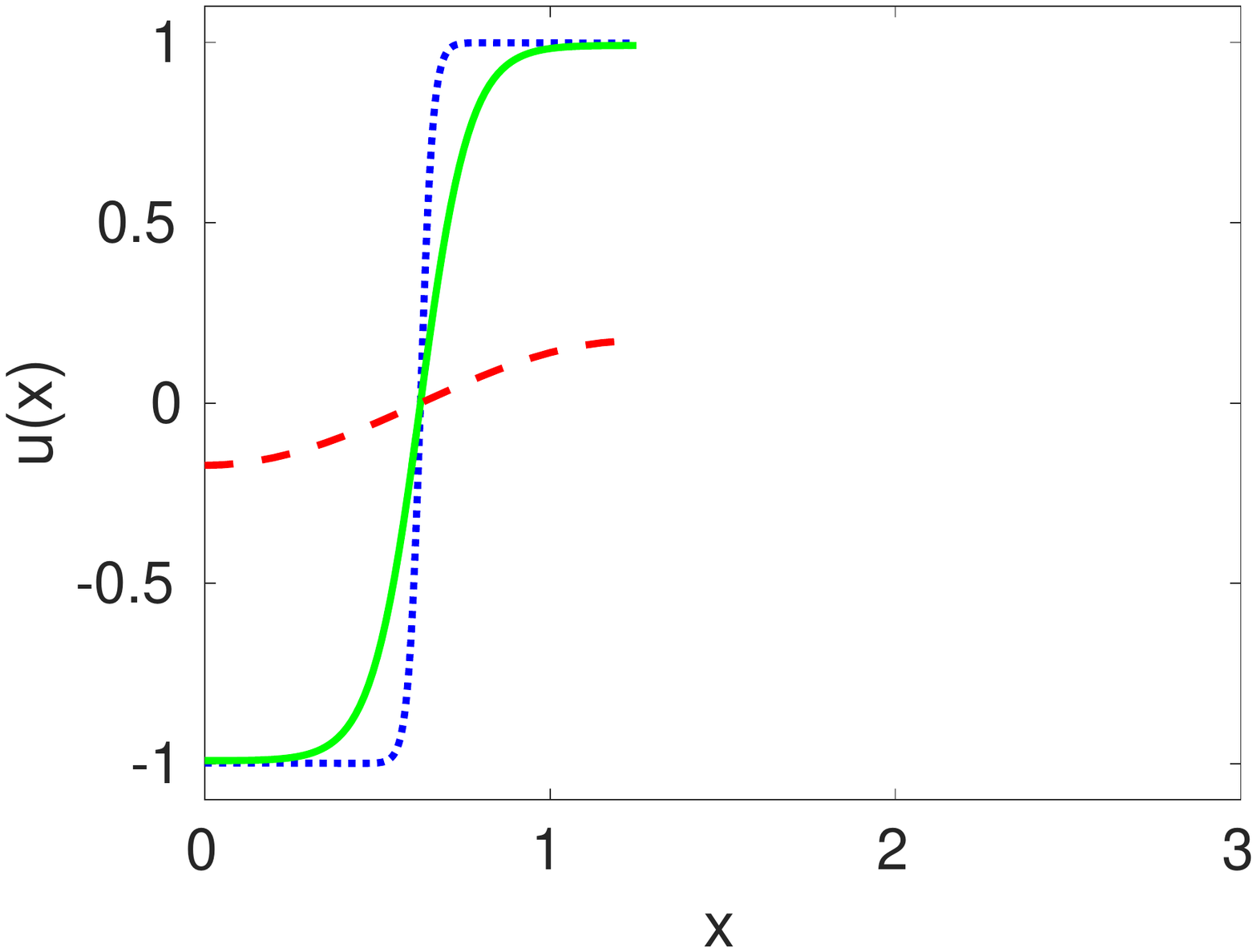}
\caption{$t=0.25$ }
\end{subfigure}
\begin{subfigure}{0.24\textwidth}
\includegraphics[clip, trim=0.5cm 6cm 0.5cm 6cm,width = \textwidth]{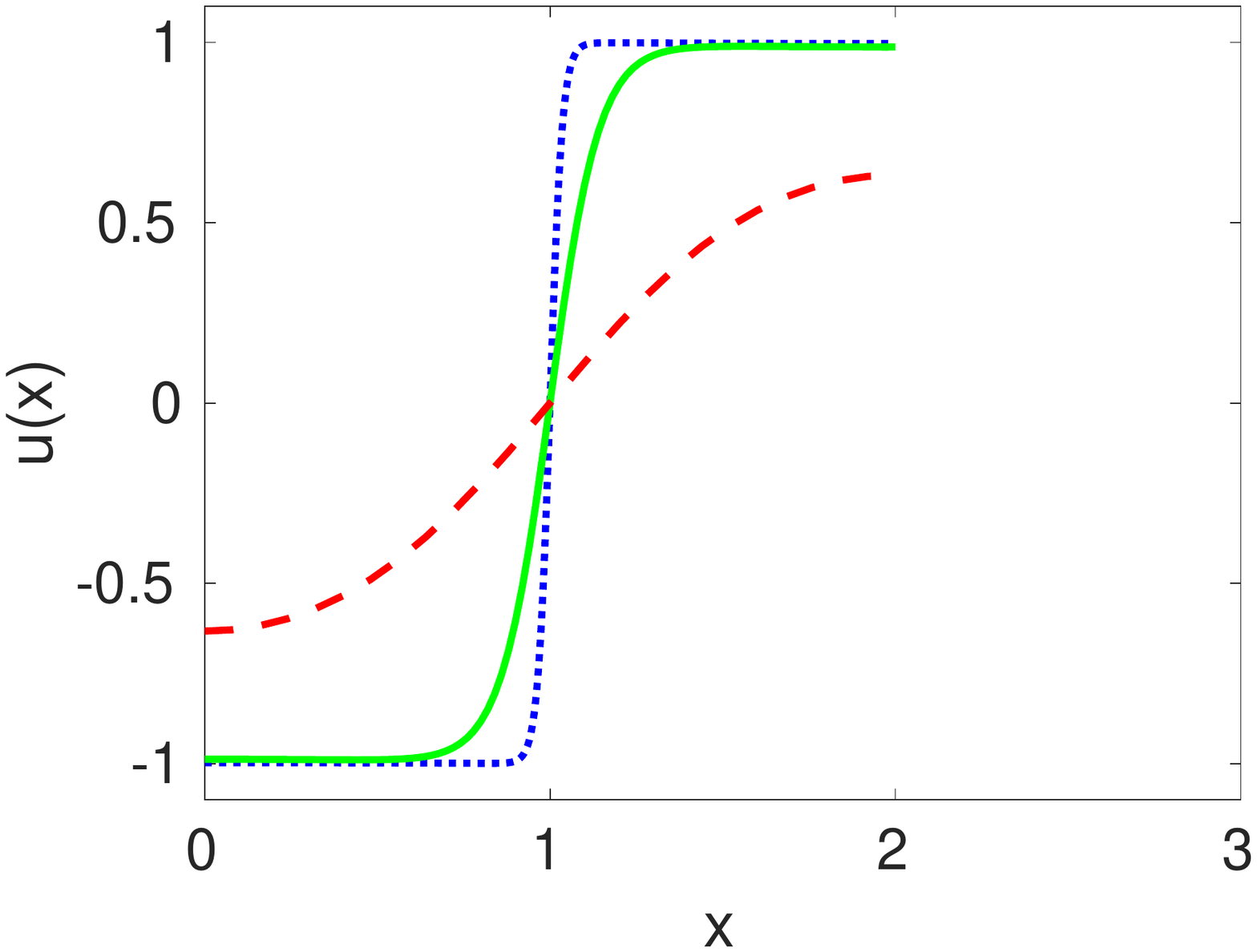}
\caption{$t=1.0$ }
\end{subfigure}
\begin{subfigure}{0.24\textwidth}
\includegraphics[clip, trim=0.5cm 6cm 0.5cm 6cm,width = \textwidth]{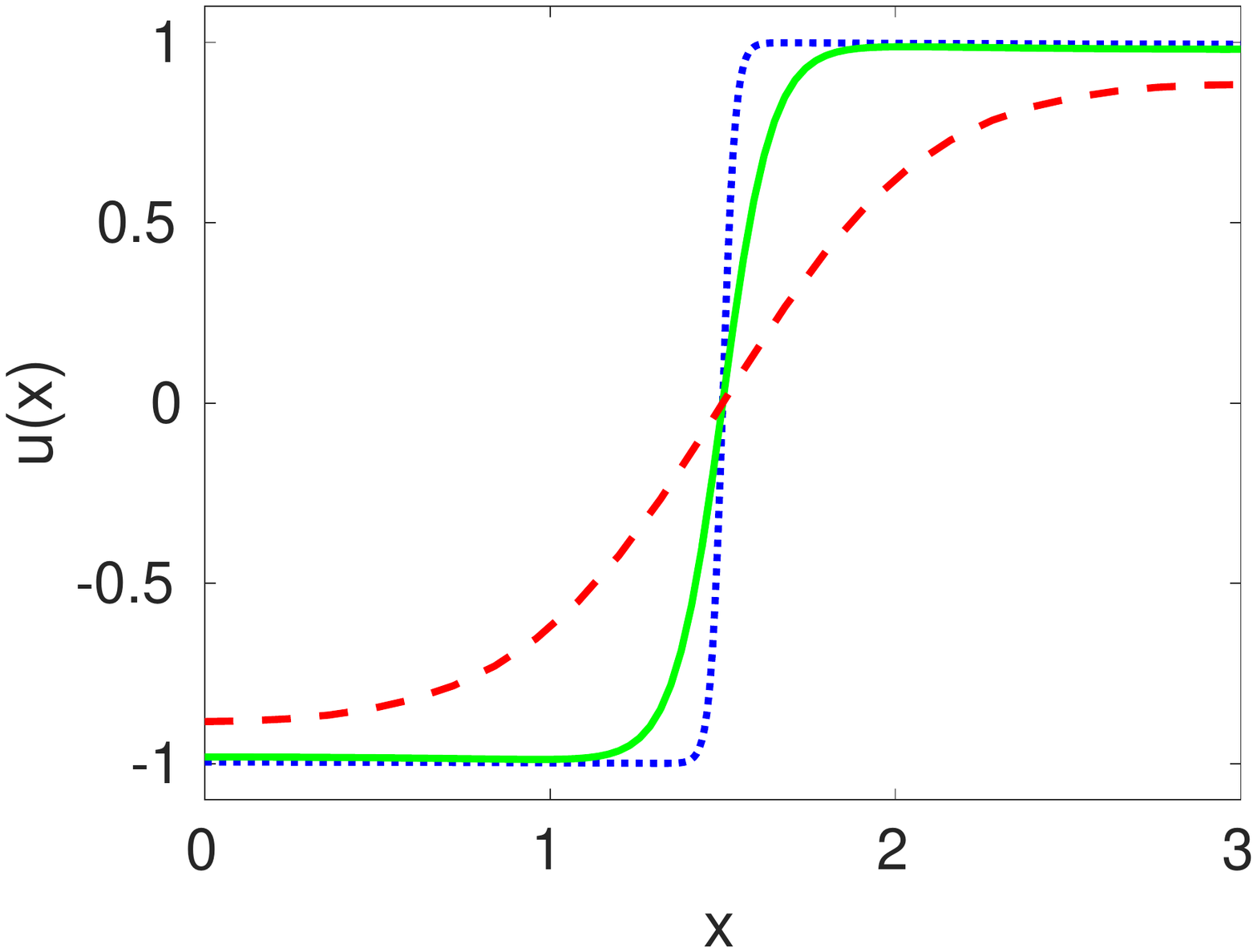}
\caption{$t=2.0$ }
\end{subfigure}
\begin{subfigure}{0.24\textwidth}
\includegraphics[clip, trim=0.5cm 6cm 0.5cm 6cm,width = \textwidth]{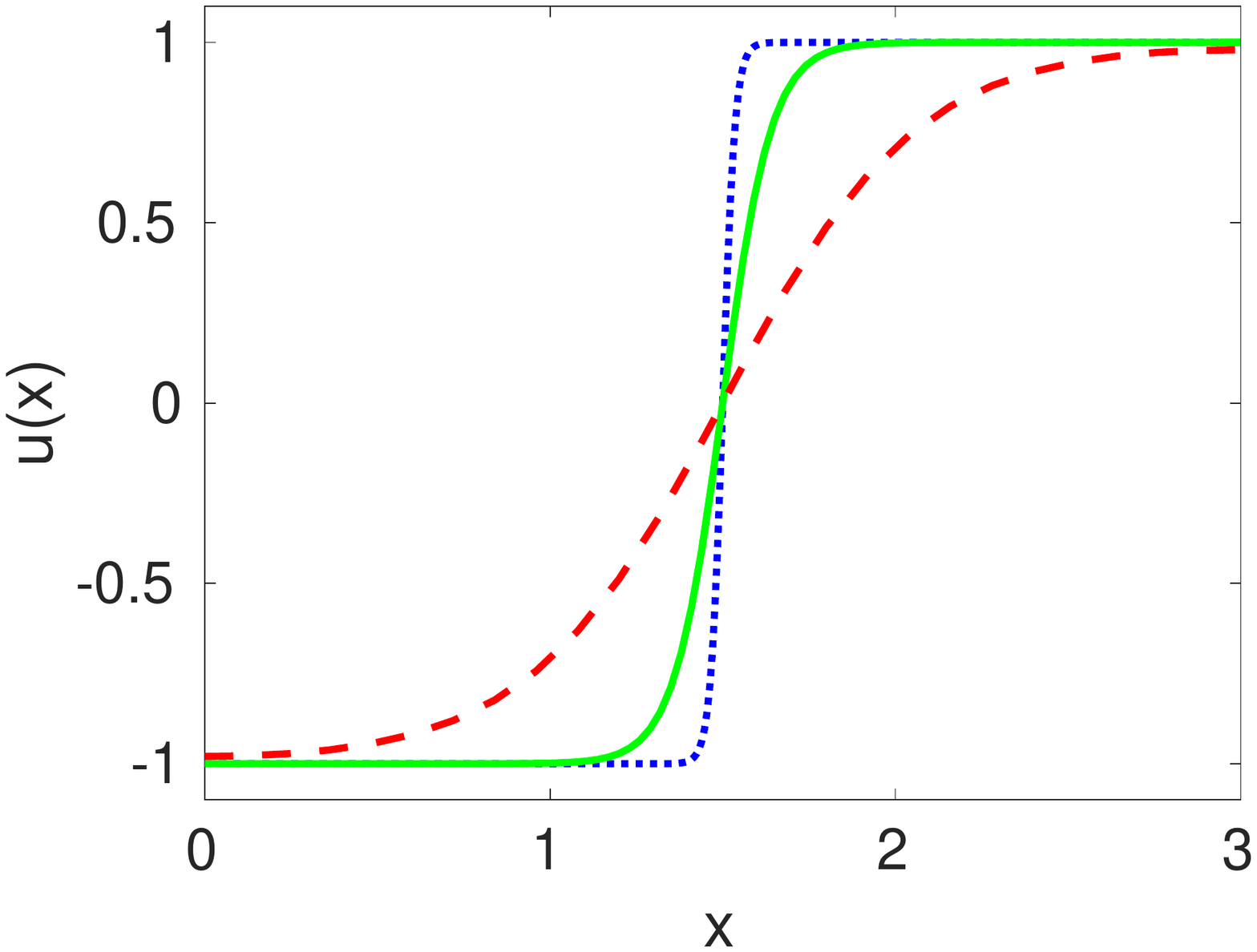}
\caption{$t=10.0$ }
\end{subfigure}
\caption{Stretching domain example as described in Section \ref{sec:deform_int}. Phase field for $\eps = 0.4$ (red), $\eps=0.1$ (green), and $\eps=0.025$ (blue). Simulation data are in Table \ref{tbl:deform_int} on the left.}
\label{stretchpm1}
\end{figure}

\begin{figure}
\begin{subfigure}{0.24\textwidth}
\includegraphics[clip, trim=0.5cm 6cm 0.5cm 6cm,width = \textwidth]{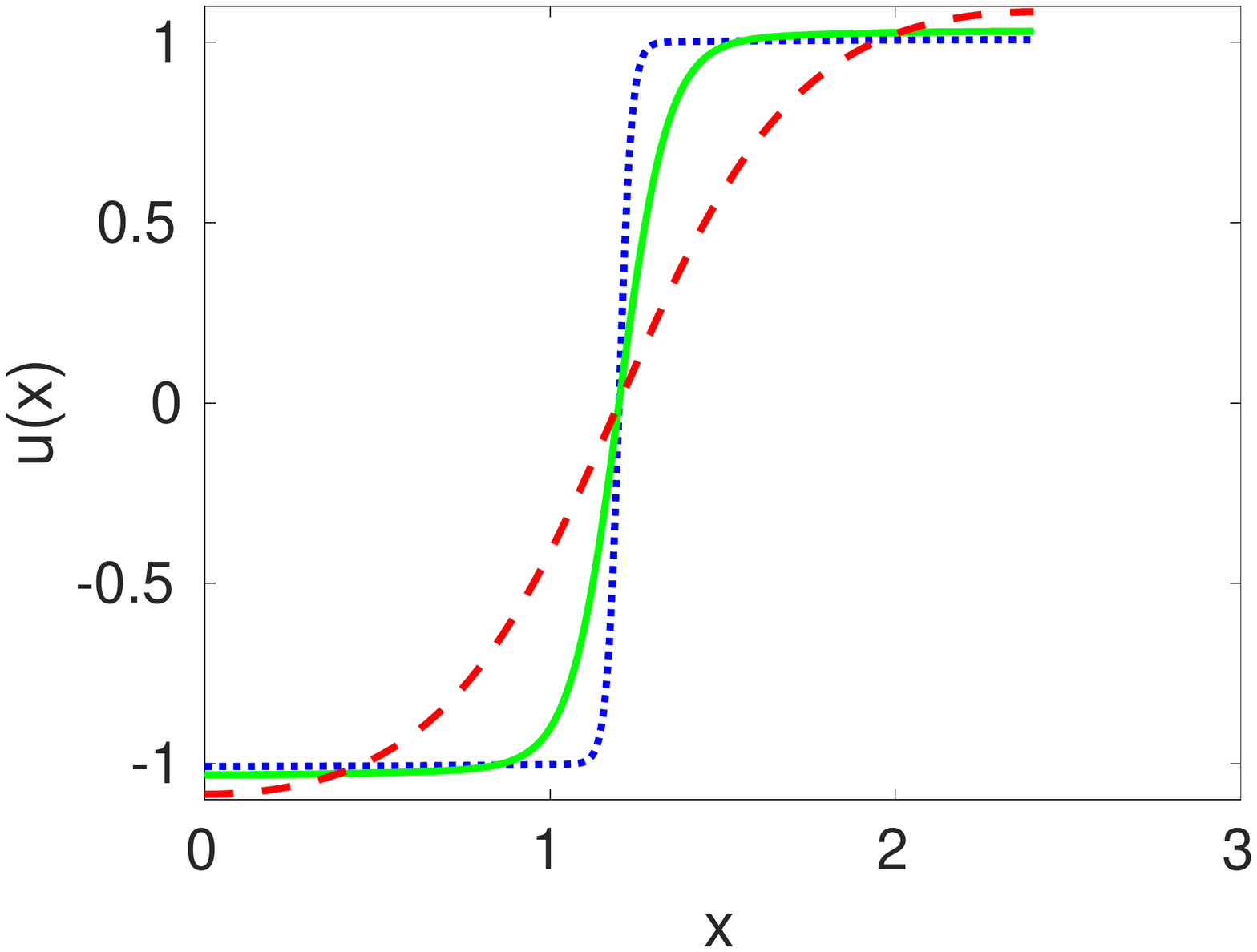}
\caption{$t=0.25$ }
\end{subfigure}
\begin{subfigure}{0.24\textwidth}
\includegraphics[clip, trim=0.5cm 6cm 0.5cm 6cm,width = \textwidth]{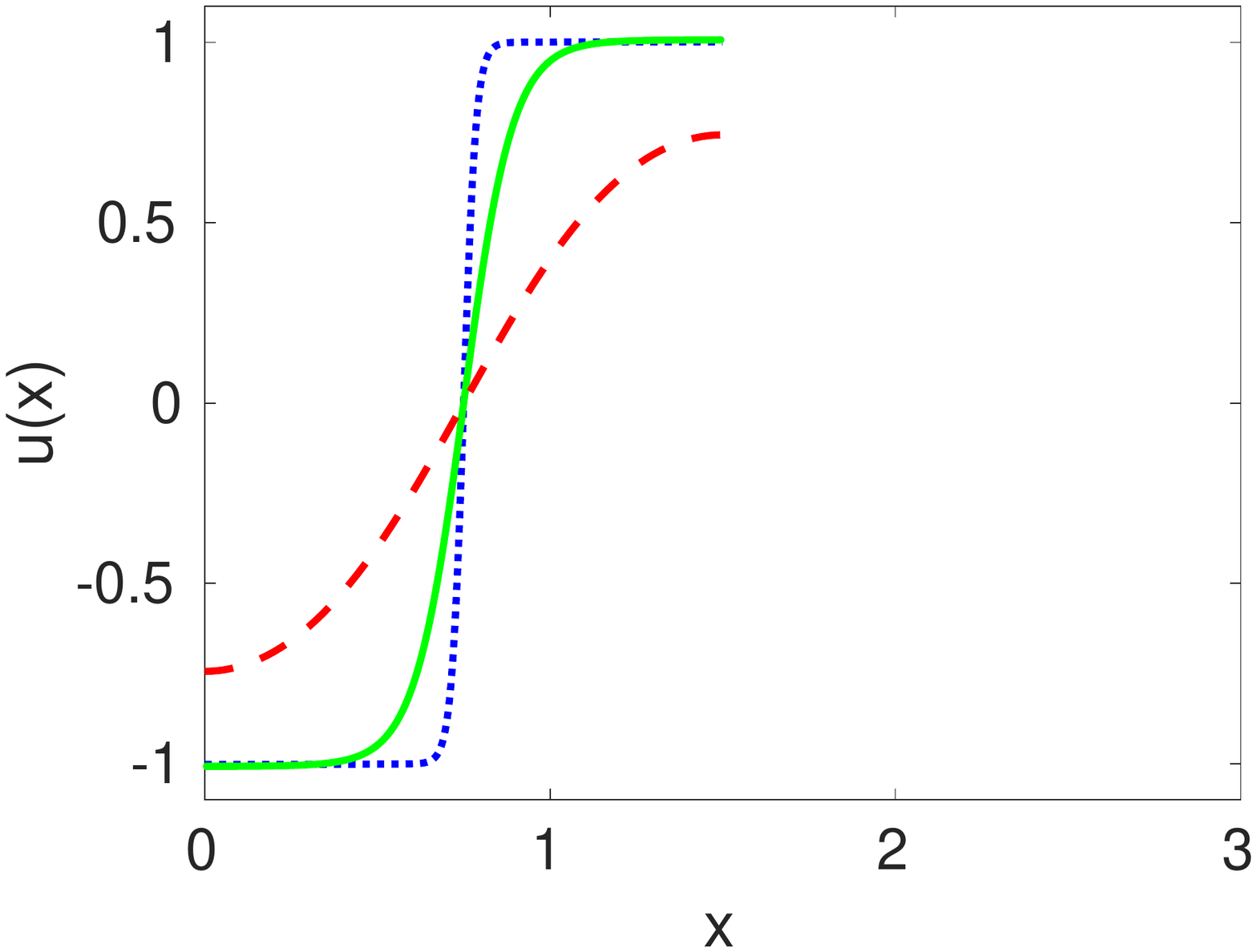}
\caption{$t=1.0$ }
\end{subfigure}
\begin{subfigure}{0.24\textwidth}
\includegraphics[clip, trim=0.5cm 6cm 0.5cm 6cm,width = \textwidth]{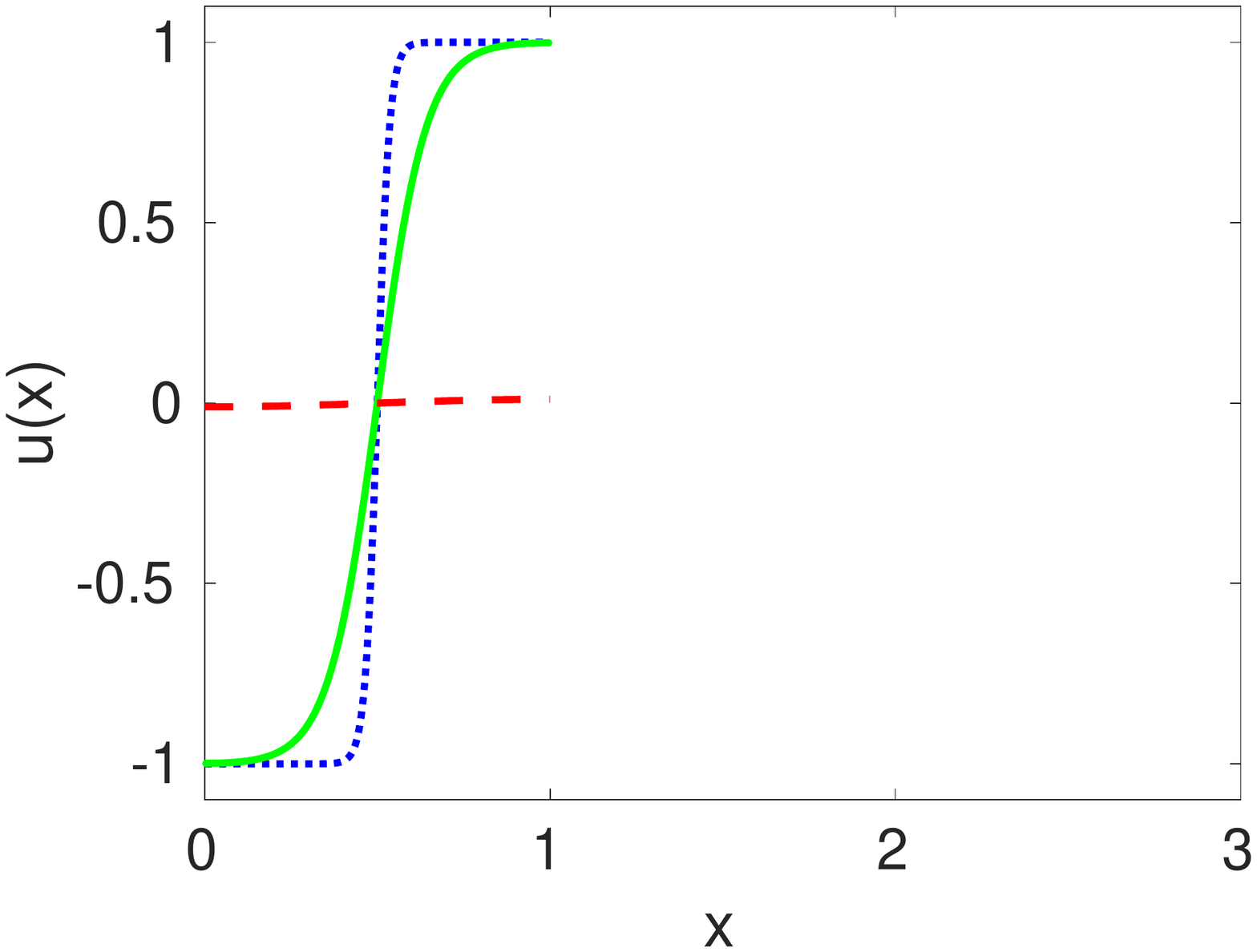}
\caption{$t=2.0$ }
\end{subfigure}
\begin{subfigure}{0.24\textwidth}
\includegraphics[clip, trim=0.5cm 6cm 0.5cm 6cm,width = \textwidth]{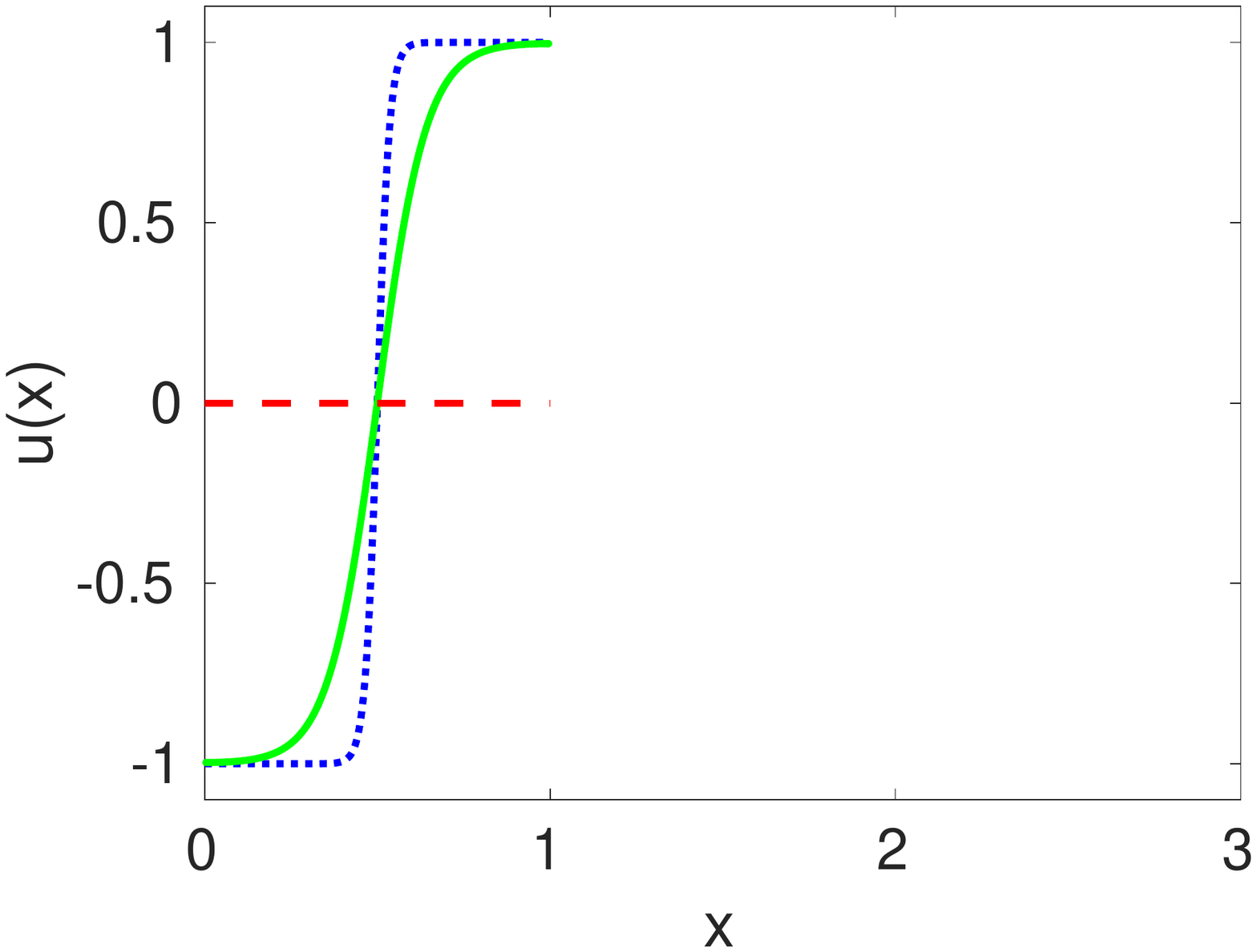}
\caption{$t=10.0$ }
\end{subfigure}
\caption{Compressing domain example as described in Section \ref{sec:deform_int}. Phase field for $\eps = 0.4$ (red), $\eps=0.1$ (green), and $\eps=0.025$ (blue). Simulation data are in Table \ref{tbl:deform_int} in the middle.}
\label{shrinkpm1}
\end{figure}

\begin{figure}
\begin{subfigure}{0.24\textwidth}
\includegraphics[clip, trim=0.5cm 6cm 0.5cm 6cm,width = \textwidth]{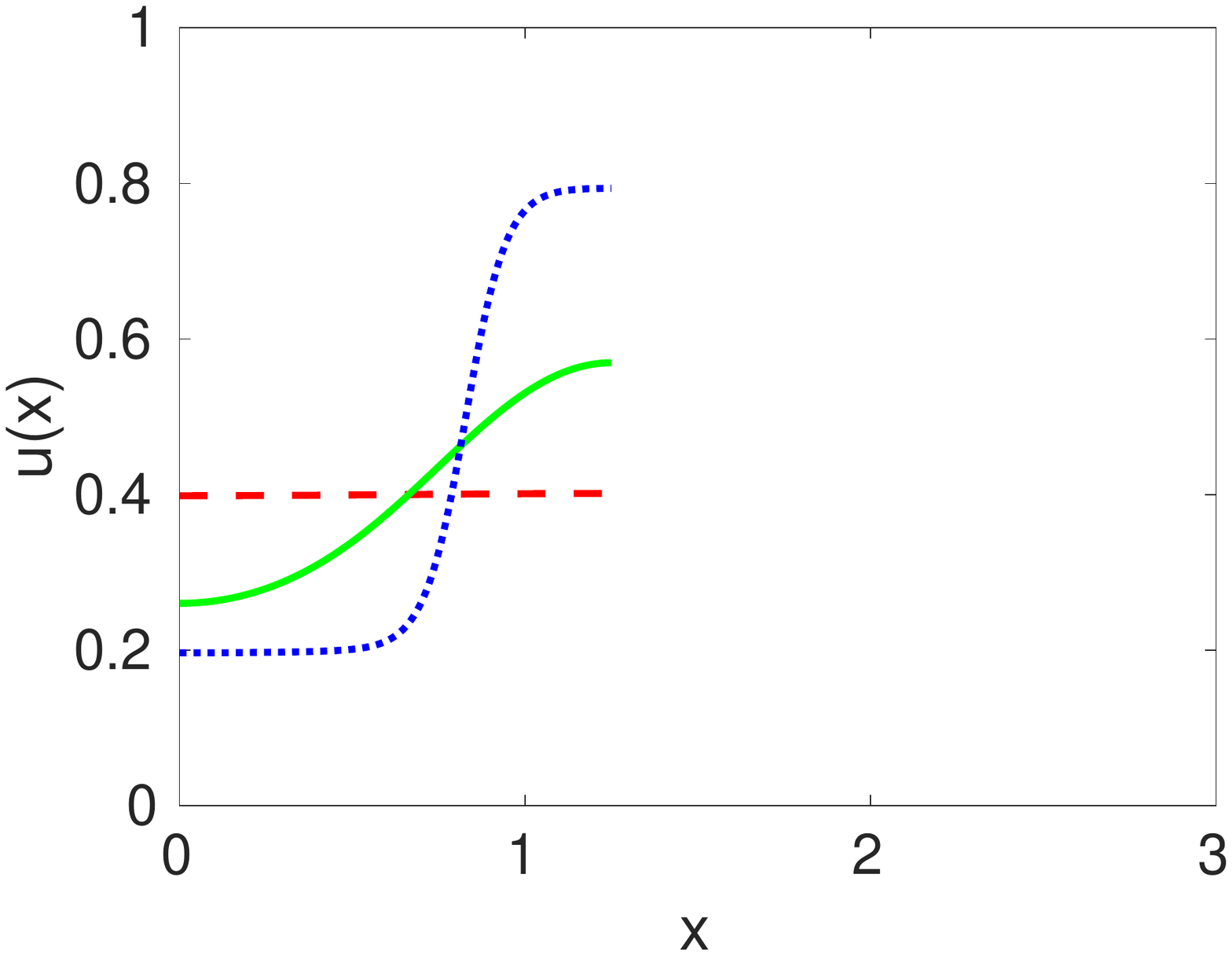}
\caption{$t=0.25$ }
\end{subfigure}
\begin{subfigure}{0.24\textwidth}
\includegraphics[clip, trim=0.5cm 6cm 0.5cm 6cm,width = \textwidth]{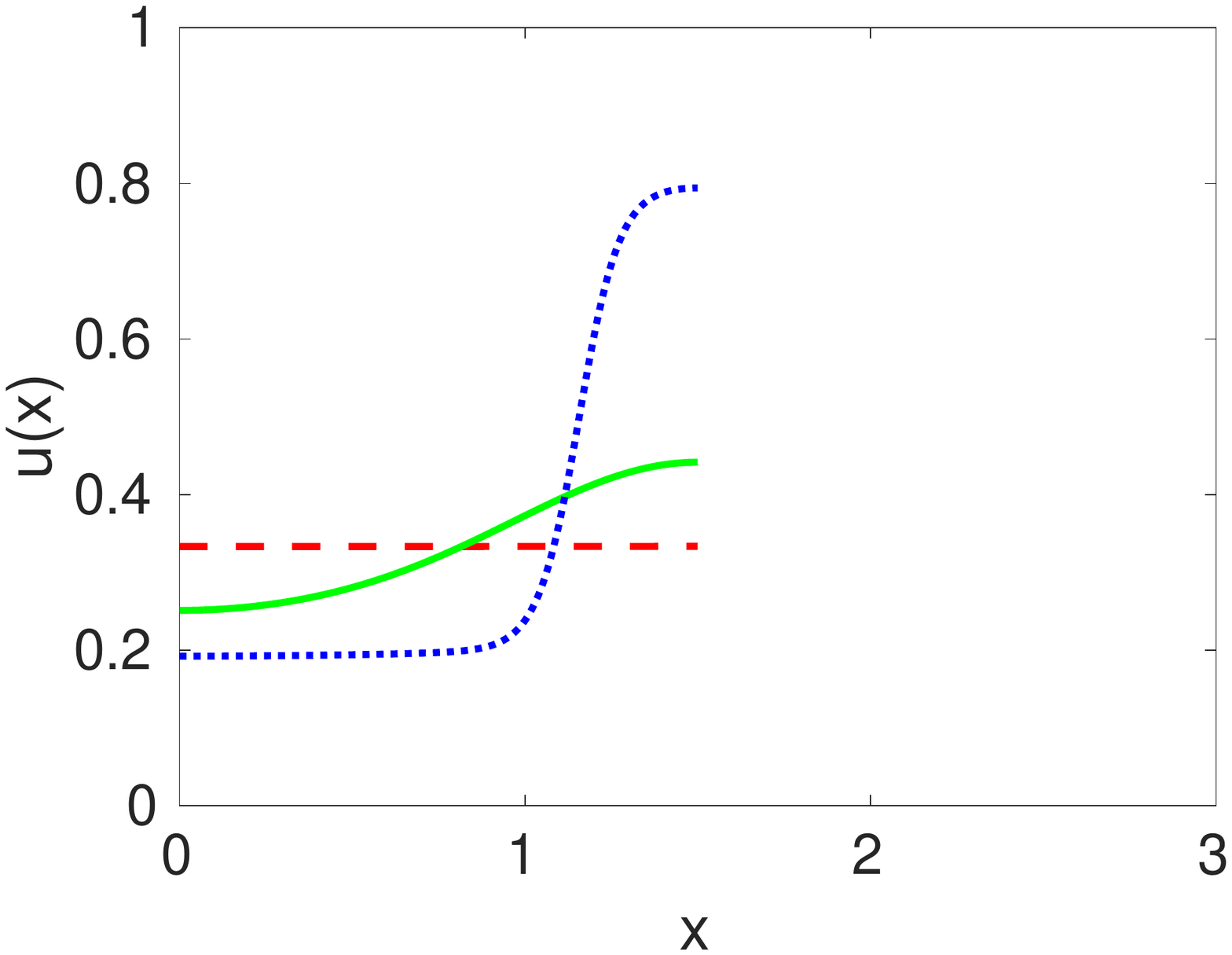}
\caption{$t=0.5$ }
\end{subfigure}
\begin{subfigure}{0.24\textwidth}
\includegraphics[clip, trim=0.5cm 6cm 0.5cm 6cm,width = \textwidth]{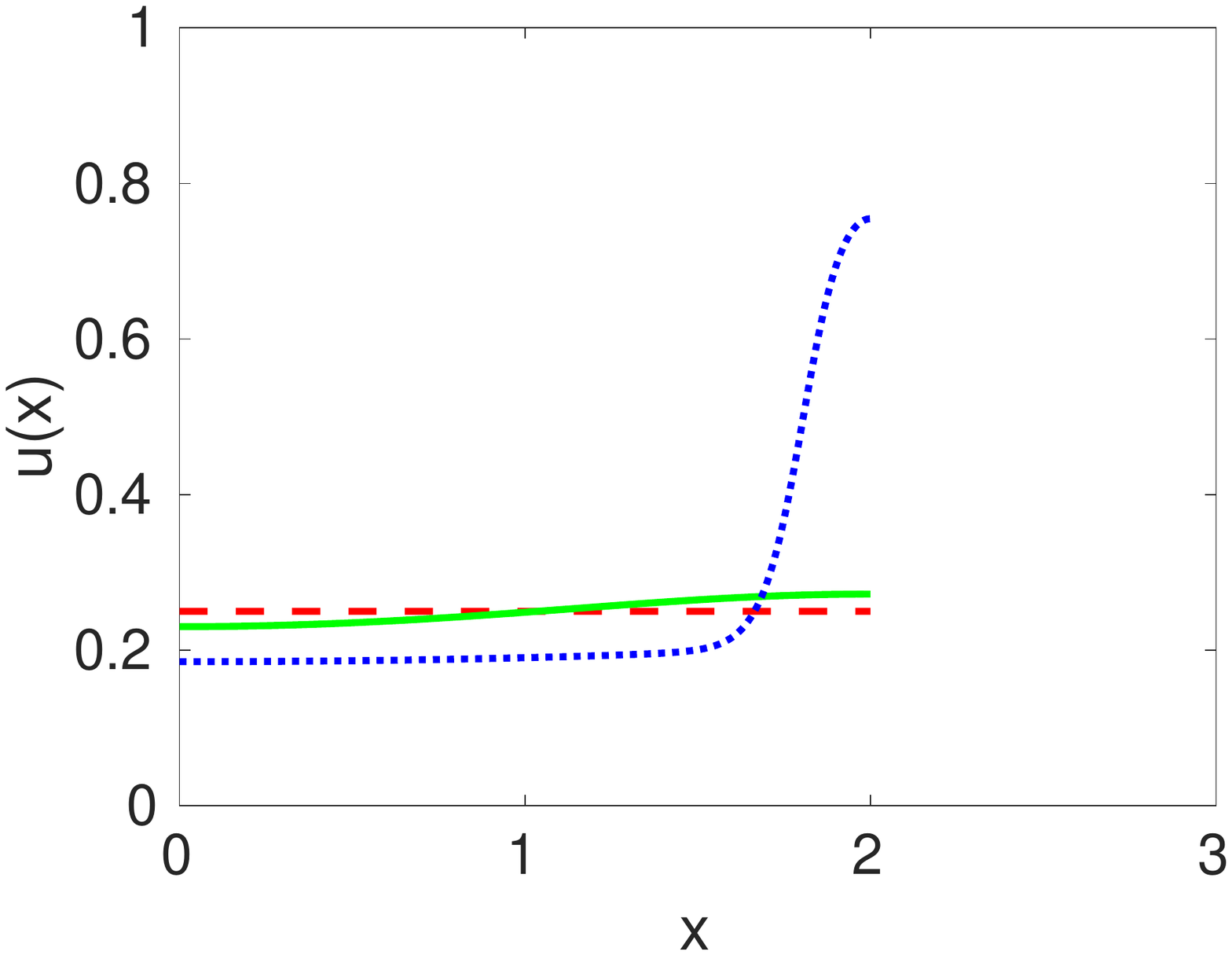}
\caption{$t=1.0$ }
\end{subfigure}
\begin{subfigure}{0.24\textwidth}
\includegraphics[clip, trim=0.5cm 6cm 0.5cm 6cm,width = \textwidth]{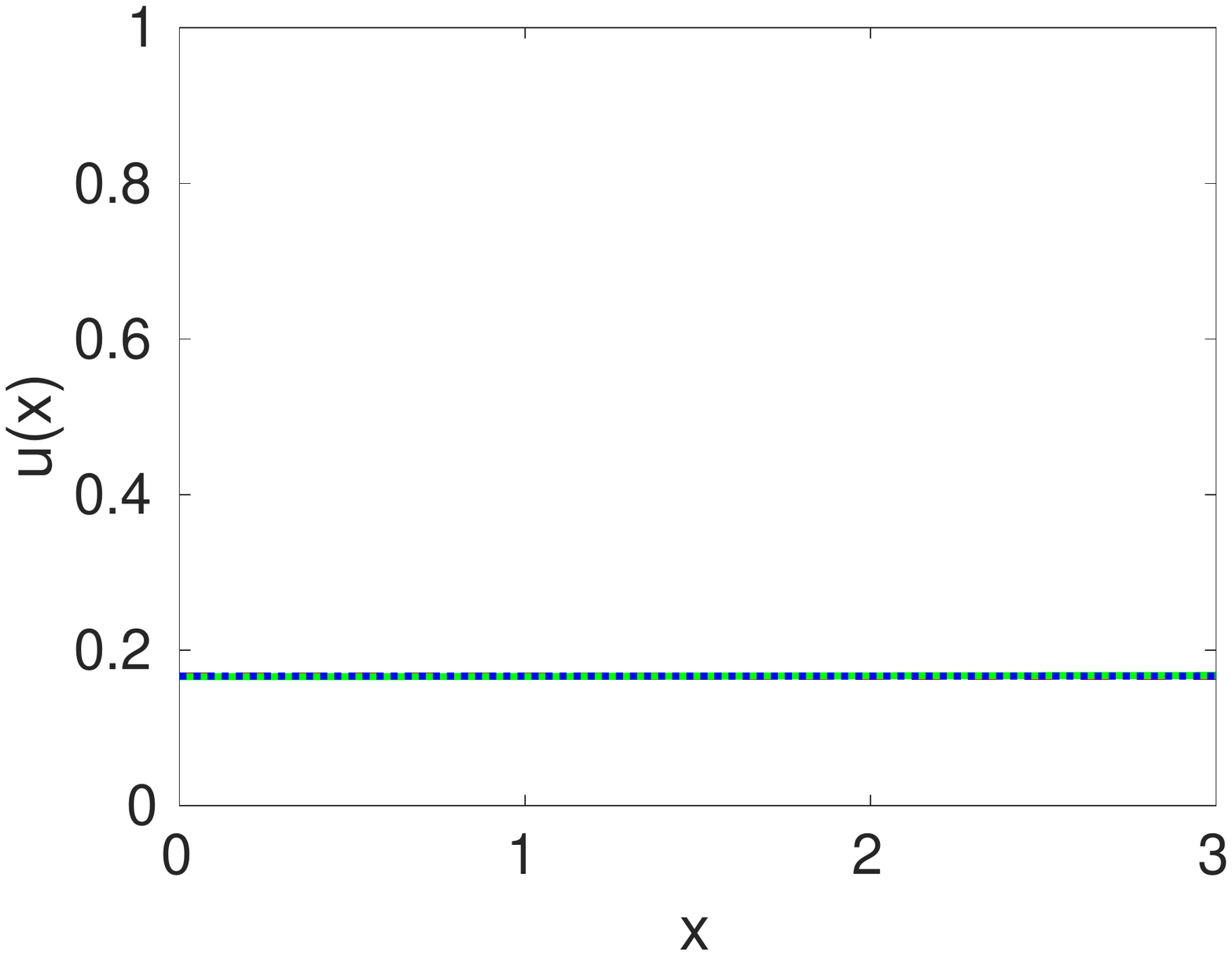}
\caption{$t=2.0$ }
\label{overstretch}
\end{subfigure}
\caption{Stretching domain example with positive minima of $F$ as described in Section \ref{sec:deform_int}. Phase field for $\eps = 0.4$ (red), $\eps=0.1$ (green), and $\eps=0.025$ (blue). Simulation data are in Table \ref{tbl:deform_int} on the right.}
\label{stretchpos}
\end{figure}

We first pick $u_a = -1$ and $u_b = 1$ and consider a phase transition at the centre of an interval. Then the interval is homogeneously stretched or compressed for a while, i.e., $\nabla_{\Gamma(t)} \cdot v(t)$ is constant in space. For the solution to \eqref{eq:FBP_bulk}, \eqref{eq:FBP_int} one will expect that the interface position moves in the direction of the deformation and, in the long term, ends up in the centre of the deformed interval. 

For our diffuse interface simulation we deform the domain as specified in Table \ref{tbl:deform_int} on the left and in the middle, respectively. At time $t=2$ the interval has reached the final length and we then further relax the profile of $u$ on the then stationary domain. 

In a first set of simulations we started with equilibrium $\tanh$ profiles defined in \eqref{eq:num_eq_profile} and shifted them such that they were symmetric with respect to the centre of the interval. In the short term, the advection leads to stretched or compressed profiles, respectively, but the effect becomes smaller the smaller $\eps$ is. In the long term, the profiles relax back to equilibrium profiles at the expected position which they do the faster the smaller $\eps$ is. 

In an attempt to investigate the robustness of the convergence behaviour we picked the profiles specified in Table \ref{tbl:deform_int} which is independent of $\eps$ for a second set of simulations. The results for different values of $\eps$ are displayed in Figure \ref{stretchpm1} and Figure \ref{shrinkpm1}, respectively, and indeed display the same long-term behaviour with one exception: for the largest $\eps = 0.4$ in the compression case the advection effect is so strong that we observe mixing of the phases, i.e., $u = 0$, in the long term. We remark that as only the largest $\eps$ is affected this is a \emph{finite $\eps$} effect which does not contradict the asymptotic result. 

We also examined the stretching example for a potential with minima at $u_a = 0.2$ and $u_b = 0.8$, see Table \ref{tbl:deform_int} on the right for the data and Figure \ref{stretchpos} for the simulation results. In this case the expansion makes the phase transitions vanish and leads to flat profiles which takes the longer the smaller $\eps$ is. Note that at time $t=2.0$ the (nearly) constant profiles of $u$ are slightly below the lower minimum $u_a = 0.2$ of the double-well potential. It thus does not satisfy the setting for the asymptotic analysis as specified at the beginning of Section \ref{sec:solregime}. Indeed, for the related sharp interface model the initial mass is
\[
 \mmm{M}(0) = \int_{\Gamma^+(t=0)} u_b + \int_{\Gamma^-(t=0)} u_a = \int_0^{0.5} u_b + \int_{0.5}^1 u_a = 0.5.
\]
If there was a solution to the sharp interface model which involves a phase transition its mass would satisfy 
\[
 \mmm{M}(2) = \int_{\Gamma^+(t=2)} u_b + \int_{\Gamma^-(t=2)} u_a \geq \int_0^3 u_a = 0.6 > \mmm{M}(0),
\]
which contradicts the mass conservation discussed around \eqref{eq:FBP_masscons}.

\subsection{Bulk Effects}
\label{sec:1D_bulk}

In the following examples we report on other effects due to velocity fields which do not have constant divergences and show some interesting behaviour.

First, recall that constant functions, different from $u_a$ and $u_b$, are unstable stationary solutions to the Cahn-Hilliard equation which also holds true for the ESCH equation with $v=0$. In our first example we start from a constant initial condition, $u_{IC}=\frac{u_a+u_b}{2}$, with $u_a>0$. We pick $v \geq 0$ as specified in the left column of Table \ref{tbl:1D_bulk}. Since mass is conserved, the advective effect of the velocity is expected to increase the mass where the domain is compressed and reduce the mass where the domain is stretched and, thus, is expected to induce a phase separation. Note that the boundary points of the domain $[0,1]$ do not move but internal movements take place, more precisely, stretching in $(0,0.5)$ and compressing in $(0.5,1)$.

\begin{figure}
\begin{subfigure}{0.24\textwidth}
\includegraphics[clip, trim=0.5cm 6cm 0.5cm 6cm,width = \textwidth]{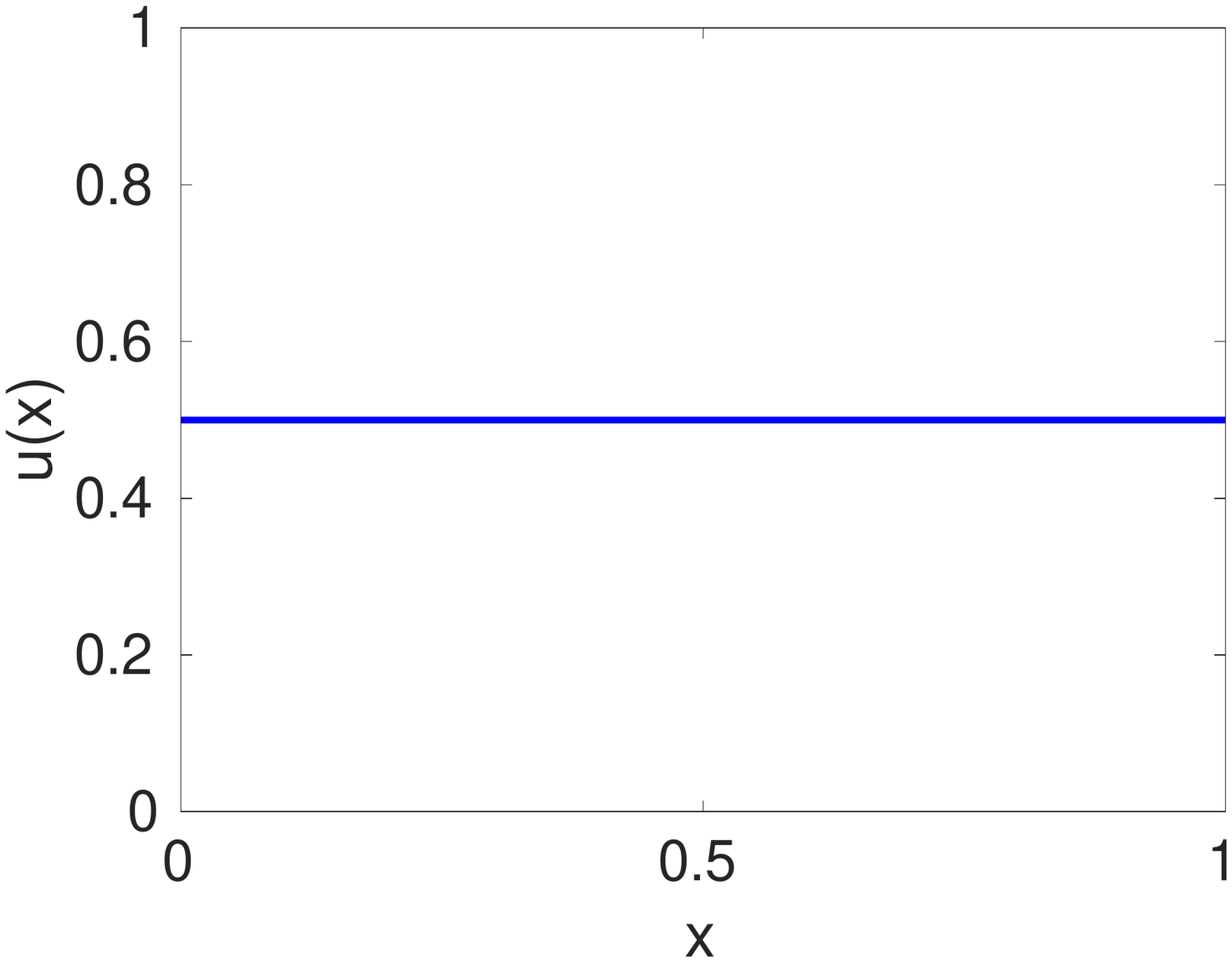}
\caption{$t=0.0$ }
\end{subfigure}
\begin{subfigure}{0.24\textwidth}
\includegraphics[clip, trim=0.5cm 6cm 0.5cm 6cm,width =\textwidth]{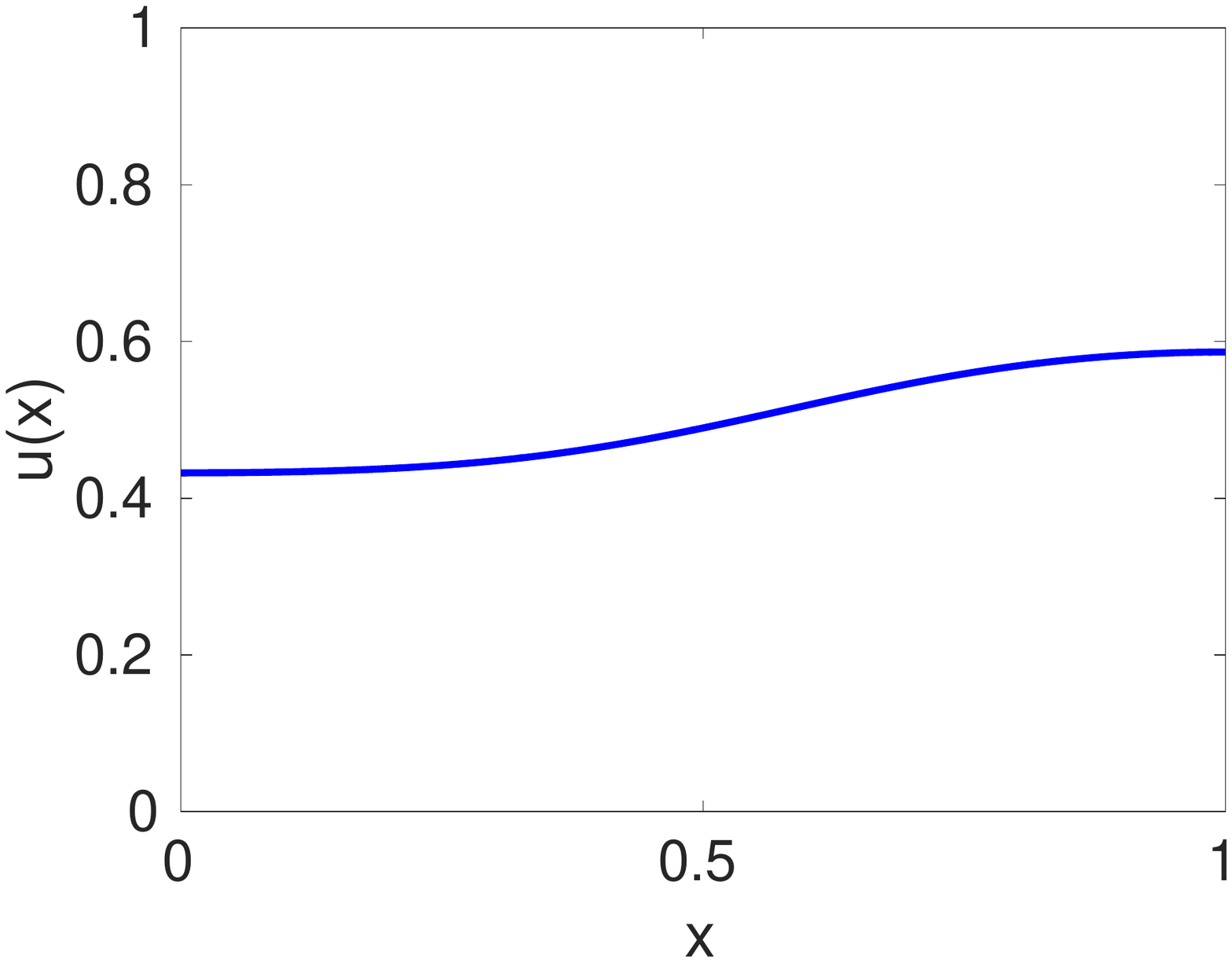}
\caption{$t=0.066$ }
\end{subfigure}
\begin{subfigure}{0.24\textwidth}
\includegraphics[clip, trim=0.5cm 6cm 0.5cm 6cm,width =\textwidth]{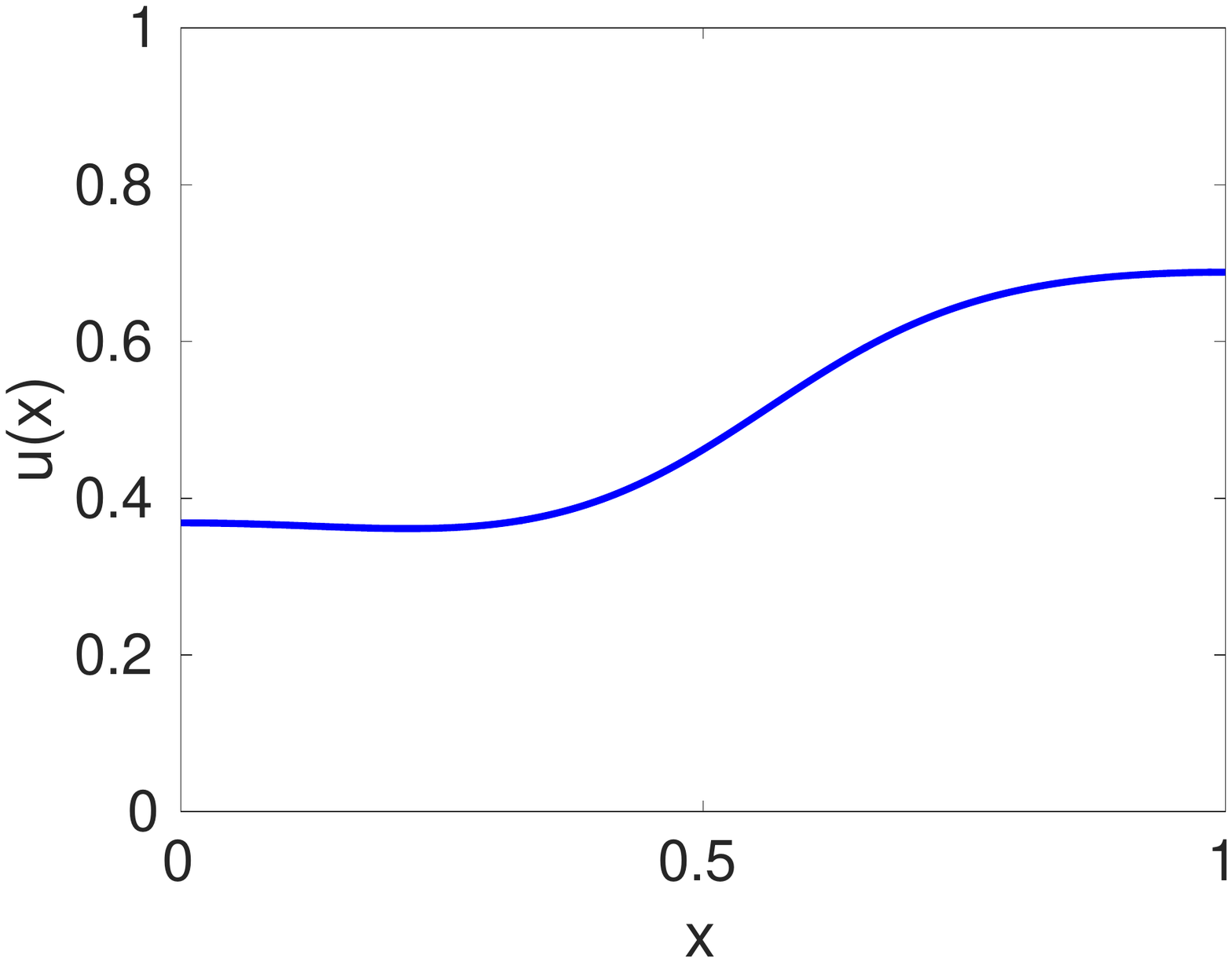}
\caption{$t=0.099$ }
\end{subfigure}
\begin{subfigure}{0.24\textwidth}
\includegraphics[clip, trim=0.5cm 6cm 0.5cm 6cm,width = \textwidth]{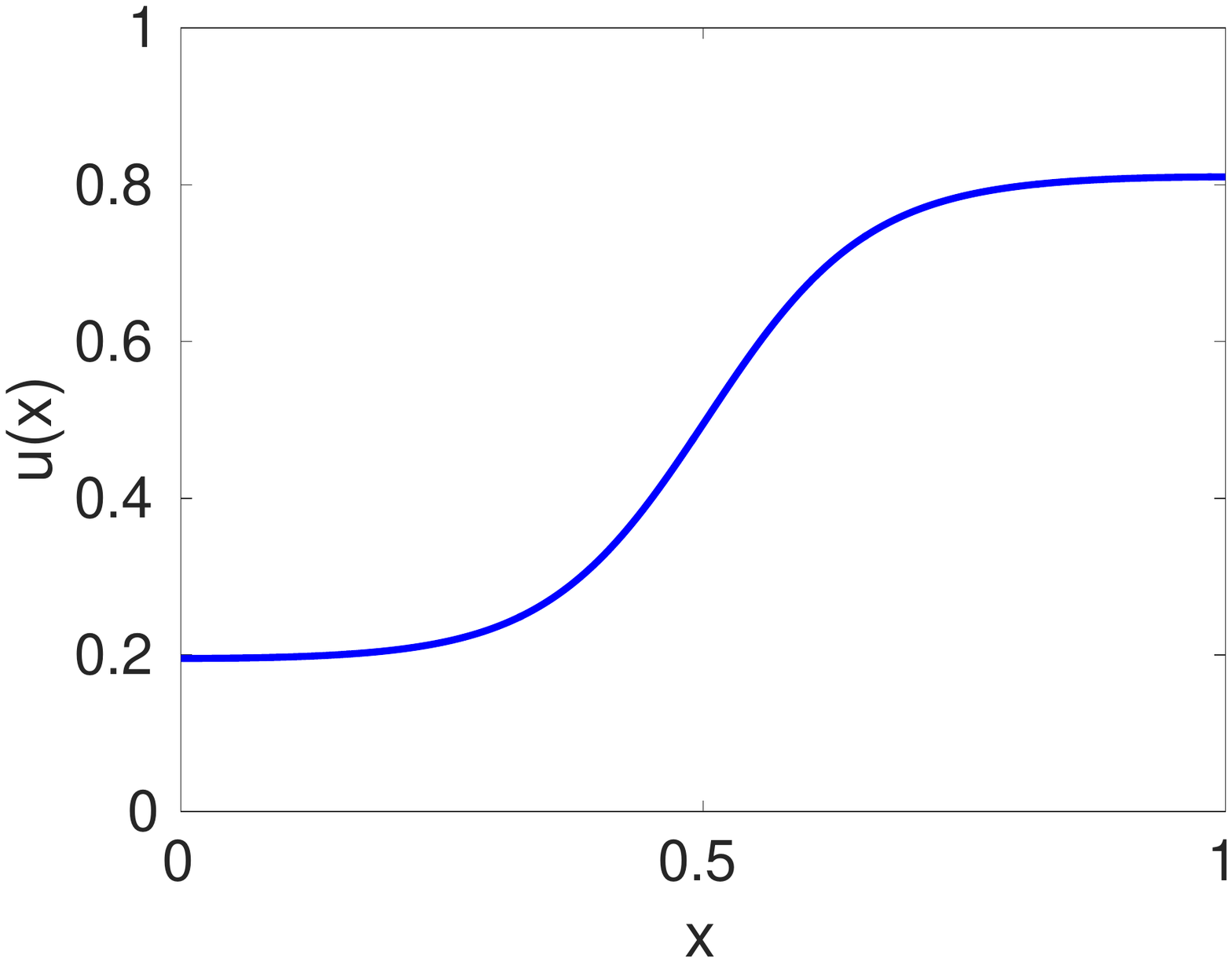}
\caption{$t=0.198$ }
\end{subfigure}
\caption{Generation of a phase interface by perturbing a flat initial profile as described in Section \ref{sec:1D_bulk}. $\eps = 0.033$, other simulation data are in Table \ref{tbl:1D_bulk} on the left.} \label{interfaceform}
\end{figure}

In Figure \ref{interfaceform} we see how the flat initial profile is perturbed by the advective effect of the velocity field such that a phase transition is obtained. The simulation data are given in Table \ref{tbl:1D_bulk} on the left. We remark that, in some cases, the velocity field from this example, has no destabilising effect. For instance, if $u_b = -u_a$ and $u_{IC} = 0$ then the solution remains constant at $u=0$ for all times.

In another example, initially, a phase interface is located at 0.25 within the initial domain $[0,1]$. We then extend the interval but such that $v=0$ in $[0,0.5]$ and $v(x,t) \neq 0$ only if $x > 0.5$, see Table \ref{tbl:1D_bulk} on the right for the details. 

Regarding the sharp interface model, \eqref{eq:FBP_bulk} implies that $w$ is no longer harmonic. Hence, the jump term $[M \sgrad w] \cdot \mu_\Lambda$ in \eqref{eq:FBP_int} changes and is expected to be non-zero. We thus expect a motion of the phase interface, $\Lambda$, in the direction of the stretching despite the surface velocity, $v$, being zero in the region containing the interface.

\begin{figure}
\includegraphics[clip, trim=0.5cm 6cm 0.5cm 6cm,width = 0.32\textwidth]{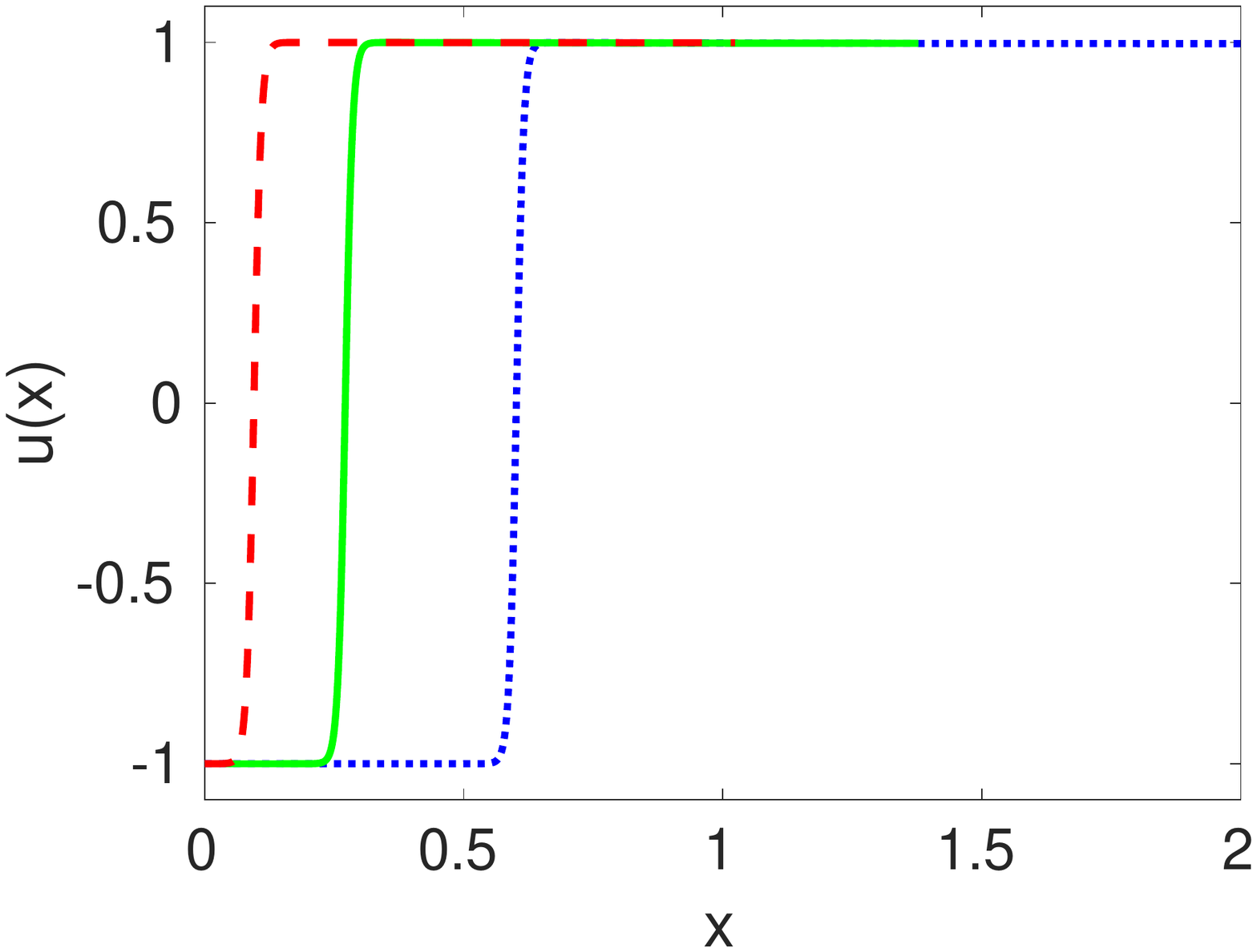}
\caption{Phase interface motion induced by bulk effects away from the interface as described in Section \ref{sec:1D_bulk}. $\eps =0.01$, $t=0.1$ (red), $t=1.0$ (green), $t=1.8$ (blue). Simulation data are in Table \ref{tbl:1D_bulk}.}
\label{bulkproc}
\end{figure}

In Figure \ref{bulkproc} we can see that there is indeed a motion induced by the non-trivial bulk problem. In addition, once the phase interface gets beyond the point 0.5, its velocity can be seen to increase. This is in accordance with the last equation in \eqref{eq:FBP_int} as $v \geq 0$ there.

\begin{table}
 \tabcolsep=3pt
 \begin{tabular}{lll}
  \hline
  Parameter & Data for Figures \ref{interfaceform}  & Data for Figure \ref{bulkproc} \\
  \hline
  $u_a$, $u_b$; $\bar{M}; T$ & 0.2,  0.8;  1;  0.2 & -1,  1;  1;  2\\
  $\Gamma(t)$ & $[0,1]$ & $[0,\cot ^{-1}(1.83-t)+0.5]$ \\
  $v(x,t)$, $x \in \Gamma(t)$ & $\sin(\pi x) $ & $\left\{\begin{array}{ll}\sin^2( x-\frac{1}{2}) & x \geq \frac{1}{2} \\ 0 & x < \frac{1}{2} \end{array} \right.$ \\
  $u_0(x)$ & $0.5$  & $ \tanh(\frac{x-0.25}{\eps})$ \\
  \hline
 \end{tabular}
 \caption{Simulation data for Section \ref{sec:1D_bulk}.} \label{tbl:1D_bulk} 
\end{table}

\subsection{A Solution on a Sphere with Tangential Mass Transport}
\label{sec:MovingSphere}

Considering 2D hypersurfaces in 3D allows us to demonstrate the geometric behaviours of solutions to the ESCH equation and to discuss effects due to the geodesic curvature which appears in \eqref{eq:FBP_int}.

In this example we consider a sphere with a tangential velocity field $v$ so that the shape doesn't change. As in \cite{R15} we look for solutions to \eqref{eq:FBP_bulk}, \eqref{eq:FBP_int} which are rotationally symmetric and, thus, are constant in the azimuthal angle $\phi \in (0,2\pi)$ but only depend on the polar angle $\theta \in (0, \pi)$, i.e., 
\begin{equation} \label{eq:ansatz_w}
w(x(\theta,\phi),t) = W(\theta,t) \quad \mbox{where } x(\theta,\phi) = \left(\sin\theta\cos\phi,\,\sin\theta\sin\phi,\,\cos\phi\right)^T.
\end{equation}
The difference to \cite{R15} is the presence of the velocity field $v$. We pick a velocity field which transports mass from the north to the south pole, $v(x(\theta,\phi),t) = \bar{v} \sin(\theta) x_{\theta}(\theta,\phi)$ with some $\bar{v} > 0$. One can easily show that $\sgrad \cdot v = 2 \bar{v} \cos (\theta)$. 

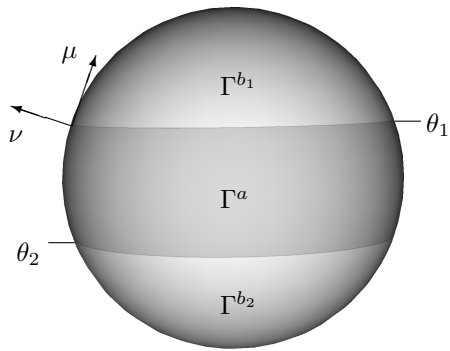
\begin{figure}
\vspace{20mm}
\begin{picture}(100,100)
\put(-50,-40){\includegraphics[scale=0.5]{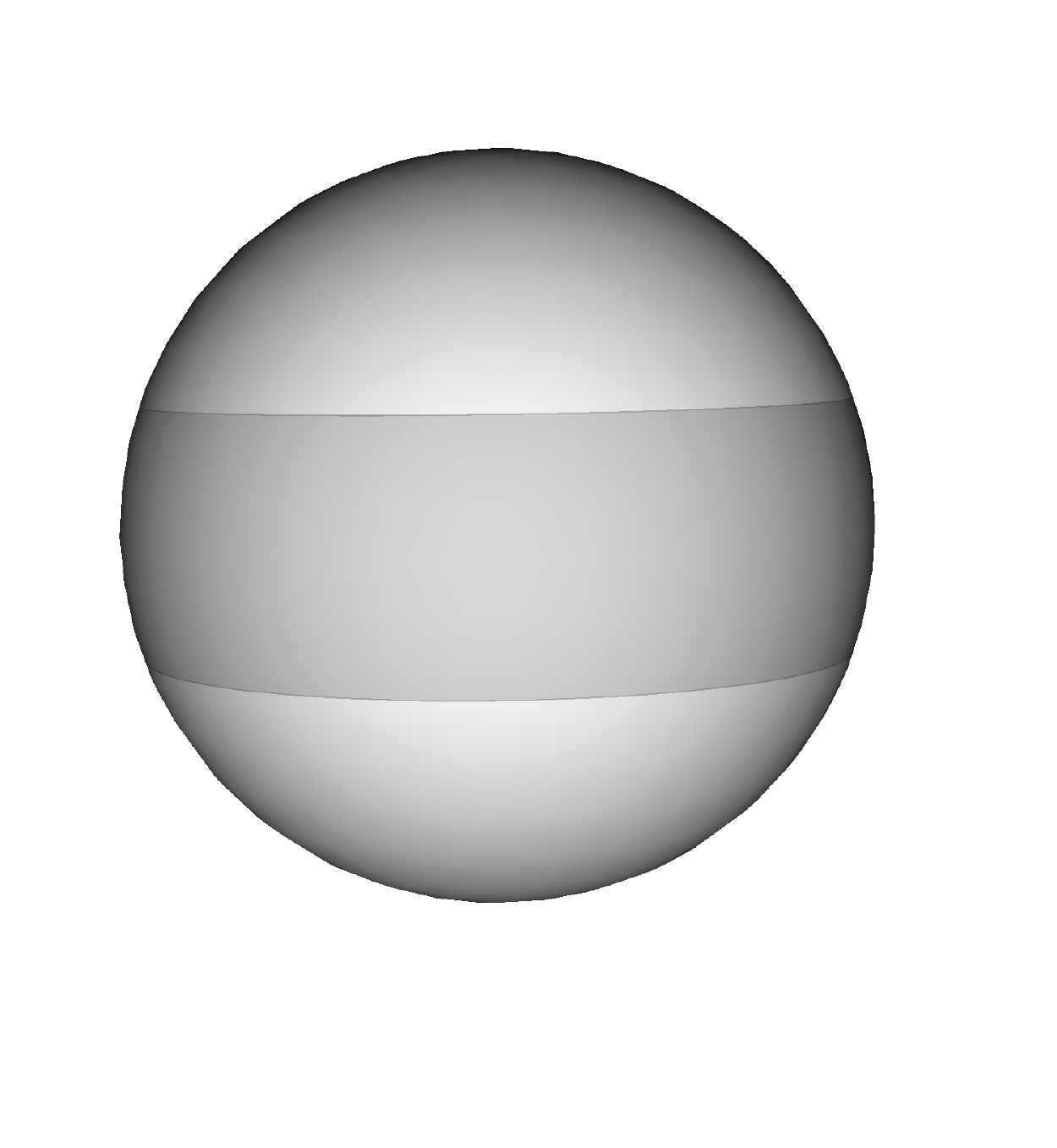}}
\put(30,15){$\Gamma^{b_2}$}
\put(30,98){$\Gamma^{b_1}$}
\put(30,55){$\Gamma^a$}
\put(108,84){$\theta_1$}
\put(106,88){\line(-1,0){10.}}
\put(-47,35){$\theta_2$}
\put(-35,42){\line(1,0){10.}}
\put(-26,86){\vector(-3,1){23.}}
\put(-26,86){\vector(1,3){9.}}
\put(-50,80){$\nu$}
\put(-30,112){$\mu$}
\end{picture}
\caption{Setup for the example in Section \ref{sec:MovingSphere}.}
\label{SetupMovingSphere}
\end{figure}

We then consider two distinct regions around the poles, where $u=u_b$, which are separated by a band where $u=u_a$, see Figure \ref{SetupMovingSphere}. We will refer to the inner region as $\Gamma^a(t)$ and and the two caps as $\Gamma^{b_{1,2}}(t)$. We denote by $\theta_{{1,2}}(t)$ the polar angle of the boundaries between $\Gamma^{b_{1,2}}(t)$ and $\Gamma^a(t)$, respectively.

With the ansatz \eqref{eq:ansatz_w} the Laplace-Beltrami operator applied to $w$ becomes $$\slap w = \frac{1}{\sin\theta}\frac{\partial}{\partial \theta} \left( \sin\theta \frac{\partial W}{\partial \theta} \right).$$ 
The general solution $W^{(i)}(\theta,t)$ to \eqref{eq:FBP_bulk} in $\Gamma^i(t)$, $i \in \{ b_1,b_2,a \}$, then is
\begin{equation} \label{eq:anasolsymm}
W^{(i)} (\theta,t) = c_1^{(i)}(t) \log\left[\tan\left(\frac{\theta}{2}\right)\right] - \frac{u_i \bar{v}}{\bar{M}} \cos(\theta) + c_2^{(i)}(t), \quad 
\theta \in 
\begin{cases}
 (0,\theta_1(t)), & \, i = b_1, \\
 (\theta_1(t), \theta_2(t)), & \, i = a, \\
 (\theta_2(t),\pi), & \, i = b_2,
\end{cases}
\end{equation}
with functions $c^{(i)}_k(t)$, $k = 1,2$, which will be determined by the interface conditions.

Assuming a smooth solution in $\Gamma^{b_{1,2}}(t)$ the gradient has to be zero at the poles which implies that $c_1^{(b_1)}(t) = c_1^{(b_2)}(t) = 0$. We now use the second equation of \eqref{eq:FBP_int} and that the geodesic curvature of the phase interface is equal to $(-1)^{k+1} \cot(\theta_{{k}}(t))$, $k=1,2$:
\begin{equation*}
c_2^{(b_k)}(t) = (-1)^{k+1} S \cot (\theta_k(t)) + \frac{\bar{v}}{\bar{M}} u_b\cos(\theta_k(t)), \quad k=1,2.
\end{equation*}
We can use the same boundary condition on each boundary of $\Gamma^a(t)$ in order to determine $c_1^{(a)}(t)$ and $c_2^{(a)}(t)$. We only include $c_1^{(a)}(t)$ below as the formula for $c_2^{(a)}(t)$ is not needed to progress:
\begin{equation*}
c_1^{(a)}(t) = \frac{S \left[ \cot\left((\theta_1(t)\right) + \cot\left(\theta_2(t)\right) \right] + \tfrac{\bar{v}}{\bar{M}}u_a (\cos(\theta_1(t)) - \cos(\theta_2(t)))}{\log\left[ \tan(\theta_1(t)/2)\right] - \log\left[ \tan(\theta_2(t)/2)\right]}
\end{equation*}

Having expressed the solution \eqref{eq:anasolsymm} in terms of the $\theta_k(t)$ we can use the third equation of \eqref{eq:FBP_int} in order to derive a system of ODEs for the $\theta_k(t)$, k=1,2 (note that $v_\Lambda(\theta_k(t)) \cdot \mu(\theta_k(t)) = (-1)^k \theta_k'(t)$):
\begin{equation*}
\theta_1'(t) =  \frac{\bar{M} \tilde{c}_1^{(a)}(\theta_1(t),\theta_2(t))}{(u_b-u_a) \sin(\theta_1(t))}, \quad \theta_2'(t) =  \frac{\bar{M} \tilde{c}_1^{(a)}(\theta_1(t),\theta_2(t))}{(u_b-u_a) \sin(\theta_2(t))}. 
\end{equation*}
where $\tilde{c}_1^{(a)}(\theta_1(t),\theta_2(t))=c_1^{(a)}(t)$.

We choose the quartic potential $F(u) = \frac{1}{4} (u^2-1)^2$, i.e., $u_a = -1$, $u_b = 1$, for which $S = \frac{\sqrt{2}}{3}$ in \eqref{eq:FBP_int} and for the initial condition of the sharp interface problem set $\theta_1(0) = 0.8$ and $\theta_2(0) = 2.1$ so that $\Gamma^{b_2}(0)$ is slightly bigger than $\Gamma^{b_1}(0)$. 

For the initial condition of the Cahn-Hilliard equation we use
\begin{equation*}
u_0(\theta) = 
\left\{\begin{array}{cc}
 \tanh\left(\frac{0.8-\theta}{\eps\sqrt{2}}\right), & \theta<1.45,  \\
 \tanh\left(\frac{\theta-2.1}{\eps\sqrt{2}}\right), & \theta\geq 1.45.
\end{array}\right. 
\end{equation*}

\begin{figure}
\begin{subfigure}{0.24\textwidth}\centering
\includegraphics[width = \textwidth,trim = 200 600 200 45, clip]{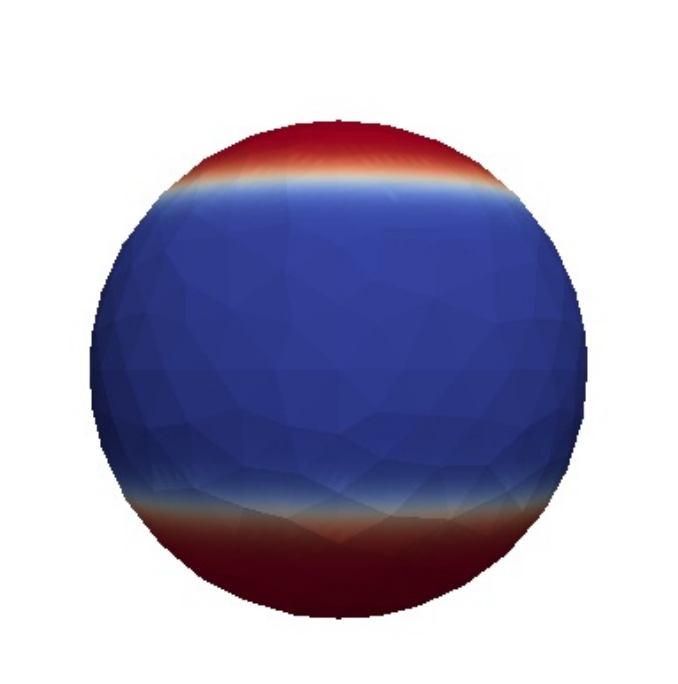}
\caption{$t=0$  }
\end{subfigure}
\begin{subfigure}{0.24\textwidth}\centering
\includegraphics[width = \textwidth,trim = 200 600 200 45, clip]{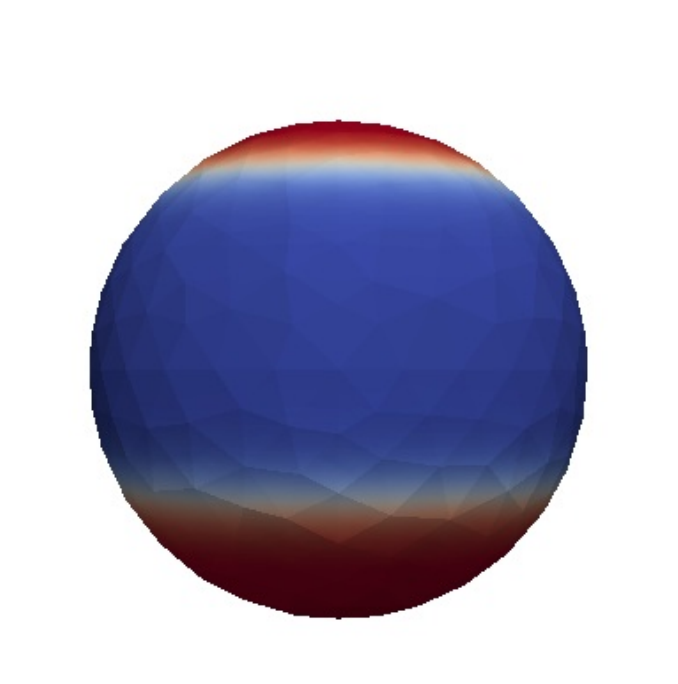}
\caption{$t=0.05$,  }
\end{subfigure}
\begin{subfigure}{0.24\textwidth}\centering
\includegraphics[width = \textwidth,trim = 200 600 200 45, clip]{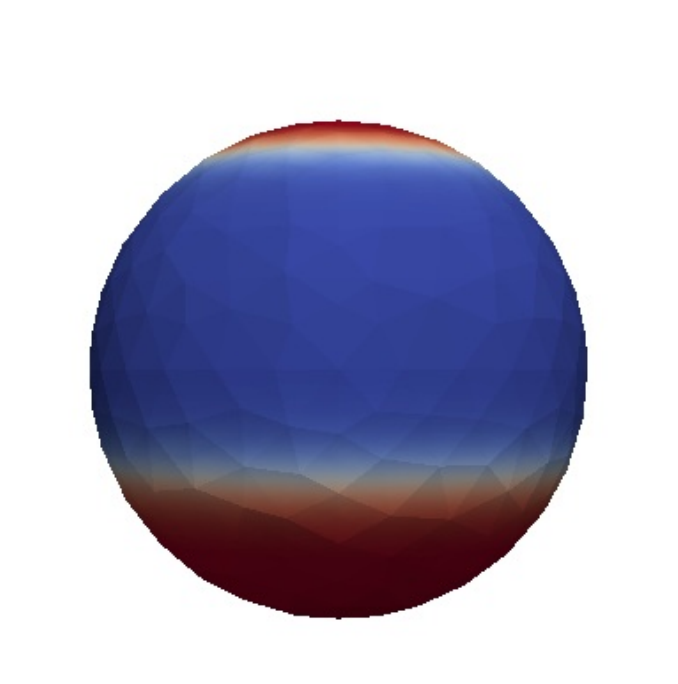}
\caption{$t=0.1$ }
\end{subfigure}
\begin{subfigure}{0.24\textwidth}\centering
\includegraphics[width = \textwidth,trim = 200 600 200 45, clip]{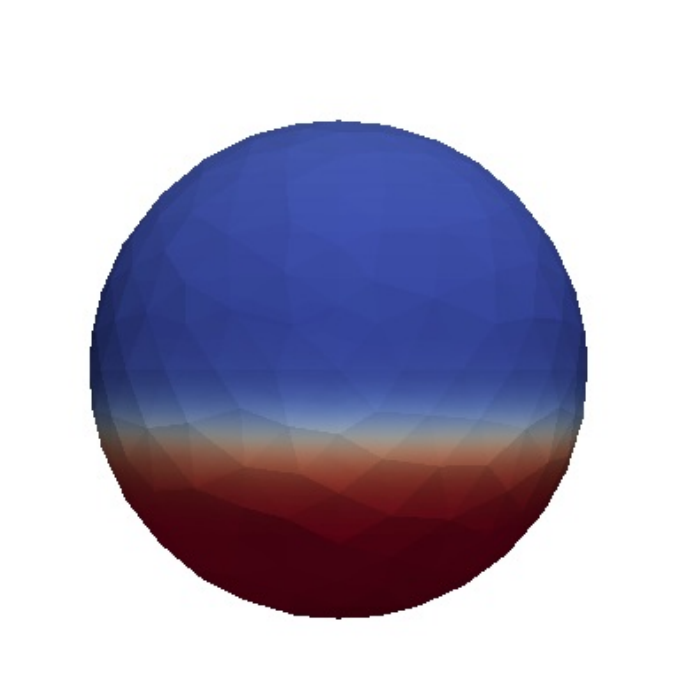}
\caption{$t=0.15$ }
\end{subfigure}
\caption{Coarsening example on the sphere as described in Section \ref{sec:MovingSphere}, $\eps = 0.1$, $\bar{v} = 0$, $\bar{M} = 5$.} \label{coarse}
\end{figure}

\begin{figure}
\begin{subfigure}{0.24\textwidth}\centering
\includegraphics[width = \textwidth,trim = 200 600 200 45, clip]{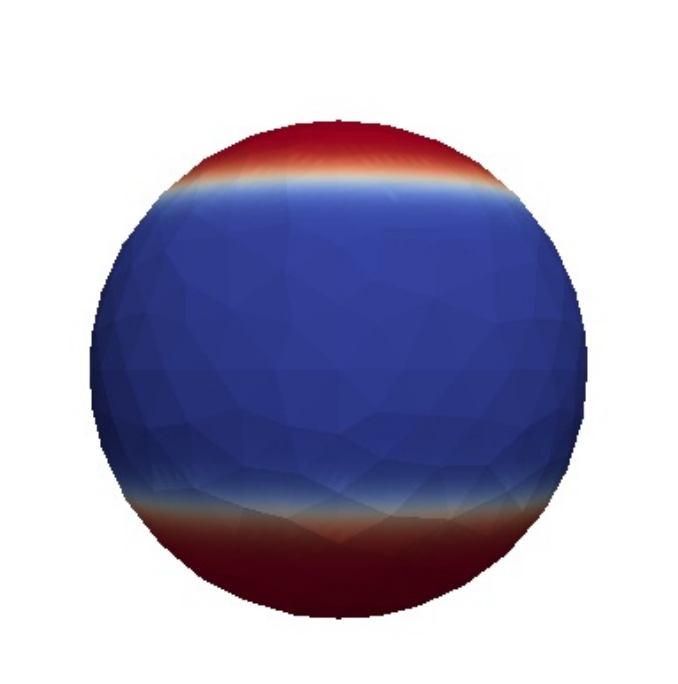}
\caption{$t=0$  }
\end{subfigure}
\begin{subfigure}{0.24\textwidth}\centering
\includegraphics[width = \textwidth,trim = 200 600 200 45, clip]{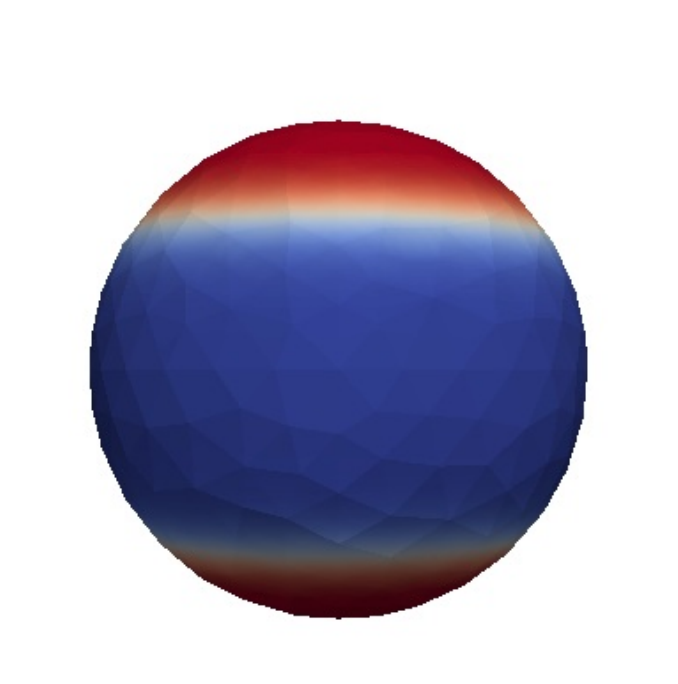}
\caption{$t=0.05$  }
\end{subfigure}
\begin{subfigure}{0.24\textwidth}\centering
\includegraphics[width = \textwidth,trim = 200 600 200 45, clip]{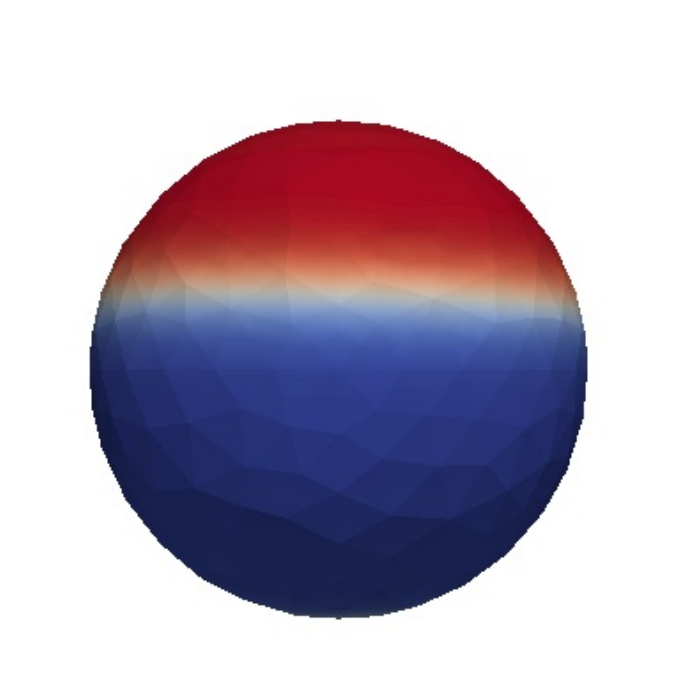}
\caption{$t=0.1$  }
\end{subfigure}
\begin{subfigure}{0.24\textwidth}\centering
\includegraphics[width = \textwidth,trim = 200 600 200 45, clip]{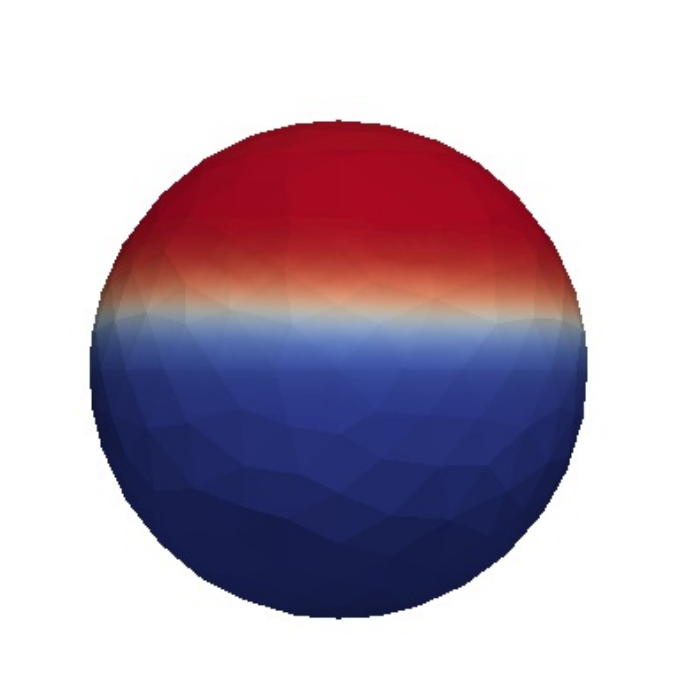}
\caption{$t=0.15$  }
\end{subfigure}
\caption{Example with tangential mass transport on the sphere as described in Section \ref{sec:MovingSphere}, $\eps = 0.1$, $\bar{v} = 10$, $\bar{M} = 5$.} \label{reverse}
\end{figure}

In the case $\bar{v} = 0$, i.e., without any mass transport, we expect the solution to coarsen to a two region solution with the area around the southern pole, $\theta = \pi$, taking the phase value $u_b$. This is indeed what we observe, see Figure \ref{coarse}. In turn, if the mass transport towards the south pole with a tangential velocity field is strong enough we expect that again a two region solution emerges but with the domain of the phase $\{ u = u_b \}$ around the northern pole, $\theta = 0$. For $\bar{v} = 10$, Figure \ref{reverse} displays that solutions indeed exhibit this behaviour.

We want to compare our solution to the sharp interface model with solutions of the Cahn-Hilliard equation by considering the energy of the system. The Ginzburg-Landau energy \eqref{GLEF} is the energy for the diffuse interface model and, as shown in \cite{N08}, converges to the energy of the sharp interface model which is proportional to the length of the phase interface:
\[ 
\mathcal{E}_{\eps} \to 2 S \; \text{length}(\Lambda) =: \mathcal{E}_0
\] 
which here amounts to
\begin{equation} \label{SIM:Energy}
\mathcal{E}_0(t) = \frac{4 \sqrt{2} \pi}{3} \left[ \sin(\theta_1(t)) + \sin(\theta_2(t)) \right].
\end{equation}

\begin{figure}
\begin{subfigure}{0.48\textwidth}
\includegraphics[width = \textwidth]{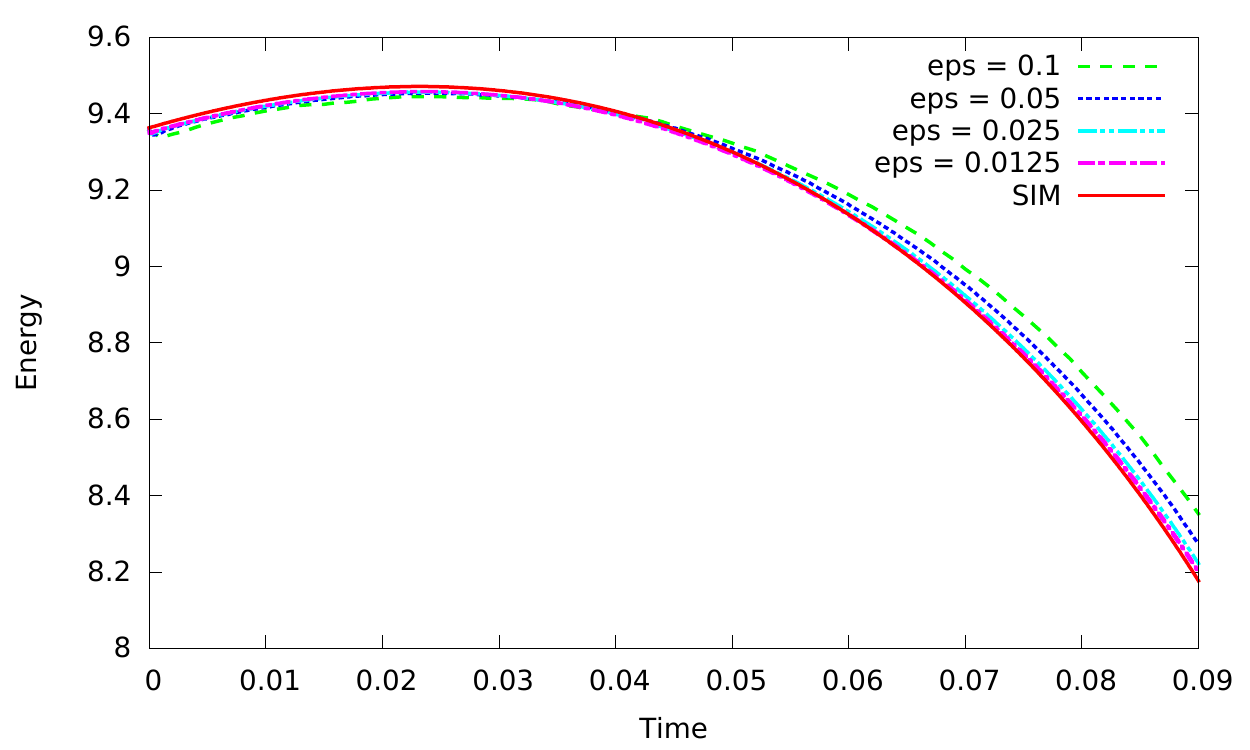}
\end{subfigure}
\begin{subfigure}{0.48\textwidth}
\includegraphics[width = \textwidth]{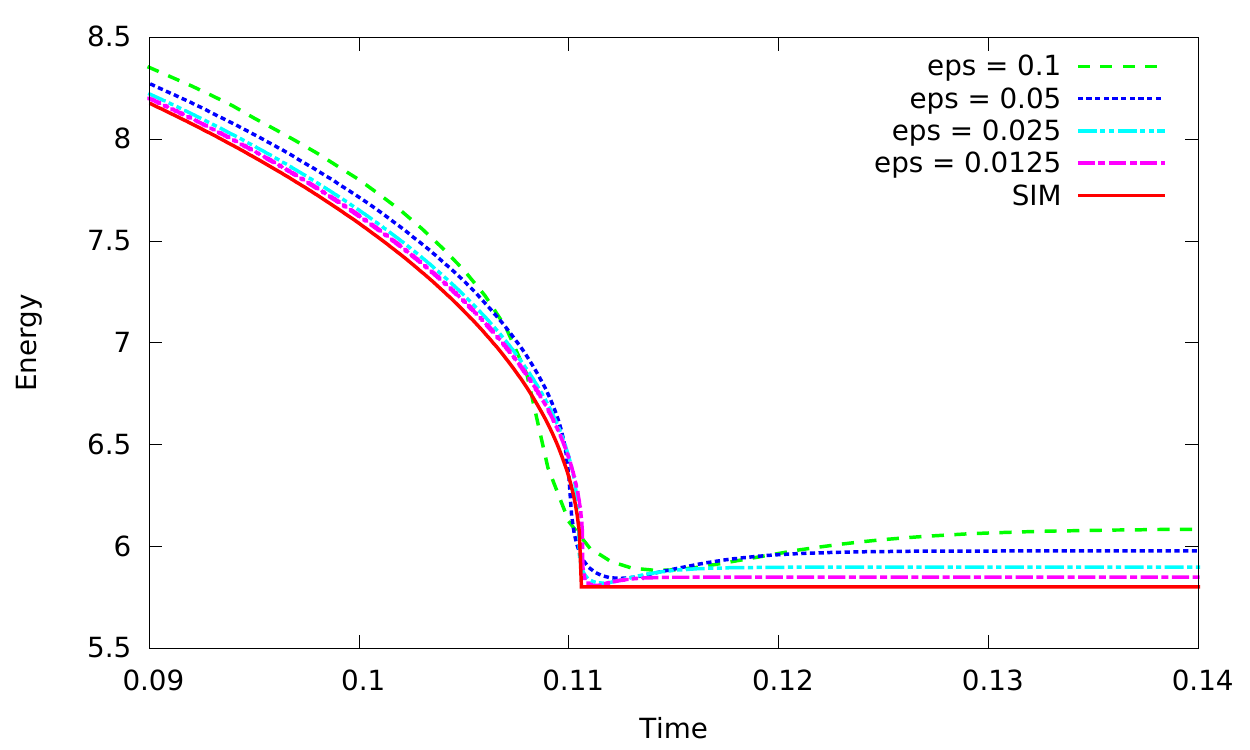}
\end{subfigure}
\caption{Energy plots for the example in Section \ref{sec:MovingSphere} with $\bar{v} = 10$, $\bar{M} = 5$. We compare the Ginzburg-Landau energy, \eqref{GLEF}, with the sharp interface energy, \eqref{SIM:Energy}.}
\label{sharp}
\end{figure}

In Figure \ref{sharp} we display the evolution of the energies \eqref{GLEF} for several values of $\eps$ as well as the limiting energy \eqref{SIM:Energy}. Around the time $0.11$ the solution to the sharp interface model becomes singular as then $\theta_2(t) \to \pi$. The asymptotic analysis is not valid around such events but we see that even then the approximation gets more accurate as $\eps \to 0$.

\subsection{Scaling Effects}
\label{sec:mob}

In our analysis we saw that different scalings of $\bar{M}$ lead to different limiting free boundary problems, namely \eqref{eq:FBP_bulk}, \eqref{eq:FBP_int} for $\bar{M} \sim \eps^0$ and \eqref{eq:HOP_bulk}, \eqref{eq:HOP_int} for $\bar{M} \sim \eps^{-1}$. In this example we present a pair of simulations to demonstrate the differing behaviour of solutions to the ESCH equation in dependence of the scaling of $\bar{M}$ in $\eps$.

\begin{figure}
\begin{subfigure}{0.32\textwidth}\centering
\includegraphics[width = 1\textwidth]{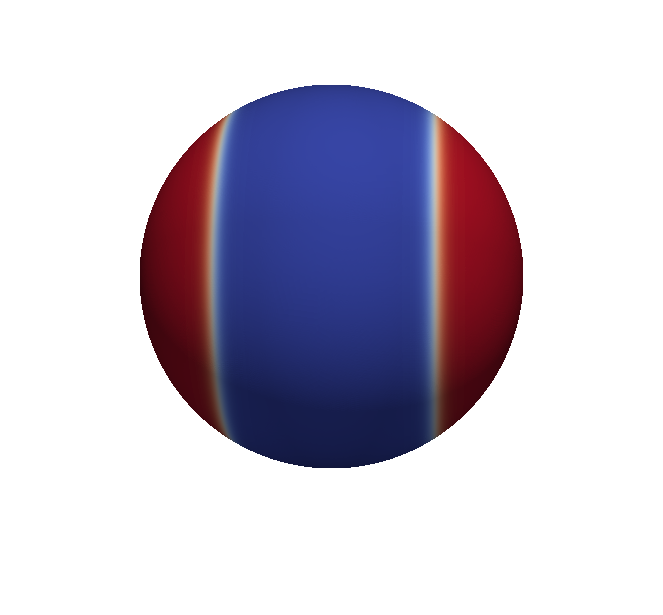}
\caption{$t=0$  }
\label{mob_initial}
\end{subfigure}
\begin{subfigure}{0.32\textwidth}\centering
\includegraphics[width =1 \textwidth]{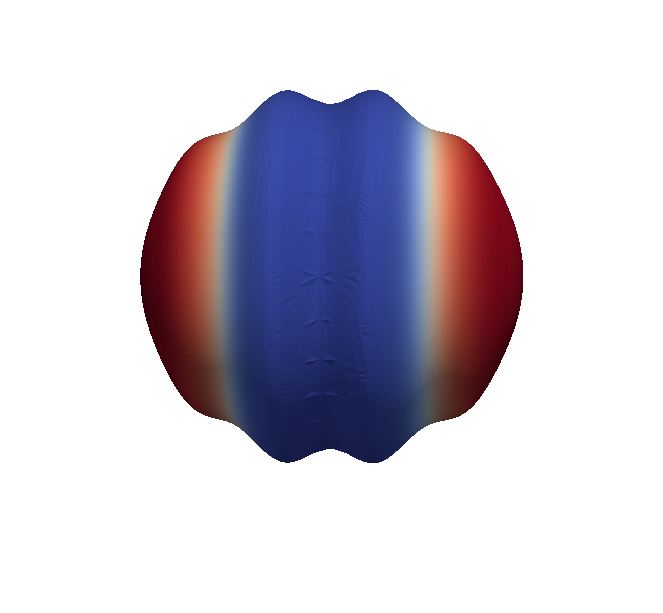}
\caption{$t=0.2$, $\bar{M} = 5$}
\label{mob_slow}
\end{subfigure}
\begin{subfigure}{0.32\textwidth}\centering
\includegraphics[width = 1\textwidth]{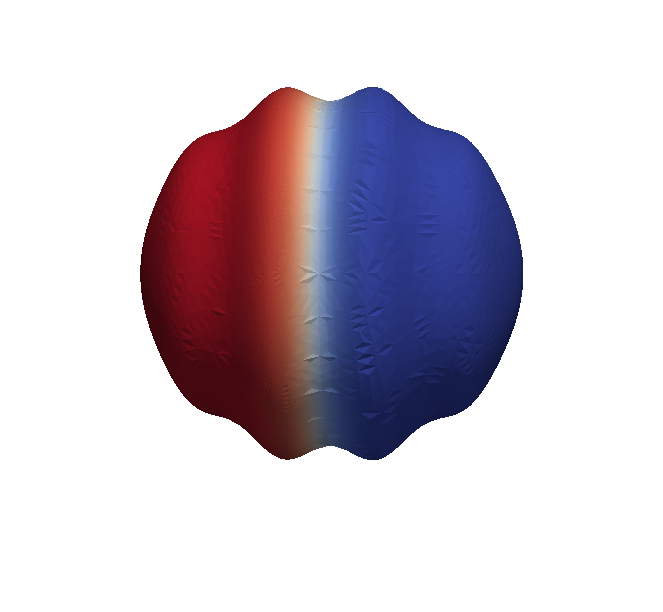}
\caption{$t=0.2$, $\bar{M} = 50$}
\label{mob_fast}
\end{subfigure}
\caption{Initial configuration (left) and different states (middle, right) achieved by varying the mobility, $\eps = 0.1$. See Section \ref{sec:mob} for other parameters and details.}
\label{fig:mob}
\end{figure}

We begin with the unit sphere with two regions of phase $u_b$ at opposite sides of the sphere separated by a band of phase $u_a$, as displayed in Figure \ref{mob_initial}. As with the previous example, the two regions are of different size so that we can expect to see the coarsening of the phase $u_b$ if the surface velocity is zero (rotating Figure \ref{coarse} through $90$ degrees would produce this solution). We choose a surface velocity to deform the sphere so as to introduce obstacles by increasing the radius of $(y,z)$-circles. More specifically, the surface $\mmm{G}_T$ is given as the image of $Q:\bbb{S}^2 \times [0,0.2] \to \bbb{R}^3$ by 
\[
Q(x,y,z,t) = (1-\tilde{t}) (x,y,z) + \tilde{t} (x,\rho(x)y,\rho(x)z), \quad \tilde{t}=\min(0.05,t)
\]
where $\rho(x) = 1-\frac{1}{2} \cos^2(2\pi x)$. For a fixed interfacial thickness parameter $\eps = 0.1$ we use two different values for the mobility, namely $\bar{M} = 5$ and, dividing by $\eps$, $\bar{M} = 50$. The other parameters are in Table \ref{tbl:mob_topo} on the left. 

Based on the observations in the previous example, the slightly larger domain of phase $u_b$ is expected to attempt to grow at the expense of the smaller domain until the latter vanishes. This phenomenon is driven by the different values of the geodesic curvature of the phase interfaces. By altering the radii of $(y,z)$-circles over time as given above the curvature of the underlying surface is varied. If a phase interface moves into the affected area, then its geodesic curvature is changed in such a way that further movement towards the equator is damped. In the case $\bar{M} = 5$ we observe (see Figure \ref{mob_slow}) that coarsening indeed is prevented and two domains of phase $b$ persist. In turn, by scaling the mobility with $\tfrac{1}{\eps}$ we increase the Cahn-Hilliard dynamics and, thus, the velocity of the phase interface. Indeed, $\bar{M} = 50$ is big enough such that the system can coarsen before the deformation can impact on the dynamics (see Figure \ref{mob_fast}).

\begin{table}
 \tabcolsep=3pt
 \begin{tabular}{lll}
  \hline	
  Parameter & Data for Figure \ref{fig:mob} & Data for Figure \ref{torus} \\
  \hline
  $u_a$, $u_b$; $\bar{M}; T$ & -1,  1;  5 or 50;  0.2 & -1,  1;  10;  0.5\\
  $u_0(x)$ & $\left\{\begin{array}{cc} \tanh\left(\frac{1-\arccos(x)}{\eps \sqrt{2}}\right) & \arccos(x)<1.55  \\
\tanh\left(\frac{\arccos(x)-2.1}{\eps \sqrt{2}}\right) & \arccos(x)\geq 1.55 \end{array}\right.$ & $\tanh\left(\frac{0.7-x_1}{\eps\sqrt{2}}\right)$ \\
  \hline
 \end{tabular}
 \caption{Simulation data for Sections \ref{sec:mob} and \ref{sec:topo}.} \label{tbl:mob_topo} 
\end{table}

\subsection{Topological Changes}
\label{sec:topo}

Topological changes can be particularly difficult to simulate for free boundary problems. In this example we display a topological change of the interface induced by the motion of the surface that would not happen in a stationary setting. 

The surface is a torus which we denote by $\mathbb{T}(R,r)$ where $R$ is the major radius and $r$ is the minor radius. We deform the torus by making $R$ and $r$ time dependent functions, specifically, $R(t) = \sqrt{2} + 1.2 \sin(2 \pi t)$ and $r(t) = 1 - 0.65 \sin(2\pi t)$, thus increasing the overall surface area in the interval $0 < t < 0.25$, decreasing the surface area in the interval $0.25 < t < 0.5$, and obtaining the same surface at final time $T=0.5$ as at $t=0$.

We consider an initial phase distribution which contains a single connected interface using the profile function as described in Table \ref{tbl:mob_topo} on the right. Note that this function is only dependent on the spatial co-ordinate $x_1$, rather than any tangential co-ordinate. This creates a relatively large initial energy, however the interfacial layers quickly relax to energetically more favourable profiles. Thus when reporting the energy of the system we start shortly after initialisation. 

\begin{figure}
\begin{subfigure}{0.24\textwidth}\centering
\includegraphics[width = \textwidth,trim = 150 500 150 30, clip]{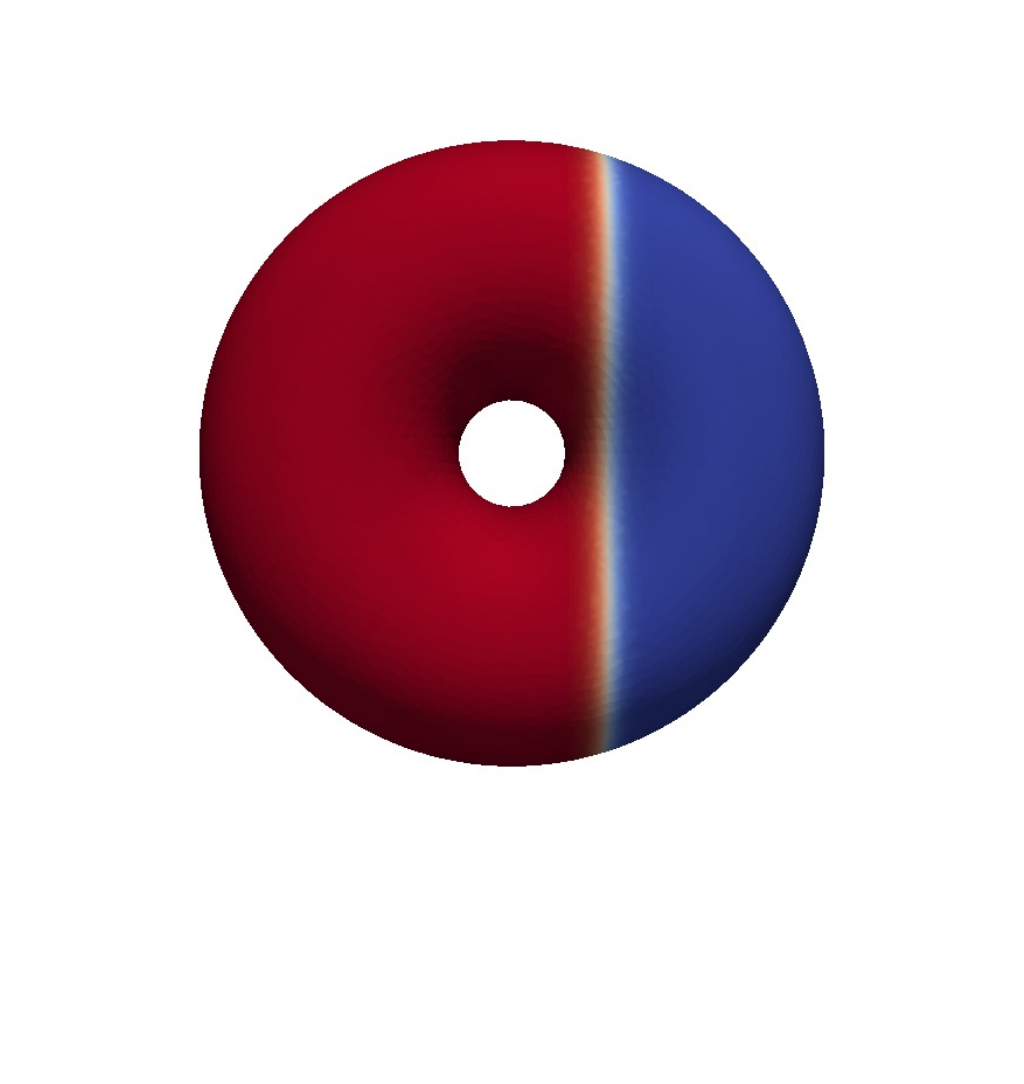}
\caption{$t=0.0$ } 
\end{subfigure}
\begin{subfigure}{0.24\textwidth}\centering
\includegraphics[width = \textwidth,trim = 150 500 150 30, clip]{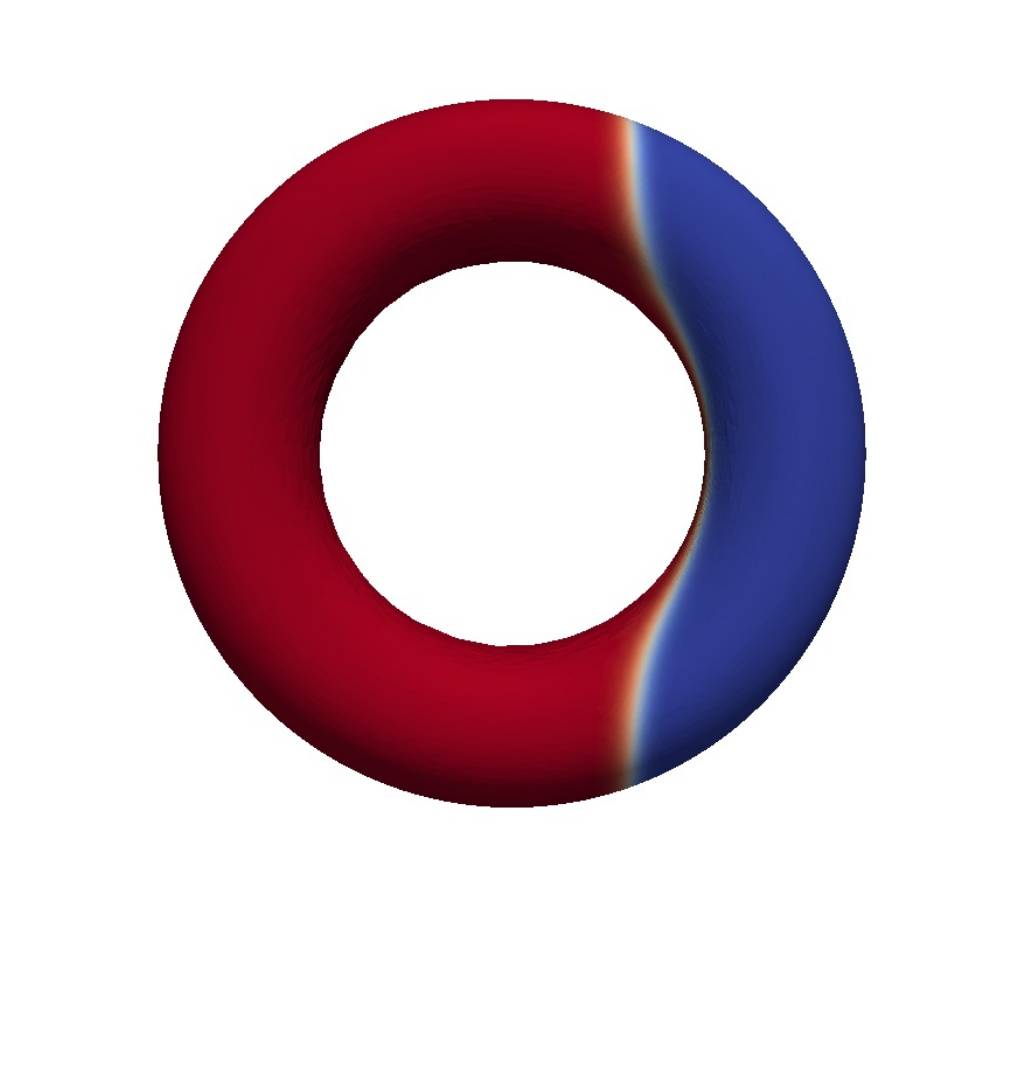}
\caption{$t=0.1$  }
\end{subfigure}
\begin{subfigure}{0.24\textwidth}\centering
\includegraphics[width = \textwidth,trim = 150 500 150 30, clip]{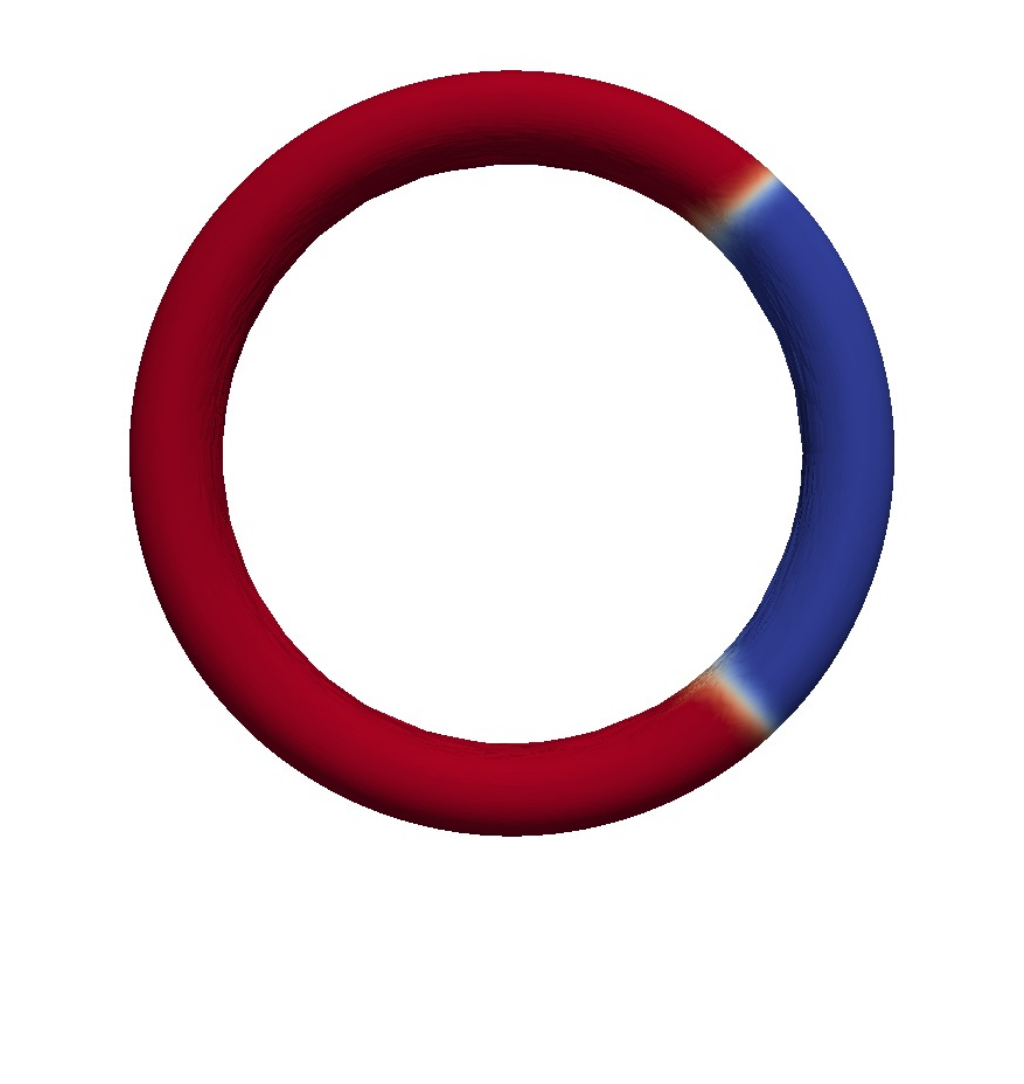}
\caption{$t=0.25$  }
\end{subfigure}
\begin{subfigure}{0.24\textwidth}\centering
\includegraphics[width = \textwidth,trim = 150 500 150 30, clip]{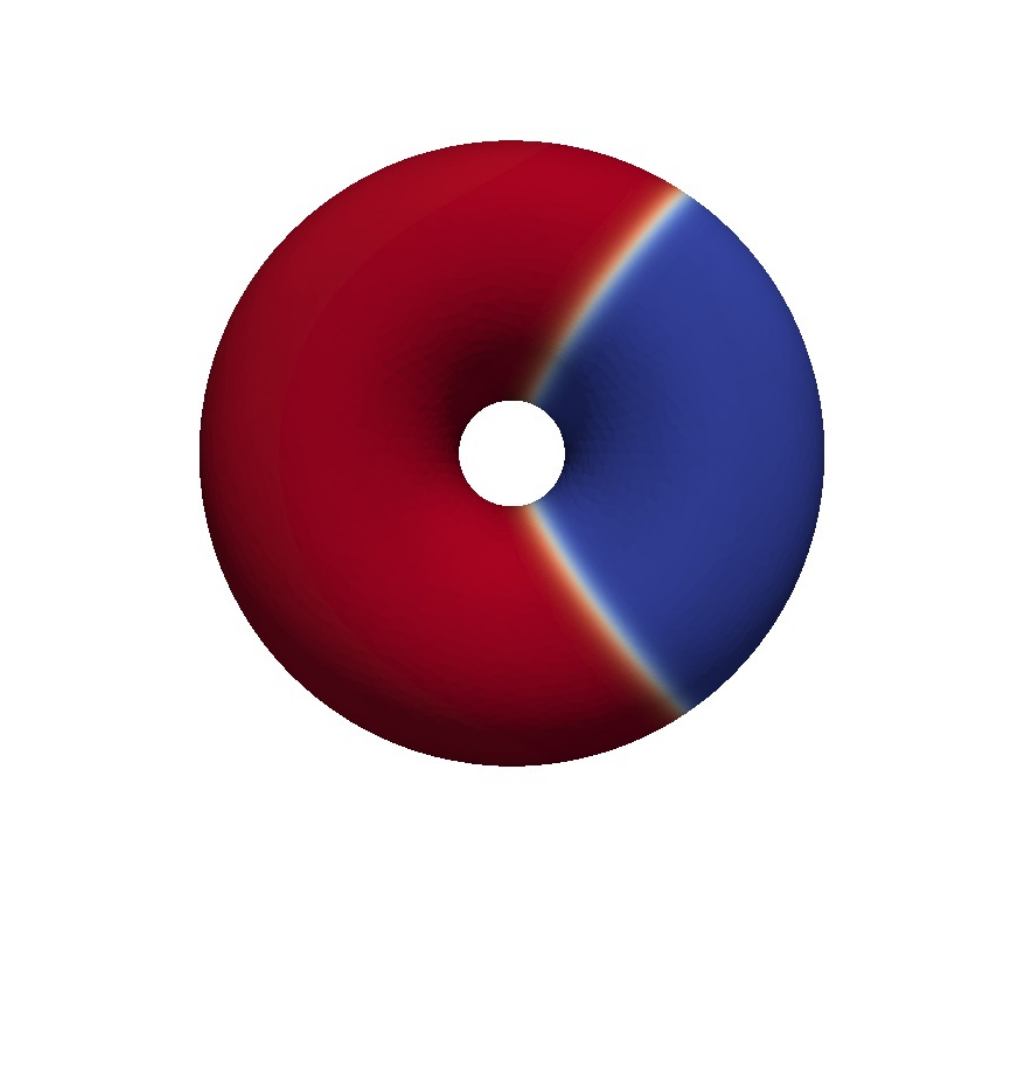}
\caption{$t=0.5$ }
\end{subfigure}
\caption{Topological change of the interface as discussed in Section \ref{sec:topo}. $\eps = 0.71$, $\bar{M}=10$.} \label{torus}
\end{figure}

\begin{figure}
\begin{subfigure}{0.24\textwidth}\centering
\includegraphics[width = \textwidth,trim = 150 500 150 30, clip]{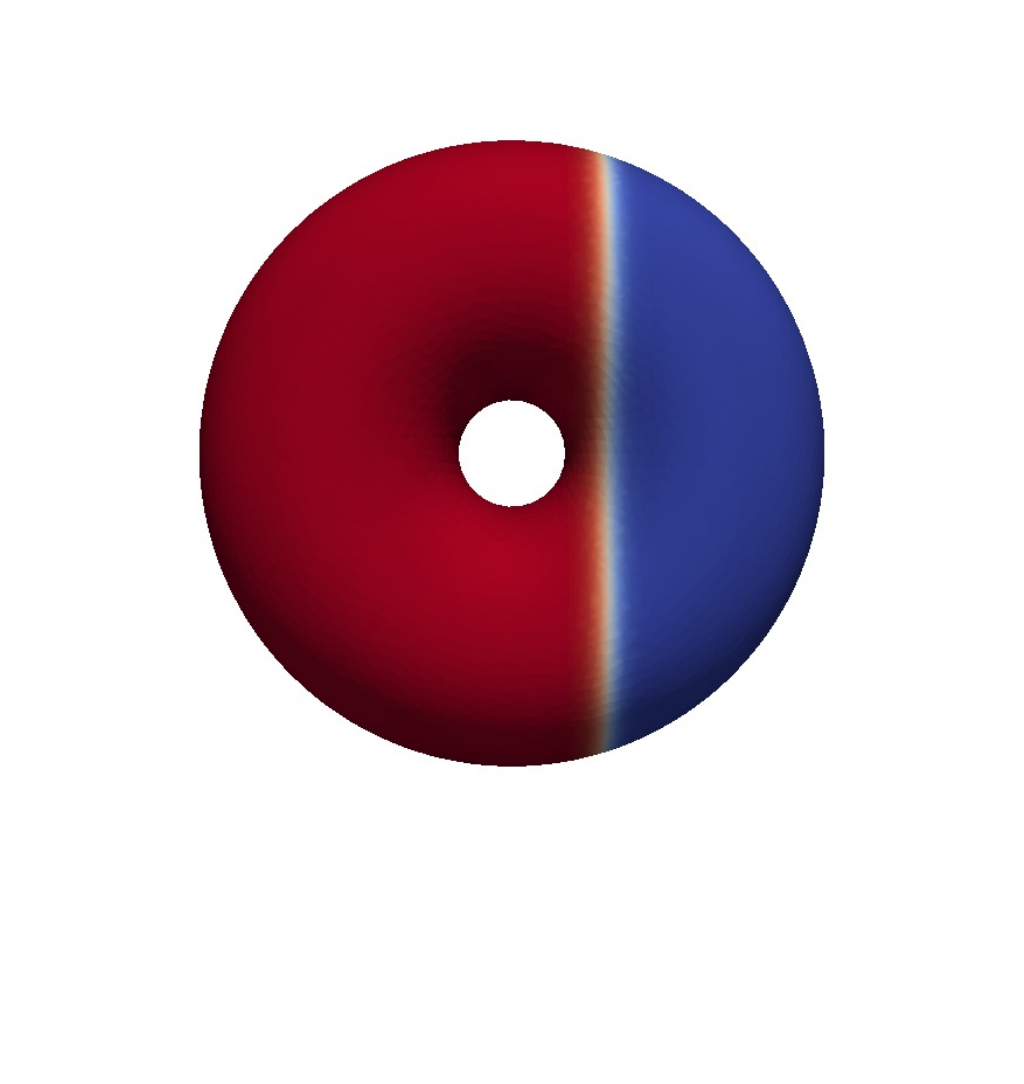}
\caption{$t=0.0$ } 
\end{subfigure}
\begin{subfigure}{0.24\textwidth}\centering
\includegraphics[width = \textwidth,trim = 150 500 150 30, clip]{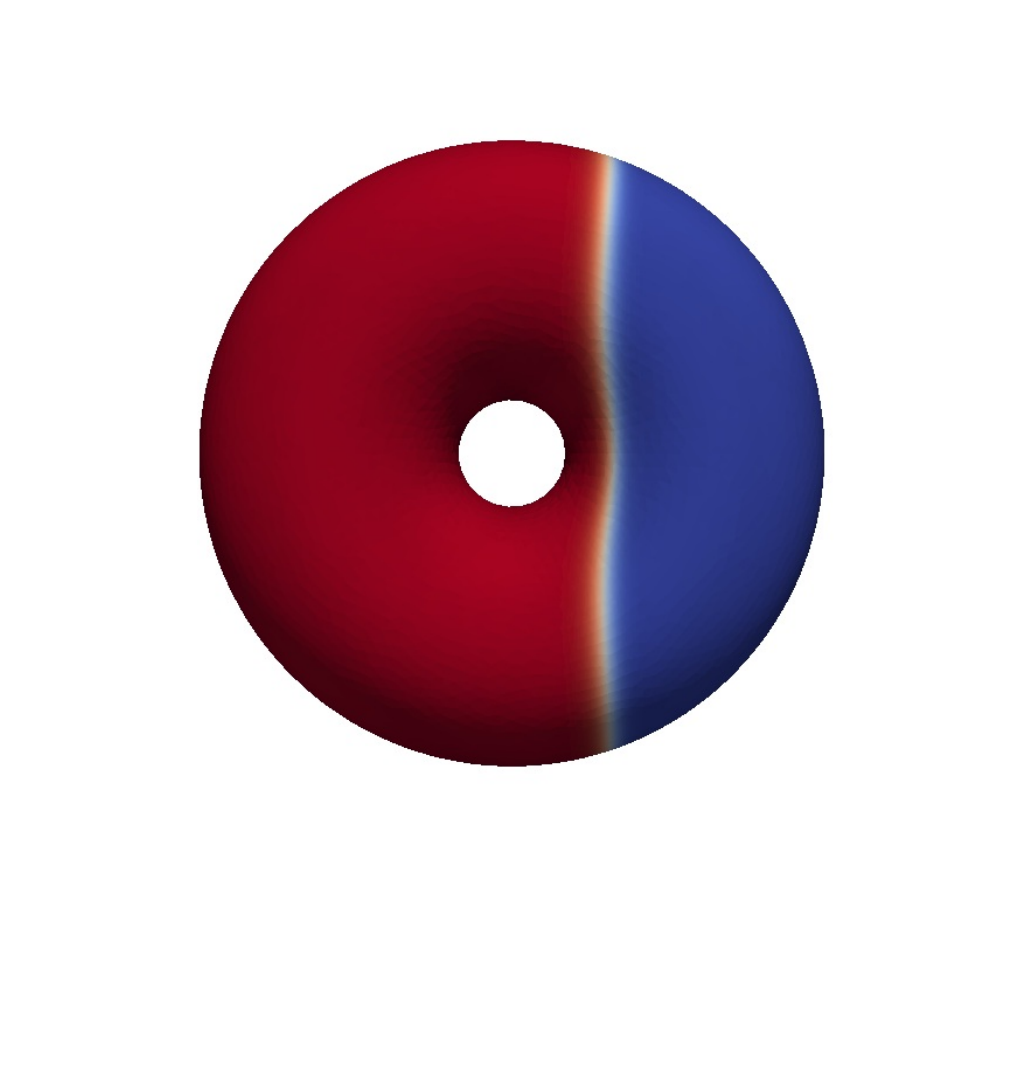}
\caption{$t=0.1$  }
\end{subfigure}
\begin{subfigure}{0.24\textwidth}\centering
\includegraphics[width = \textwidth,trim = 150 500 150 30, clip]{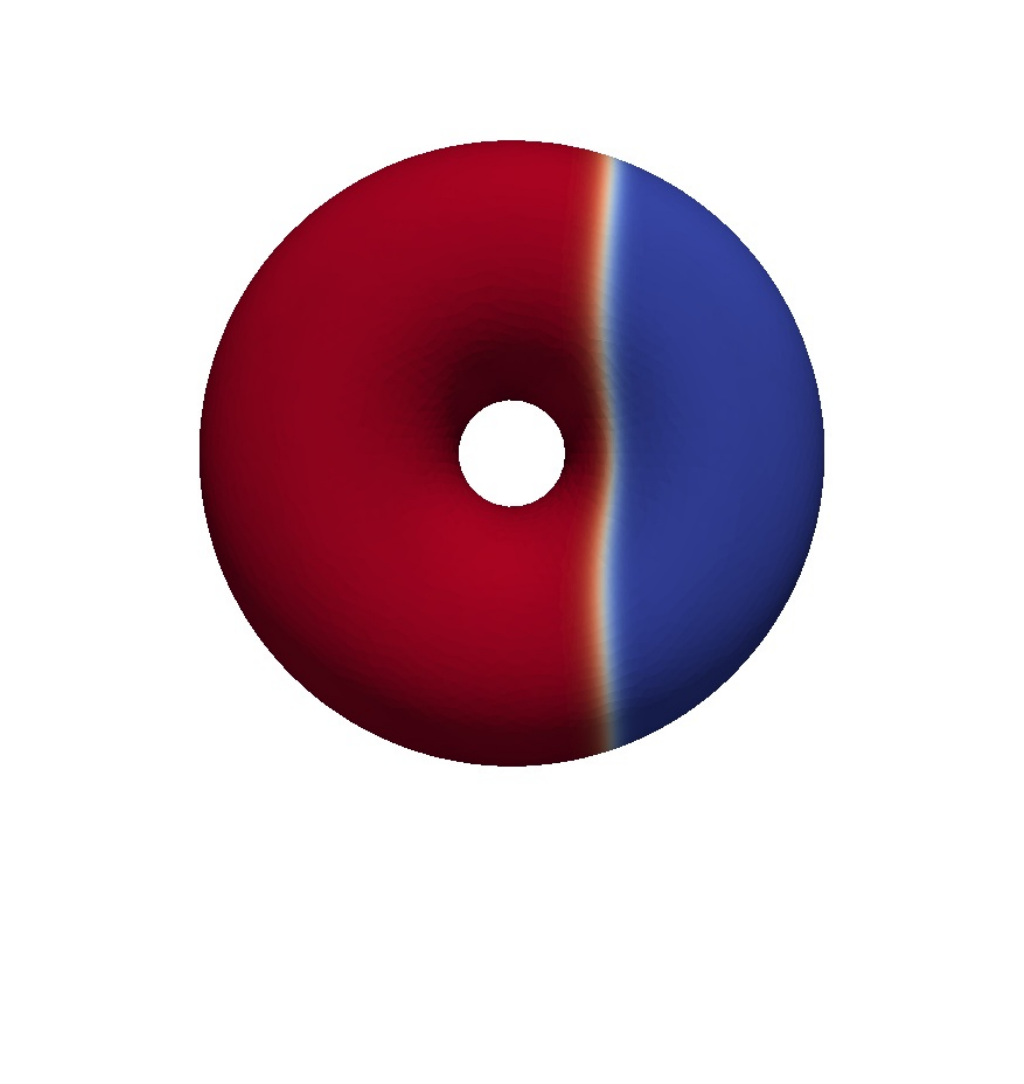}
\caption{$t=0.25$  }
\end{subfigure}
\begin{subfigure}{0.24\textwidth}\centering
\includegraphics[width = \textwidth,trim = 150 500 150 30, clip]{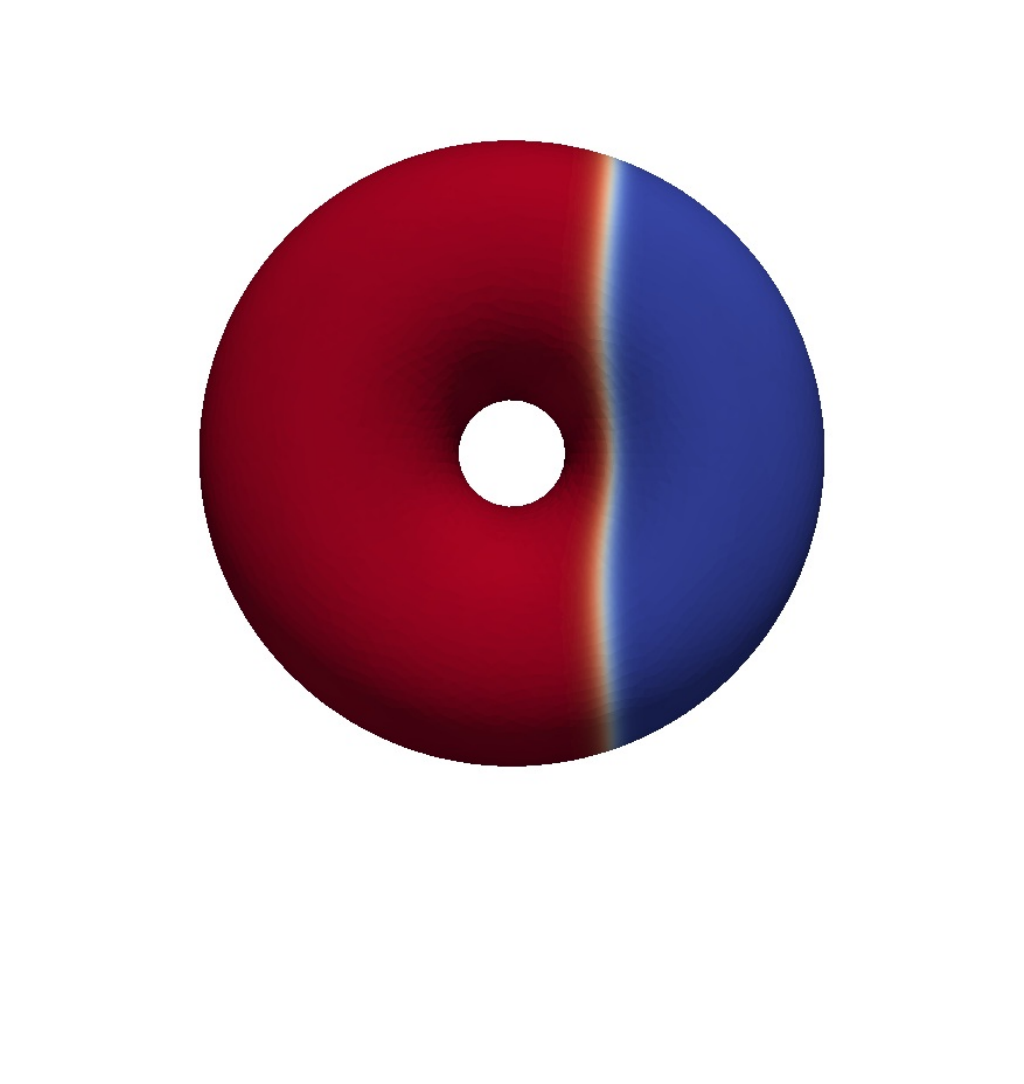}
\caption{$t=0.5$ }
\end{subfigure}
\caption{Relaxation of singular interface in case of stationary surface as discussed in Section \ref{sec:topo}. $\eps = 0.71$, $\bar{M}=10$.} \label{torus_stationary}
\end{figure}

On the stationary torus $\mathbb{T}(R(0),r(0))$ the described phase interface would evolve only so as to reduce its length but without any topological change as seen in Figure \ref{torus_stationary}. However, by changing the ratio of the two radii, the phase interface can be driven to self intersect and even to induce a topological change. In Figure \ref{torus} we display the latter solution at $4$ time steps for one specific value of $\eps=0.71$. We observe that the interfacial layer self intersects and splits up into two independent interfacial layers through the hole of the torus. These remain stable when the surface relaxes back to its original shape. 

\begin{figure}
\includegraphics[width = 0.6\textwidth]{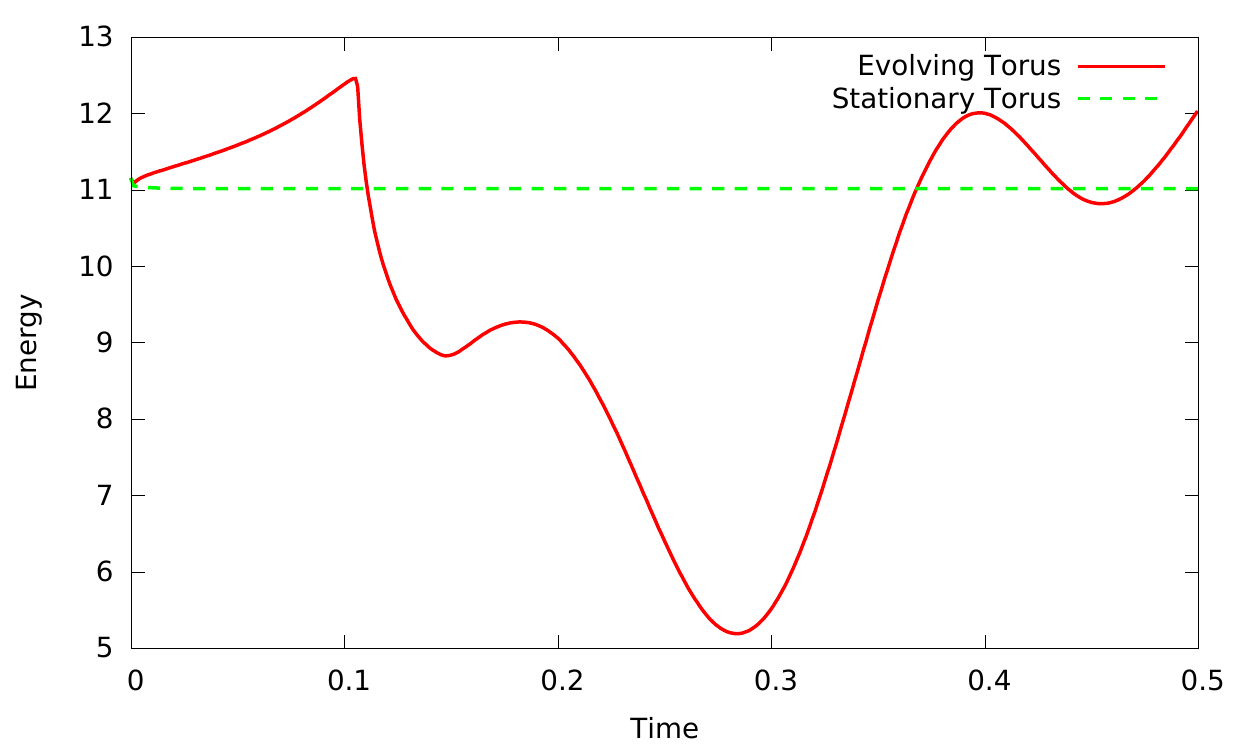}
\caption{Comparison of energy profiles in case of stationary torus and moving torus as discussed in Section \ref{sec:topo} and displayed in Figures \ref{torus} and \ref{torus_stationary}.} \label{torusgraph}
\end{figure}

In Figure \ref{torusgraph} we also include a plot showing the energy evolution of solutions for the two discussed cases. For the stationary surface we see a small drop in the energy due to relaxation and then it remains constant. In contrast the energy in the evolving setting increases initially before the rapid transition through the topological change, around $t=0.1$, as the forming two interfaces becomes energetically more favourable. When the surface returns to its original proportions the total energy is higher than that of the final resting energy in the stationary setting, indicating a local minimum.

\section{Conclusion}

The Cahn-Hilliard equation has been derived on an evolving surface and studied with respect to the asymptotic limit as the interfacial width parameter $\eps$ tends to zero. We have used the method of formally matched asymptotic expansions to derive two classes of free boundary problems. For this purpose, techniques from the literature had to be extended in order to be able to deal with the domain movement, most notably, with regards to expanding the material time derivative in the interfacial coordinates \eqref{tcov}. The two limits equate to two different scalings of the mobility parameter $\bar{M}$ in $\eps$ and consist of the equations \eqref{eq:FBP_bulk}, \eqref{eq:FBP_int} for the scaling $\bar{M} \sim \eps^0$ and \eqref{eq:HOP_bulk}, \eqref{eq:HOP_int} for $\bar{M} \sim \eps^{-1}$. 

As long as there is no degeneracy the result is rather independent of the type of double-well potential $F$ and mobility $M$. As in the case of a stationary domain the two minima of $F$ determine the phase field value in the outer region but the chemical potential is no longer harmonic, in general, but depends on local area changes of the domain due to its movement \eqref{eq:FBP_bulk}, \eqref{eq:HOP_bulk}. Otherwise, the domain evolution only manifests as an additional transport term for the phase interface movement (last equations of \eqref{eq:FBP_int} and \eqref{eq:HOP_int}, respectively). 

We also get a sensible result in the deep quench limit of the logarithmic potential to the double obstacle potential but only if the velocity field is surface divergence free. As on a flat, stationary domain in the case of a fast mobility we obtain surface diffusion but in terms of the geodesic curvature and with an additional transport term which is due to the domain movement \eqref{eq:surfdiff}. 

Numerical simulations both on curves in 2D and on surfaces in 3D support the theoretical findings. There is evidence that the asymptotic analysis indeed identifies the correct limiting free boundary problem, most notably, in the rotationally symmetric setting in Section \ref{sec:MovingSphere}. Techniques such as used in \cite{Stoth96} may translate and enable a rigorous proof of a convergence result. We have also seen how the domain movement impacts on the evolution of the phase interfaces and can lead to a different behaviour than in the case of a stationary domain, for instance, with respect to the topology of the phase interface in Section \ref{sec:topo}.

\subsection*{Acknowledgement}

This research has been supported by the British Engineering and Physical Sciences Research Council (EPSRC), Grant EP/H023364/1.

\end{document}